\documentclass[11pt]{article}
\usepackage{amssymb}
\usepackage{amsmath}
\usepackage{amsthm}
\newcommand{\cA}{{\cal A}}
\newcommand{\cB}{{\cal B}}
\newcommand{\cC}{{\cal C}}
\newcommand{\cD}{{\cal D}}
\newcommand{\cH}{{\cal H}}
\newcommand{\cE}{{\cal E}}

\newcommand{\cJ}{{\cal J}}
\newcommand{\cO}{{\cal O}}
\newcommand{\cL}{{\cal L}}
\newcommand{\cM}{{\cal M}}
\newcommand{\cN}{{\cal N}}
\newcommand{\cF}{{\cal F}}
\newcommand{\cK}{{\cal K}}
\newcommand{\cP}{{\cal P}}
\newcommand{\cQ}{{\cal Q}}
\newcommand{\cS}{{\cal S}}
\newcommand{\cT}{{\cal T}}
\newcommand{\cU}{{\cal U}}
\newcommand{\cV}{{\cal V}}
\newcommand{\cW}{{\cal W}}
\newcommand{\cX}{{\cal X}}
\newcommand{\cY}{{\cal Y}}
\newcommand{\cZ}{{\cal Z}}
\newcommand{\cSN}{{\cal SN}}

\renewcommand{\AA}{{\mathbb A}}
\newcommand{\BB}{{\mathbb B}}
\newcommand{\GG}{{\mathbb G}}
\newcommand{\NN}{{\mathbb N}}
\newcommand{\ZZ}{{\mathbb Z}}

\newcommand{\OO}{{\mathbb O}}

\newcommand{\HH}{{\mathbb H}}
\newcommand{\TT}{{\mathbb T}}

\renewcommand{\gg}{\mathfrak{g}}  

\newcommand{\gp}{\mathfrak{p}}
\newcommand{\gq}{\mathfrak{q}}

\newcommand{\gr}{\mathfrak{r}}
\newcommand{\gs}{\mathfrak{s}}

\newcommand{\on}{\operatorname}
\newcommand{\Sh}{\on{Sh}}

\newcommand{\RCov}{\on{RCov}}
\newcommand{\mult}{\on{mult}}

\newcommand{\Rep}{{\on{Rep}}}

\newcommand{\Laum}{{\on{Laum}}}
\newcommand{\Qlb}{\mathbb{\bar Q}_\ell}
\newcommand{\Gm}{\mathbb{G}_m}

\newcommand{\A}{\mathbb{A}}

\newcommand{\toup}[1]{\stackrel{#1}{\to}}
\newcommand{\hook}[1]{\stackrel{#1}{\hookrightarrow}}

\newcommand{\getsup}[1]{\stackrel{#1}{\gets}}
\newcommand{\Sp}{\on{\mathbb{S}p}}
\newcommand{\Spin}{\on{\mathbb{S}pin}}
\newcommand{\GSpin}{\on{G\mathbb{S}pin}}
\newcommand{\GSp}{\on{G\mathbb{S}p}}
\newcommand{\IC}{\on{IC}}
\newcommand{\CT}{\on{CT}}
\newcommand{\Hom}{\on{Hom}}
\newcommand{\Mod}{\on{Mod}}

\newcommand{\Sym}{\on{Sym}}
\newcommand{\SO}{\on{S\mathbb{O}}}
\newcommand{\GO}{\on{G\mathbb{O}}}
\newcommand{\Ker}{\on{Ker}}

\newcommand{\Aut}{\on{Aut}}

\newcommand{\RG}{\on{R\Gamma}}

\newcommand{\sym}{\on{sym}}
\newcommand{\triv}{\on{triv}}

\newcommand{\Pic}{\on{Pic}}
\newcommand{\uPic}{\on{\underline{Pic}}}
\newcommand{\Bun}{\on{Bun}}

\newcommand{\Bunb}{\on{\overline{Bun}} }
\newcommand{\Bunt}{\on{\widetilde\Bun}}

\newcommand{\rk}{\on{rk}}
\newcommand{\Spec}{\on{Spec}}

\newcommand{\supp}{\on{supp}}

\newcommand{\HOM}{{{\cal H}om}}
\newcommand{\END}{{{\cal E}nd}}
\newcommand{\Gr}{\on{Gr}}
\newcommand{\Grb}{\overline{\Gr}}

\newcommand{\GL}{\on{GL}}

\newcommand{\pr}{\on{pr}}
\newcommand{\id}{\on{id}}

\newcommand{\QED}{$\square$} 
\newcommand{\Fq}{\mathbb{F}_q}  
\newcommand{\Fp}{\mathbb{F}_p}  
\newcommand{\Out}{\on{\mathbb{O}ut}}
\newcommand{\iso}{{\widetilde\to}}

\newcommand{\comp}{\circ}
\newcommand{\Four}{\on{Four}}
\renewcommand{\H}{{\on{H}}}   

\newcommand{\DD}{\mathbb{D}}  
\newcommand{\D}{\on{D}}       

\newcommand{\ov}[1]{\overline{#1}}
\newcommand{\select}[1]{{\it{#1}}}

\newcommand{\und}[1]{\underline{#1}}

\renewcommand{\div}{\on{div}}

\renewcommand{\P}{{\on{P}}}

\newcommand{\<}{\langle}
\renewcommand{\>}{\rangle}

\newcommand{\Av}{\on{Av}}
\newcommand{\ev}{\mathit{ev}}

\newcommand{\Loc}{\on{Loc}}
\newcommand{\Lie}{\on{Lie}}
\newcommand{\Sph}{\on{Sph}}

\newcommand{\ttimes}{\tilde\times}

\newcommand{\act}{\on{act}}
\newcommand{\dimrel}{\on{dim.rel}}

\newcommand{\Isom}{\on{Isom}}

\renewcommand{\Im}{\on{Im}}

\newcommand{\sign}{\on{sign}}
\newcommand{\SL}{\on{SL}}
\newcommand{\St}{\on{St}}

\newcommand{\WP}{\on{WP}} 
\newcommand{\GWP}{\on{GWP}}
\newcommand{\BP}{\on{BP}}
\newcommand{\DR}{\on{DR}}
\newcommand{\LocSys}{\on{LocSys}}
\newcommand{\jac}{\on{jac}}

\newtheorem{Lm}{Lemma}
\newtheorem{Th}{Theorem}
\newtheorem{Pp}{Proposition}
\newtheorem{Cor}{Corollary}
\newtheorem{Con}{Conjecture}
\newtheorem{Slm}{Sublemma}

\theoremstyle{remark}
\newtheorem{Rem}{Remark}

\theoremstyle{definition}
\newtheorem{Def}{Definition}

\newenvironment{Prf}{\par\noindent {\it Proof }}{\QED}
\newcommand{\Step}[1]{\par\noindent{\bf Step {#1}}.}

\begin{document}

\author{Sergey Lysenko}
\title{Geometric Waldspurger periods}
\date{}
\maketitle
\begin{abstract}
\noindent{\scshape Abstract}\hskip 0.8 em 
Let $X$ be a smooth projective curve. We consider the dual reductive pair $H=\GO_{2m}$, $G=\GSp_{2n}$ over $X$, where $H$ splits on an \'etale two-sheeted covering $\pi:\tilde X\to X$. Let $\Bun_G$ (resp., $\Bun_H$) be the stack of $G$-torsors (resp., $H$-torsors) on $X$. We study the functors $F_G$ and $F_H$ between the derived categories $\D(\Bun_G)$ and $\D(\Bun_H)$, which are analogs of the classical theta-lifting operators in the framework of the geometric Langlands program.

 Assume $n=m=1$ and $H$ nonsplit, that is, $H=\pi_*\Gm$ with $\tilde X$ connected. 
We establish the geometric Langlands functoriality for this pair. Namely, we show that $F_G:\D(\Bun_H)\to\D(\Bun_G)$ commutes with Hecke operators with respect to the corresponding map of Langlands $L$-groups $^L H\to {^L G}$. 

As an application, we calculate Waldspurger periods of cuspidal automorphic sheaves on $\Bun_{\GL_2}$ and Bessel periods of theta-lifts from $\Bun_{\GO_4}$ to $\Bun_{\GSp_4}$. Based on these calculations, we give three conjectural constructions of certain automorphic sheaves on $\Bun_{\GSp_4}$ (one of them makes sense for $\cD$-modules only). 
\end{abstract} 

\bigskip
\centerline{\scshape 1. Introduction and main results}

\bigskip\noindent
1.0 This paper, which is a sequel to \cite{Ly}, is a part of two (related) research projects: i) geometric  version of the Howe correspondence (an analogue of the theta-lifting in the framework of the geometric Langlands program); ii) geometric Langlands program for $\GSp_4$. 

 We consider only the (unramified) dual reductive pair $(H=\GO_{2m}, G=\GSp_{2n})$ over a smooth projective connected curve $X$. Let $\Bun_G$ (resp., $\Bun_H$) denote the stack of $G$-torsors (resp., $H$-torsors) on $X$. Using the theta-sheaf introduced in \cite{Ly}, we define functors $F_G:\D(\Bun_H)\to\D(\Bun_G)$ and $F_H:\D(\Bun_G)\to\D(\Bun_H)$ between the corresponding derived categories, which are geometric analogs of the theta-lifting operators. Based on classical Howe correspondence (cf., for example, \cite{A},\cite{Ku},\cite{MVW},\cite{R}) and our results from \cite{Ly4}, we conjecture a precise relation between the theta-lifting functors and Hecke functors on $\Bun_G$ and $\Bun_H$ (cf. Conjecture~\ref{Con_functoriality}). For $n=m$ (resp., for $m=n+1$) the functor $F_G$ (resp., $F_H$) is expected to realize the geometric Langlands functoriality for a morphism of Langlands $L$-groups $H^L\to G^L$ (resp., $G^L\to H^L$). 
 
  We prove this conjecture for the dual pair $(\GO_2,\GL_2)$, where $\GO_2=\pi_*\Gm$ is a group scheme over $X$, here $\pi:\tilde X\to X$ is a nontrivial \'etale two-sheeted covering. If $\tilde E$ is a rank one local system on $\tilde X$ then this provides a new proof of the geometric Langlands conjecture for $\pi_*\tilde E$ independent of the existing proof due to Frenkel, Gaitsgory and Vilonen (\cite{FGV}, \cite{G}). 
   
 Let us describe our remaining main results in a form less technical then their actual formulation. Assume that the ground field $k=\Fq$ is finite of $q$ elements (with $q$ odd). Set $G=\GL_2$. Let $E$ be a rank 2 irreducible $\ell$-adic local system on $X$. Write $\Aut_E$ for the corresponding automorphic sheaf on $\Bun_G$ (cf. Definition~\ref{Def_Aut_E}). Let $f_E: \Bun_G(k)\to\Qlb$ denote the function `trace of Frobenius' of $\Aut_E$. 
 
 Let $\phi: Y\to X$ be a nontrivial \'etale two-sheeted covering. Write $\Pic Y$ for the Picard stack of $Y$. 
Let $\cJ$ be a rank one local system on $Y$ equipped with an isomorphism $N(\cJ)\,\iso\,\det E$, where $N(\cJ)$ is the norm of $\cJ$ (cf. A.1). Write $f_{\cJ}: (\Pic Y)(k)\to \Qlb$ for the corresponding character (the trace of Frobenius of the automorphic local system $A\cJ$ corresponding to $\cJ$). 
The Waldspurger period of $f_E$ is
$$
\int_{\cB\in (\Pic Y)(k)/(\Pic X)(k)} \, f_E(\phi_*\cB) f_{\cJ}^{-1}(\cB)d\cB
$$
(the function that we integrate does not change when $\cB$ is tensored by $\phi^*L$, $L\in\Pic X$), here $d\cB$ is a Haar measure. A beautiful theorem of Waldspurger says that \select{the square} of this period is equal (up to an explicit harmless coefficient) to the value of the $L$-function $L(\phi^*E\otimes \cJ^{-1}, \frac{1}{2})$ (cf. \cite{Wa}). 
  
  We prove a geometric version of this result. The role of the $L$-function in geometric setting is played by the complex 
\begin{equation}
\label{complex_L_function}
\oplus_{d\ge 0} \RG(Y^{(d)}, (\phi^*E\otimes \cJ^*)^{(d)})[d]
\end{equation}
Here $Y^{(d)}$ is the $d$-th symmetric power of $Y$, and $V^{(d)}$ denotes the $d$-th symmetric power of a local system $V$ on $X$.  The geometric Waldspurger period is
\begin{equation}
\label{complex_WP_intro}
\RG_c(\Pic Y/\Pic X, \phi_1^*\Aut_E\otimes A\cJ^{-1}),
\end{equation}
where $\phi_1: \Pic Y\to\Bun_G$ sends $\cB$ to $\phi_*\cB$. The sense of the quotient $\Pic Y/\Pic X$ is precised in 6.3.3, this stack has two connected components (the degree of $\cB$ modulo two), so (\ref{complex_WP_intro}) is naturally $\ZZ/2\ZZ$-graded. Our Theorem~\ref{Th_WP} says that there is a $\ZZ/2\ZZ$-graded isomorphism between (\ref{complex_L_function}) and the tensor square of (\ref{complex_WP_intro}), the $\ZZ/2\ZZ$-grading of (\ref{complex_L_function}) is given by the parity of $d$. If $\phi^*E$ is still irreducible then (\ref{complex_L_function}) is the exteriour algebra of the vector space $\H^1(Y, \cJ^*\otimes\phi^* E)$, which is placed in cohomological degree zero. 
  
  In the classical theory of automorphic forms there has been a philosophy that for multiplicity one models of representations the corresponding periods of Hecke eigenforms can be expressed in terms of the L-functions (of the corresponding eigenvalue local system). In addition to the Waldspurger periods, we also consider Bessel periods for $\GSp_4$ (Sect.~6) and generalized Waldspurger periods for $\GL_4$ (Sect.~7), which all illustrate this phenomenon. 
  
   Consider now the dual pair $(G,\tilde H)$, where $G=\GSp_4$ and $\tilde H$ is as follows. Let $\GO_4^0$ be given by the exact sequence $1\to\Gm\to\GL_2\times\GL_2\to\GO_4^0\to 1$, where the first map sends $x\in\Gm$ to $(x,x^{-1})$. Let $\pi: \tilde X\to X$ be an \'etale degree two covering, set $\Sigma=\Aut_X(\tilde X)$. The group $\Sigma$ act on this exact sequence permuting the two copies of $\GL_2$. Let $\tilde H$ be the group scheme on $X$, the twisting of $\GO_4^0$ by the $\Sigma$-torsor $\pi:\tilde X\to X$. The above exact sequence yields a morphism of stacks $\rho: \Bun_{2,\tilde X}\to \Bun_{\tilde H}$, here $\Bun_{2,\tilde X}$ denotes the stack of rank two vector bundles on $\tilde X$. Write $\Bun_r$ for the stack of rank $r$ vector bundles on $X$. 
   
   Let $\tilde E$ be an irreducible rank two local system on $\tilde X$, $\Aut_{\tilde E}$ be the corresponding automorphic sheaf on $\Bun_{2,\tilde X}$ (cf. Definition~\ref{Def_Aut_E}). Assume given a rank one local system $\chi$ on $X$ equipped with an isomorphism $\det\tilde E\,\iso\, \pi^*\chi$. Then $\Aut_{\tilde E}$ descends naturally to a perverse sheaf $K_{\tilde E,\chi,\tilde H}$ on $\Bun_{\tilde H}$.  Now assume $\tilde X$ connected. For the theta-lifting functor $F_G:\D(\Bun_{\tilde H})\to\D(\Bun_G)$ our Theorem~\ref{Th_BP} calculates the Bessel periods of 
$K:=F_G(K_{\tilde E,\chi,\tilde H})$. 
   
   At the level of functions the Bessel periods are defined, for example, in \cite{BFF}. In geometric setting, let $P\subset G$ be the Siegel parabolic, $\nu_P:\Bun_P\to\Bun_G$ be the natural map. The stack $\Bun_P$ of $P$-torsors on $X$ classifies collections: $L\in\Bun_2$, $\cA\in\Bun_1$ and an exact sequence $0\to\Sym^2 L\to ?\to\cA\to 0$ on $X$. Let $\cS_P$ be the stack classifying $L\in\Bun_2$, $\cA\in\Bun_1$ together with a map $\Sym^2 L\to\cA\otimes\Omega$. Here $\Omega$ is the canonical line bundle on $X$. So, $\cS_P$ and $\Bun_P$ are dual (generalized) vector bundles over $\Bun_2\times\Bun_1$, and one has the Fourier transform functor $\Four_{\psi}:\D(\Bun_P)\to\D(\Bun_{\cS_P})$ between the corresponding derived categories of $\ell$-adic sheaves. 
   
  Let $\phi: Y\to X$ be a nontrivial \'etale two-sheeted covering. Let $e: \Pic Y\to\cS_P$ be the map sending $\cB\in\Pic Y$ to the pair $L=\phi_*\cB, \, \cA\otimes\Omega=N(\cB)$ with natural symmetric form $\Sym^2 L\to\cA\otimes\Omega$ (cf. Section~6.1.1). Let $\cJ$ be a rank one local system on $Y$ equipped with $N(\cJ)\,\iso\,\chi$. The complex
$$
A\cJ\otimes e^*\Four_{\psi}(\nu_P^*K)
$$
descends naturally with respect to the map $\Pic Y\to\Pic Y/\Pic X$. The Bessel period of $K$ is the ($\ZZ/2\ZZ$-graded) complex
\begin{equation}
\label{BP_introduction_new}
\RG_c(\Pic Y/\Pic X, \, A\cJ\otimes e^*\Four_{\psi}(\nu_P^*K))
\end{equation}
Our Theorem~~\ref{Th_BP} says that (up to a shift) there is a $\ZZ/2\ZZ$-graded isomorphism between
\begin{equation}
\label{L_function_into_for_Bessel_periods}
\mathop{\oplus}\limits_{d\ge 0} \RG(Y^{(d)}, (\cJ\otimes \phi^*(\pi_*\tilde E^*))^{(d)})[d]
\end{equation}
and the tensor square of (\ref{BP_introduction_new}). The $\ZZ/2\ZZ$-grading on (\ref{L_function_into_for_Bessel_periods}) is given by the parity of $d$. 
  
   From Conjecture~\ref{Con_functoriality} it would follow that $K$ is an automorphic sheaf on $\Bun_G$ corresponding to the local system $E_{\check{G}}$ given by the pair $(\pi_*\tilde E^*, \chi^{-1})$ with symplectic form $\wedge^2(\pi_*\tilde E^*)\to \chi^{-1}$ (cf. also Conjecture~\ref{Con_cusp}). Here $\check{G}$ stands for the Langlands dual group. As predicted by the general philosophy on multiplicity one models, the complex (\ref{L_function_into_for_Bessel_periods}) makes sense for all $\check{G}$-local systems on $X$, this allows us to formulate a conjectural answer for the Bessel periods of all automorphic sheaves on $\Bun_G$ (cf. Conjectures~\ref{Con_BP} and \ref{Con_BP_refined}).  
   
   The geometry suggests that one should be able to recover an automorphic sheaf on $\Bun_G$ from the knowledge of all its Bessel periods (including those for ramified two-sheeted coverings $\phi: Y\to X$). To formulate the corresponding conjecture we switch from $\ell$-adic sheaves to $\cD$-modules (for Section~8 only), as it requires the Fourier-Laumon transform, which is not known in $\ell$-adic setting. 
   
   We also propose one more conjectural construction of automorphic sheaves on $\Bun_G$ as theta-lifting from $\GO_6$ (cf. Conjecture~\ref{Con_autom_GSp4_from_GO6}). 
   
\bigskip\noindent
1.1 {\scshape General notation\ } Let $k$ denote an algebraically closed field of characteristic $p>2$, all the schemes (or stacks) we consider are defined over $k$. Fix a prime $\ell\ne p$. For a scheme (or stack) $S$ write $\D(S)$ for the derived category of $\ell$-adic \'etale sheaves on $S$, and $\P(S)\subset \D(S)$ for the category of perverse sheaves. 
 
 Fix a nontrivial character $\psi: \Fp\to\Qlb^*$ and denote by
$\cL_{\psi}$ the corresponding Artin-Shreier sheaf on $\A^1$. 
Since we are working over an algebraically closed field, we systematically ignore Tate twists. 

 If $V\to S$ and $V^*\to S$ are dual rank $n$ vector bundles over a stack $S$, we normalize the Fourier trasform $\Four_{\psi}: \D(V)\to\D(V^*)$ by 
$\Four_{\psi}(K)=(p_{V^*})_!(\xi^*\cL_{\psi}\otimes p_V^*K)[n]$,  
where $p_V, p_{V^*}$ are the projections, and $\xi: V\times_S V^*\to \A^1$ is the pairing. 

 Let $X$ be a smooth projective connected curve. Write $\Omega$ for the canonical line bundle on $X$. For a smooth scheme of finite type $S$ and a locally free $\cO_S$-module $\cL$ write $\cL^{\star}=\cL^*\otimes\Omega_S$, where $\Omega_S$ is the canonical line bundle on $S$. 
 For a morphism of stacks $f:Y\to Z$ we denote by $\dimrel(f)$ the function of a connected component $C$ of $Y$ given by $\dim C-\dim C'$, where $C'$ is the connected component of $Z$ containing $f(C)$. 
 
 Write $\Bun_k$ for the stack of rank $k$ vector bundles on $X$. For $k=1$ we also write $\Pic X$ for the Picard stack $\Bun_1$ of $X$. We have a line bundle $\cA_k$ on $\Bun_k$ with fibre $\det\RG(X,V)$ at $V\in\Bun_k$. View it as a $\ZZ/2\ZZ$-graded placed in  degree $\chi(V)\!\mod 2$. Our conventions about $\ZZ/2\ZZ$-grading are those of (\cite{Ly}, 3.1).
 
\medskip\noindent
1.2 {\scshape Other results and ideas of proofs}

\medskip\noindent
1.2.1 {\scshape Theta-sheaf} \  Let $G_k$ denote the sheaf of automorphisms of $\cO_X^k\oplus \Omega^k$ preserving the natural symplectic form $\wedge^2(\cO_X^k\oplus \Omega^k)\to \Omega$. The stack $\Bun_{G_k}$ of $G_k$-bundles on $X$ classifies $M\in\Bun_{2k}$ equipped with a symplectic form $\wedge^2 M\to \Omega$. We have a $\mu_2$-gerbe $\Bunt_{G_k}\to\Bun_{G_k}$, where $\Bunt_{G_k}$ is the stack of metaplectic bundles on $X$. In \cite{Ly} we have introduced the theta-sheaf $\Aut=\Aut_g\oplus\Aut_s$ on $\Bunt_{G_k}$ (cf. Section~2.1 for precise definitions). We refer to $\Aut_g$ (resp., to $\Aut_s$) as the generic (resp., special) part of $\Aut$. We write $_X\Aut=\Aut$ when we need to express the dependence on $X$. 
 
 Let $P_k\subset G_k$ be the Siegel parabolic preserving the Lagrangian
subsheaf $\cO_X^k\subset \cO_X^k\oplus\Omega^k$. Write $\nu_k:\Bun_{P_k}\to\Bun_{G_k}$ for the projection, where $\Bun_{P_k}$ is the stack of $P_k$-bundles on $X$. We extend $\nu_k$ to a map $\tilde\nu_k: \Bun_{P_k}\to\Bunt_{G_k}$ (cf. Section~2.1). 
 
 The stack $\Bun_{P_k}$ classifies $L\in\Bun_k$ together with an exact sequence $0\to \Sym^2 L\to ?\to \Omega\to 0$ of $\cO_X$-modules. Let
$^0\Bun_{P_k}\subset \Bun_{P_k}$ be the open substack given by $\H^0(X,\Sym^2 L)=0$. 

 In (\cite{Ly}, Definition~3) we have introduced the complex $S_{P,\psi}$ on $\Bun_{P_k}$ by some explicit construction (cf. 2.1 for details). It was shown in (\select{loc.cit.}, Proposition~7) that there is an isomorphism over $^0\Bun_{P_k}$
$$
\gr_k: S_{P,\psi}\,\iso\,\tilde\nu_k^*\Aut [\dimrel(\nu_k)]
$$
We show that $\gr_k$ extends naturally to an isomorphism over $\Bun_{P_k}$ (cf. Proposition~\ref{Con_everywhere_iso}). 

 Let $\pi: \tilde X\to X$ be an \'etale degree 2 covering, $\Sigma=\Aut_X(\tilde X)=\{1, \sigma\}$ the automorphisms group of $\tilde X$ over $X$. Let $\cE$ be the $\Sigma$-anti-invariants in $\pi_*\cO$, so $\cE$ is equipped with a trivialization $\kappa: \cE^2\,\iso\,\cO$. Let $\cE_0$ be the $\Sigma$-anti-invariants in $\pi_*\Qlb$, it is equipped with $\cE_0^2\,\iso\,\Qlb$. Let $g$ denote the genus of $X$.
 
 Write $\Bun_{G_n,\tilde X}$ for the stack classifying rank $2n$ vector bundles $W$ on $\tilde X$ with symplectic form $\wedge^2 W\to\Omega_{\tilde X}$. Let $\pi_n: \Bun_{G_n,\tilde X}\to\Bun_{G_{2n}}$ be the map sending the above point to $\pi_*W$ equipped with natural symplectic form $\wedge^2 (\pi_*W)\to\Omega$. Let $\Bunt_{G_n,\tilde X}$ denote the corresponding stack of metaplectic bundles. The map $\pi_n$ extends to a map
$\tilde\pi_n:\Bunt_{G_n,\tilde X}\to\Bunt_{G_{2n}}$ (cf. 3.5). 
We establish a canonical isomorphism  
$$
{_{\tilde X}\Aut}\,\iso\, \tilde\pi_n^*\Aut[\dimrel(\tilde\pi_n)]
$$
preserving the generic and special parts (cf. Proposition~\ref{Pp_Weil_and_coverings}).

\medskip 
\noindent
1.2.2  {\scshape Theta-lifting functors} 

\smallskip\noindent
Let $n,m\in\NN$ and $G=\GSp_{2n}$. Let $\HH=\GO_{2m}^0$ denote the connected component of unity of the split orthogonal similitude group $\GO_{2m}$ over $\Spec k$.  Pick a maximal torus and a Borel subgroup $\TT_{\HH}\subset\BB_{\HH}\subset\HH$. 
We pick an involution $\tilde\sigma\in\OO_{2m}(k)$ preserving $\TT_{\HH}$ and $\BB_{\HH}$ such that $\tilde\sigma\notin\SO_{2m}$. So, for $m\ge 2$ (and $m\ne 4$) it induces the unique nontrivial automorphism of the Dynkin diagram of $\HH$. Consider the corresponding $\Sigma$-action on $\GO_{2m}^0$ by conjugation. Let $\tilde H$ be the group scheme on $X$, the twisting of $\GO_{2m}^0$ by the $\Sigma$-torsor $\pi:\tilde X\to X$. 

  The stack $\Bun_G$ of $G$-torsors on $X$ classifies $M\in\Bun_{2n}, \cA\in\Bun_1$ with symplectic form $\wedge^2 M\to\cA$. The stack $\Bun_{\tilde H}$ of $\tilde H$-torsors on $X$ classifies $V\in\Bun_{2m},  \, \cC\in\Bun_1$, a nondegenerate symmetric form $\Sym^2 V\to\cC$, and a compatible trivialization $\gamma: \cC^{-m}\otimes\det V\,\iso\,\cE$. This means that the composition 
$$
\cC^{-2m}\otimes(\det V)^2\,\toup{\gamma^2}\,\cE^2\,\iso\,\cO
$$ 
is the isomorphism induced by $V\,\iso\, V^*\otimes\cC$. 

 Let $\RCov^0$ denote the stack classifying a line bundle $\cU$ on $X$ together with a trivialization $\cU^{\otimes 2}\,\iso\,\cO$. Its connected components are indexed by 
$\H^1_{et}(X, \ZZ/2\ZZ)$.
 
 Let $\Bun_H$ be the stack classifying $V\in\Bun_{2m}, \cC\in\Bun_1$ and
a symmetric form $\Sym^2 V\to \cC$ such that the corresponding trivialization $(\cC^{-m}\otimes\det V)^2\,\iso\,\cO$ lies in the component of $\RCov^0$ given by $(\cE,\kappa)$. Note that 
$$
\Bun_{\tilde H}\,\iso\, \Spec k\times_{\RCov^0} \Bun_H,
$$
where the map $\Spec k\to \RCov^0$ is given by $(\cE,\kappa)$. 
The projection $\rho_H:\Bun_{\tilde H}\to \Bun_H$ is a $\mu_2$-torsor. 
 
 Write $\Bun_H^d\subset\Bun_H$ for the open substack given by $\deg\cC=d$, and similarly for $\Bun_{\tilde H}^d$. 
 Set
$$
\Bun_{G,H}=\Bun_H\times_{\Pic X} \Bun_G, 
$$
where  the map $\Bun_H\to \Pic X$ sends $(V,\cC, \Sym^2 V\to\cC)$ to $\Omega\otimes\cC^{-1}$. The map $\Bun_G\to\Pic X$ sends $(M, \wedge^2 M\to\cA)$ to $\cA$. We have an isomorphism
$\cC\otimes \cA\,\iso\, \Omega$ for a point of $\Bun_{G,H}$. 
Let
$$
\tau: \Bun_{G,H}\to \Bun_{G_{2nm}}
$$ 
be the map sending a point as above to $V\otimes M$ with symplectic form $\wedge^2(V\otimes M)\to\Omega$. We extends $\tau$ to a map $\tilde\tau: \Bun_{G,H}\to \Bunt_{G_{2nm}}$ (cf. 3.2). Let $\Bun_{G,\tilde H}$ be the stack obtained from $\Bun_{G,H}$ by the base change $\Bun_{\tilde H}\to\Bun_H$. 

 Viewing $\tilde\tau^*\Aut$ as a kernel of integral operators, we define functors
$F_G: \D(\Bun_H)\to\D(\Bun_G)$ and $F_H: \D(\Bun_G)\to\D(\Bun_H)$, set also $F_{\tilde H}=\rho_H^*\comp F_H$ (cf. 3.2). 
 
\medskip\noindent
1.2.3  {\scshape The pair $\GO_2, \GL_2$} 

\medskip\noindent  
Assume $n=m=1$ and $\tilde X$ connected. In this case we prove Conjecture~\ref{Con_functoriality} for the functor $F_G$. To do so, we first prove Theorem~\ref{Th_Hecke_property} saying how the action of Hecke operators on $\tilde\tau^*\Aut$ with respect to $G$ is expressed in terms of the similar action with respect to $\tilde H$. This is a global geometric analog of a particular case of the theorem of Rallis (\cite{R}, cf. also \cite{Ly4}). 

 We also show that both
$\tilde\tau^*\Aut_g[\dimrel(\tilde\tau)]$ and $\tilde\tau^*\Aut_s[\dimrel(\tilde\tau)]$ are (self-dual) irreducible perverse sheaves on each connected component of $\Bun_{G,\tilde H}$ (cf. Proposition~\ref{Pp_irreducibility}), and the functor $F_G:\D(\Bun_{\tilde H})\to\D(\Bun_G)$ commutes with the Verdier duality. 
 
  If $\tilde E$ is a rank one local system on $\tilde X$, let $K_{\tilde E}$ denote the automorphic sheaf on $\Bun_{\tilde H}\,\iso\,\Pic\tilde X$ corresponding to $\tilde E$. Then $F_G(K_{\tilde E})$ is an automorphic sheaf on $\Bun_G$ corresponding to the local system $E=(\pi_*\tilde E)^*$. We check that (up to a tensoring by a 1-dimensional vector space) the sheaf $F_G(K_{\tilde E})$ coincides with the perverse sheaf $\Aut_E$ constructed via Whittaker models in \select{loc.cit.} (cf. Proposition~\ref{Pp_lift_for_GL2}). 
 
 Theorem~\ref{Th_Hecke_property} also allows us to calculate the following Rankin-Selberg type convolution (we need it for our proof of Theorem~\ref{Th_WP}). Let $E$ be an irreducible rank 2 local system on $X$, $E_1$ be a rank one local system on $X$. We denote by $\Aut_{E_1\oplus \Qlb}$ the corresponding geometric Eisenstein series (cf. 4.3). Our Theorem~\ref{Th_Rankin_Selberg_GO_2_GL_2} provides an explicit calculation of $F_{\tilde H}(\Aut_{E_1\oplus\Qlb}\otimes\Aut_E)$. The method of its proof is inspired by (\cite{Ly3}). We don't know if this Rankin-Selberg convolution was known before in classical theory of automorphic forms.  
 
\bigskip\noindent
1.2.4 {\scshape Waldspurger periods}   
 
\medskip\noindent
Let us explain how we calculate the Waldspurger periods (Theorem~\ref{Th_WP}). Mainly, we follow the approach of Waldspurger (\cite{Wa}), but there are some new phenomena in geometric settings.

 Let $n=2$, so $G=\GL_2$. Let $E$ be an irreducible rank two local system on $X$, $\Aut_E$ be the corresponding automorphic sheaf on $\Bun_G$. Take both $\pi:\tilde X\to X$ and $\tilde H=\GO_4^0$ split. Remind the perverse sheaf $K_{\pi^*E,\det E,\tilde H}$ on $\Bun_{\tilde H}$
from Section~1.0. First, we identify $K_{\pi^*E,\det E,\tilde H}$ with the theta-lift $F_{\tilde H}(\Aut_{E^*})$ from $\Bun_G$ (cf. Proposition~\ref{Pp_Shimizu}). This is a geometric version of a Theorem of Shimizu (\cite{Wa}). 

 Then we consider the diagram
$$
\begin{array}{ccccc}
\Pic Y & \getsup{m} & \Pic Y\times\Pic Y & \toup{\phi_1\times\phi_1} & \Bun_2\times\Bun_2\\
 & \nwarrow\lefteqn{\scriptstyle \gp_{R_{\phi}}} & \downarrow &&\downarrow\\
 && \Bun_{R_{\phi}} & \toup{\gq_{R_{\phi}}} & \Bun_{\tilde H},
\end{array}
$$  
where $m$ is the tensor product map (followed by the automorphism sending $\cB\in\Pic Y$ to $\cB^*\otimes\Omega_Y$), and the vertical arrows are, roughly speaking, the quotients by the action of $\Pic X$, where $\cL\in\Pic X$ sends $(L_1,L_2)\in\Bun_2\times\Bun_2$ to $L_1\otimes\cL, L_2\otimes\cL^*$. Remind that $\phi_1$ sends $\cB\in\Pic Y$ to $\phi_*\cB\in\Bun_2$. 

 The key step is Theorem~\ref{Th_special_case_BP} that calculates the complex $(\gp_{R_{\phi}})_!\gq_{R_{\phi}}^*F_{\tilde H}(\Aut_{E^*})$ explicitely in terms of $E$ and $\phi:Y\to X$ (in our actual formulation of Theorem~\ref{Th_special_case_BP} the covering $\pi:\tilde X\to X$ may be nonsplit). We derive Theorem~\ref{Th_special_case_BP} from the properties of the theta-lifting between $\GO_2$ and $\GL_2$ (Proposition~\ref{Pp_lift_for_GL2}) combined with our Rankin-Selberg convolution result (Theorem~\ref{Th_Rankin_Selberg_GO_2_GL_2}).  

 Let us indicate at this point that the existence of the geometric Waldspurger periods of automorphic sheaves (the fact that condition ($C_W$) in Definition~\ref{Def_WP} holds) is a consequence of an intriguing acyclicity result (Theorem~\ref{Th_ULA}, Section~6.1.2). It says that the Hecke property of a given automorphic sheaf $\cS$ on $\Bun_2$ already implies that $\cS$ is universally locally acyclic (ULA) over `the moduli of spectral curves'. This allows to control perversity of the complexes $\cK_E$ in Theorem~\ref{Th_WP} (and similarly for Bessel periods in Theorem~\ref{Th_BP}). 
  
  We also formulate conjectural answers for the geometric Waldspurger periods (Conjecture~\ref{Con_WP}) and the geometric Bessel periods (Conjecture~\ref{Con_BP}) of all automorphic sheaves. Besides, we verify Conjecture~\ref{Con_WP} for geometric Eisenstein series on $\Bun_2$ (cf. Proposition~\ref{Pp_Eis_WP}). 

\bigskip\noindent
1.2.5 {\scshape Case $H=\GO_6$}  
 
\medskip\noindent
Assume $m=3$ and $\tilde X$ split, so $\tilde H=\GO_6^0$. Let $n=2$, so $G=\GSp_4$. Let $E_{\check{G}}$ be a $\check{G}$-local system on $X$ viewed as a pair $(E,\chi)$, where $E$ (resp., $\chi$) is a rank 4 (resp., rank one) local system on $X$ with symplectic form $\wedge^2 E\to\chi$. Assume $E$ irreducible. Remind that $\check{G}$ is a subgroup of $\GSpin_6$, which is the Langlands dual to $\tilde H$. 
We define the perverse sheaf $K_{E,\chi,\tilde H}$ on $\Bun_{\tilde H}$ corresponding to a $\GSpin_6$-local system $(E,\chi)$. We conjecture that 
\begin{equation}
\label{complex_intro_fromGO6}
F_G(K_{E^*, \chi^*,\tilde H})
\end{equation}
is an automorphic sheaf on $\Bun_G$ corresponding to $E_{\check{G}}$ (cf. Conjecture~\ref{Con_autom_GSp4_from_GO6}). 
We show that the geometric Bessel periods of (\ref{complex_intro_fromGO6}) are essentially the generalized Waldspurger periods of $K_{E^*, \chi^*,\tilde H}$ (cf. Proposition~\ref{Pp_GWP_and_BP}).

\bigskip\bigskip
\centerline{\scshape 2. Theta-sheaf}

\bigskip\noindent
2.1 Let $G_k$ denote the sheaf of automorphisms of $\cO_X^k\oplus \Omega^k$ preserving the natural symplectic form $\wedge^2(\cO_X^k\oplus \Omega^k)\to \Omega$.
The stack $\Bun_{G_k}$ of $G_k$-bundles on $X$ classifies $M\in\Bun_{2k}$ equipped with a symplectic form $\wedge^2 M\to \Omega$. Write $\cA_{G_k}$ for the line bundle on $\Bun_{G_k}$ with fibre $\det\RG(X,M)$ at $M$. We view it as a $\ZZ/2\ZZ$-graded line bundle (purely of degree zero).  Denote by $\Bunt_{G_k}\to\Bun_{G_k}$ the $\mu_2$-gerbe of square roots of $\cA_{G_k}$. 

 Remind the definition of the theta-sheaf $\Aut$ on $\Bunt_{G_k}$ from \cite{Ly}. Let $_i\Bun_{G_k}\subset \Bun_{G_k}$ be the locally closed substack given by $\dim\H^0(X,M)=i$ for $M\in\Bun_{G_k}$. Let $_i\Bunt_{G_k}$ denote the preimage of $_i\Bun_{G_k}$ in $\Bunt_{G_k}$. 
 
  Write $_i\cB$ for the line bundle on $_i\Bun_{G_k}$ whose fibre at $M\in{_i\Bun_{G_k}}$ is $\det\H^0(X, M)$. View $_i\cB$ as $\ZZ/2\ZZ$-graded placed in degree $\dim\H^0(X, M)$ modulo $2$. For each $i\ge 0$ we have a canonical $\ZZ/2\ZZ$-graded isomorphism 
$$
_i\cB^2\,\iso\, \cA\mid_{_i\Bunt_{G_k}}
$$ 
It yields a two-sheeted covering 
$_i\rho: {_i\Bun_{G_k}}\to {_i\Bunt_{G_k}}$ locally trivial in \'etale topology. Define a local system $_i\Aut$ on $_i\Bunt_{G_k}$ by
$$
_i\Aut=\Hom_{S_2}(\sign, \; {_i\rho_!}\Qlb)
$$
The perverse sheaf $\Aut_g\in \P(\Bunt_{G_k})$ (resp., $\Aut_s\in \P(\Bunt_{G_k})$) is defined as the intermediate extension of
$_0\Aut[\dim\Bun_G]$ (resp., of $_1\Aut[\dim\Bun_G-1]$) under $_i\Bunt_{G_k}\hook{} \Bunt_{G_k}$. The theta-sheaf $\Aut$ is defined by
$$
\Aut=\Aut_g\oplus\Aut_s
$$

 Let $P_k\subset G_k$ be the Siegel parabolic preserving the Lagrangian
subsheaf $\cO_X^k\subset \cO_X^k\oplus\Omega^k$. Write $Q_k$ for the Levi quotient of $P_k$, so $Q_k\,\iso\, \GL_k$ canonically.

 Write $\nu_k:\Bun_{P_k}\to\Bun_{G_k}$ for the projection. As in (\select{loc.cit.}, 5.1), we extend it to a map $\tilde\nu_k:\Bun_{P_k}\to\Bunt_{G_k}$ defined as follows. The stack $\Bun_{P_k}$ classifies $L\in\Bun_k$ together with an exact sequence $0\to \Sym^2 L\to ?\to \Omega\to 0$ of $\cO_X$-modules. The induced exact sequence $0\to L\to M\to L^*\otimes\Omega\to 0$ yields an isomorphism of $\ZZ/2\ZZ$-graded lines
$$
\det\RG(X,M)\,\iso\, \det\RG(X,L)\otimes\det\RG(X, L^*\otimes\Omega)\,\iso\,  \det\RG(X,L)^{\otimes 2}
$$
The map $\tilde\nu_k$ sends the above point to $(\cB, M)$, where $\cB=\det\RG(X, L)$ is equipped with the above isomorphism $\cB^2\,\iso\, \det\RG(X,M)$.  
 
   Remind the definition of the complex $S_{P,\psi}$ on $\Bun_{P_k}$ (\select{loc.cit.}, Definition~3). Denote by $\cV\to\Bun_k$ the stack whose fibre over $L\in\Bun_k$ is $\Hom(L,\Omega)$. Write $\cV_2\to\Bun_k$ for the stack whose fibre over $L\in\Bun_k$ is $\Hom(\Sym^2 L, \Omega^2)$. We have a projection $\pi_2:\cV\to \cV_2$ sending $s\in\Hom(L,\Omega)$ to $s\otimes s\in\Hom(\Sym^2 L, \Omega^2)$. We set
$$
S_{P,\psi}\,\iso\, \Four_{\psi}(\pi_{2 !}\Qlb)[\dimrel],
$$
where $\Four_{\psi}: \D(\cV_2)\to\D(\Bun_P)$ denotes the Fourier transform functor,  and $\dimrel$ is the function of a connected component of $\Bun_k$ given by $\dimrel(L)=\dim\Bun_k-\chi(L)$, $L\in\Bun_k$.   

 The group $S_2$ acts on $\cV$ sending $(L, s: L\to\Omega)$ to $(L, -s)$. This gives rise to a $S_2$-action on $S_{P,\psi}$. 
By (\select{loc.cit.}, Remark~3), the $S_2$-invariants of $S_{P,\psi}$ are $S_{P,\psi,g}$ (resp., $S_{P,\psi,s}$) over the connected component of $\Bun_{P_k}$ with $\chi(L)$ even (resp., odd).  
  
  Let $^0\Bun_{P_k}\subset \Bun_{P_k}$ be the open substack given by $\H^0(X,\Sym^2 L)=0$. 
By (\select{loc.cit.}, Proposition~7) there is an isomorphism\footnote{the isomorphism $\gr_k$ is not canonical: once $\sqrt{-1}\in k$ is chosen, $\gr_k$ is well-defined up to a sign.}
\begin{equation}
\label{iso_gr_k}
\gr_k: S_{P,\psi}\,\iso\,\tilde\nu_k^*\Aut [\dimrel(\nu_k)]
\end{equation}
over $^0\Bun_{P_k}$, here $\dimrel(\nu_k)=\dim\Bun_{P_k}-\dim\Bun_{G_k}$ is a function of a connected component of $\Bun_{P_k}$. From (\select{loc.cit.}, Sect.~2) it may be deduced that, in the case of a finite base field, the function `trace of Frobenius' of $S_{P,\psi}$ descends with respect to $\tilde\nu_k:\Bun_{P_k}\to\Bunt_{G_k}$ over the whole of $\Bun_{P_k}$. We claim that it is also true in the geometric setting.
\begin{Pp} 
\label{Con_everywhere_iso}
The isomorphism $\gr_k$ extends naturally to an isomorphism over $\Bun_{P_k}$.
\end{Pp}
 
\noindent 
\select{Proof}

\Step 1  For an effective divisor $D$ on $X$ denote by $_{D,P}\Bun_{G_k}$ the stack classifying $M\in\Bun_{G_k}$ together with a $P_k$-structure on $M\mid_D$. A point of $_{D,P}\Bun_{G_k}$ is given by $M\in\Bun_{G_k}$ together with a lagrangian $\cO_D$-submodule $L_D\subset M\mid_D$. Denote by $p_D: \Bun_{P_k}\to {_{D,P}\Bun_{G_k}}$ the map
sending $(L\subset M)\in\Bun_{P_k}$ to $(M, L_D)$ with $L_D=L\mid_D$. 
Let $\nu_D: {_{D,P}\Bun_{G_k}}\to\Bun_{G_k}$ be the projection. 

Pick a point $x\in X$ and a nonnegative integer $i$. Set $D=ix$. Let $a_D:\Bun_{P_k}\to\Bun_{P_k}$ be the map sending an exact sequence $0\to \Sym^2 L\to ?\to\Omega\to 0$ to its push-forward with respect to the map $\Sym^2 L\hook{} \Sym^2 (L(D))$. Since $\Hom(\Omega, \, \Sym^2 (L(D))/\Sym^2 L)$ acts freely and transitively on a fibre of $a_{D}$, the map $a_D$ is an affine fibration of rank $k(k+1)i$. 

 We are going to establish a canonical isomorphism 
$$
(a_D)_! S_{P,\psi}\,\iso\, S_{P,\psi}[-k^2 i]
$$ 
To do so, write $\pi_2^{D}: \cV^{D}\to \cV^{D}_2$ for the map obtained from $\pi_2:\cV\to \cV_2$ by the base change $\Bun_k\to\Bun_k$ sending $L$ to $L(D)$. So, a fibre of $\cV^{D}_2\to\Bun_k$ is $\Hom(\Sym^2 (L(D)), \Omega^2)$, and we have a natural map 
$$
^ta_{D}: \cV^{D}_2\to \cV_2,
$$ 
the transpose of $a_{D}$. Since $(^t a_{D})^* \pi_{2 !}\Qlb\,\iso\, (\pi_2^{D})_!\Qlb$ canonically,  our assertion follows from the standard properties of the Fourier transform functor. 
 
\smallskip 
\Step 2  Denote by $_x\cH_{G_k}$ the Hecke stack classifying $M,M'\in\Bun_{G_k}$ together with an isomorphism of $G_k$-torsors $M\,\iso\, M'\mid_{X-x}$. Let $T_k\subset Q_k=\GL_k$ denote the maximal torus of diagonal matrices, its coweight lattice identifies with $\ZZ^k$. The preimage of the standard Borel subgroup of $Q_k$ in $P_k$ is a Borel subgroup of $G_k$, this also fixes the set of simple roots of $G_k$.  Set $\omega_k=(1,\ldots,1)$, where 1 appears $k$ times. This is a dominant coweight of $G_k$ orthogonal to all the roots of $Q_k$.  

 Denote by $_{D,P}\cH_{G_k}$ the stack classifying $(M, M', \;  M'\,\iso\,M\mid_{X-x})\in {_x\cH_{G_k}}$ such that $M'$ is in the position $i\omega_k$ with respect to $M$ at $x$, $L_D\subset M\mid_D$ with
$(M, L_D)\in {_{D,P}\Bun_{G_k}}$ satisfying $L_D\cap (M'/M(-D))=0$.
The latter intersection is taken inside $M(D)/M(-D)$, it makes sense because $L_D\subset M/M(-D)$ and 
$$
M(-D)\subset M'\subset M(D)
$$

 Denote by $a_{\cH,D}: {_{D,P}\cH_{G_k}}\to {_{D,P}\Bun_{G_k}}$ the map sending the above point to $(M, L_D)$. We have a diagram, where the square is cartesian 
\begin{equation}
\label{diag_first_Pp_one}
\begin{array}{ccccc}
\Bun_{P_k} & \toup{p_{\cH,D}} &  _{D,P}\cH_{G_k} & \toup{\nu_{\cH,D}} & \Bun_{G_k}\\
\downarrow\lefteqn{\scriptstyle a_D} &&  \downarrow\lefteqn{\scriptstyle a_{\cH,D}}\\
\Bun_{P_k} & \toup{p_D} & {_{D,P}\Bun_{G_k}} & \toup{\nu_D} & \Bun_{G_k}
\end{array}
\end{equation}
Here $p_{\cH,D}$ is the map sending  $(L'\subset M')$ to $(M', M, L_D)$, where $(L\subset M)$ is the image of $(L'\subset M')$ under $a_D$, and $L_D=L\mid_D$.  The map $\nu_{\cH,D}: {_{D,P}\cH_{G_k}}\to\Bun_{G_k}$ sends the above point to $M'$. 

 Consider the diagram
$$
_{D,P}\tilde\cH_{G_k} \toup{\tilde a_{\cH,D}} {_{D,P}\Bunt_{G_k}} \toup{\tilde\nu_D} \Bunt_{G_k}
$$ 
obtained from $_{D,P}\cH_{G_k} \toup{a_{\cH,D}} {_{D,P}\Bun_{G_k}} \toup{\nu_D} \Bun_{G_k}$ by the base change $\Bunt_{G_k}\to\Bun_{G_k}$. 
Now $\nu_{\cH,D}$ lifts to a map 
$$
\tilde\nu_{\cH,D}: {_{D,P}\tilde\cH_{G_k}}\to \Bunt_{G_k}
$$ 
defined as follows.  A point of $_{D,P}\tilde \cH_{G_k}$ is given by $(\cB, M, M', L_D\subset M\mid_D)$, where $\cB$ is a 1-dimensional ($\ZZ/2\ZZ$-graded purely of degree zero) vector space equipped with a $\ZZ/2\ZZ$-graded isomorphism $\cB^2\,\iso\, \det\RG(X, M)$. The map $\tilde\nu_{\cH,D}$ sends this point to $(\cB', M')$, where $\cB'=\cB\otimes\det_k (L_D)^{-1}$ is equipped with an isomorphism
\begin{equation}
\label{iso_Pp_one_for_gerbs}
(\cB')^2\,\iso\,\det\RG(X,M')
\end{equation}
that we are going to define. 

 For a vector bundle $\cN$ on $X$ write $\cN_D=\cN/\cN(-D)$. 
We have an exact sequence of $\cO_D$-modules 
$$
0\to (M\cap M')/M(-D)\to
M'/M(-D)\to (M+M')/M\to 0
$$ 
Note also that $(M+M')/M$ is the orthogonal complement of $(M\cap M')/M(-D)$ with respect to the perfect pairing of $\cO_D$-modules $M_D\otimes M(D)_D\to \Omega(D)_D$ given by the symplectic form. 
By our assumptions, $(M\cap M')/M(-D)\subset M_D$ is a lagrangian $\cO_D$-submodule such that 
$$
((M\cap M')/M(-D)) \oplus L_D\,\iso\, M_D\,\iso\, (((M+M')/M)\otimes\cO(-D))\oplus L_D
$$
The exact sequences of $\cO_X$-modules $0\to M(-D)\to M\to M_D\to 0$ and $0\to M(-D)\to M'\to M'/M(-D)\to 0$ yield $\ZZ/2\ZZ$-graded isomorphisms
$$
\det\RG(X, M')\,\iso\, \frac{\det\RG(X,M)}{\det_k (L_D)}\otimes
\frac{\det_k (M(D)_D)}{\det_k (L_D\otimes \cO(D))}\,\iso\, 
\frac{\det\RG(X,M)}{\det_k (L_D)^{\otimes 2}}
$$
giving rise to (\ref{iso_Pp_one_for_gerbs}). 


 To summarize, the diagram (\ref{diag_first_Pp_one}) is refined to the following commutative diagram
\begin{equation}
\label{diag_Pp_one_to_summarize}
\begin{array}{ccccccc}
&& \Bun_{P_k} & \getsup{a_D} & \Bun_{P_k}\\
& \swarrow\lefteqn{\scriptstyle \tilde\nu_k} & \downarrow\lefteqn{\scriptstyle \tilde p_D} && \downarrow\lefteqn{\scriptstyle \tilde p_{\cH,D}} & \searrow\lefteqn{\scriptstyle \tilde\nu_k}\\
\Bunt_{G_k} & \getsup{\tilde\nu_D} & _{D,P}\Bunt_{G_k} & \getsup{\tilde a_{\cH,D}} & _{D,P}\tilde\cH_{G_k} & \toup{\tilde\nu_{\cH,D}} & \Bunt_{G_k},
\end{array}
\end{equation}
where the middle square is cartesian. Here $\tilde p_D$ is the product map $(p_D\times \tilde\nu_k)$, and $\tilde p_{\cH,D}$ is the product map $(p_{\cH,D}\times \tilde\nu_k)$. 

\smallskip
\Step 3 Set $^{0,i}\Bun_P=a_D(^0\Bun_P)$, this is an open substack of $\Bun_P$. For $i\le j$ we have $^{0,i}\Bun_P\subset {^{0,j}\Bun_P}$ and the union of all $^{0,i}\Bun_P$ equals $\Bun_P$. We are going to extend $\gr_k$ to each $^{0,i}\Bun_P$ in a compatible way.

 Now (\ref{iso_gr_k}) and the diagram (\ref{diag_Pp_one_to_summarize}) yield an isomorphism over $^{0,i}\Bun_{P_k}$
\begin{equation}
\label{iso_used_Pp_one_for_descent}
S_{P,\psi}[-k^2 i]\,\iso\, (a_{D})_! S_{P,\psi}\,\iso\, \tilde p_D^*(\tilde a_{\cH,D})_!(\tilde\nu_{\cH,D})^* \Aut[\dimrel], 
\end{equation}
where $\dimrel=\dim\Bun_{P_k}+\dimrel(a_D)-\dim\Bun_{G_k}$ and $\dimrel(a_D)=k(k+1)i$.  

Restricting (\ref{iso_used_Pp_one_for_descent}) to the open substack $^0\Bun_P\subset {^{0,i}\Bun_P}$ and applying  (\ref{iso_gr_k}) once again, we get an isomorphism of (shifted) perverse sheaves over $^0\Bun_P$
\begin{equation}
\label{iso_Pp_one_step_3}
\tilde p_D^* \tilde\nu_D^*\Aut \,\iso\, 
\tilde p_D^*(\tilde a_{\cH,D})_!(\tilde\nu_{\cH,D})^* \Aut [k(k+1)i+k^2 i]
\end{equation}
 
\Step 4  Denote by $^0_D\Bun_P\subset {^0\Bun_P}$ the open substack given by $\H^0(X, (\Sym^2 L)(D))=0$. Let us show that the map 
$p_D: {^0_D\Bun_{P_k}}\to {_{D,P}\Bun_{G_k}}$ is smooth.

 Set $\gp=\Lie P_k$ and $\gg=\Lie G_k$. Let $\cF_{P_k}$ be a $k$-point of $\Bun_{P_k}$ given by $(L\subset M)$. Let $K$ denote the kernel of  the composition 
$$\gg_{\cF_{P_k}}\to  (\gg_{\cF_{P_k}})\mid_D\to \Hom_{\cO_D}(L_D, M_D/L_D)
$$

  Remind the following notion. For a 1-morphism $\Spec k\toup{x}\cX$ to a stack $\cX$ the tangent groupoid to $x$ is the category, whose objects are pairs $(x_1,\alpha_1)$, where $x_1$ is a 1-morphism
$\Spec k[\epsilon]/\epsilon^2 \to\cX$ and $\alpha$ is a 2-morphism $x\to \bar x_1$. Here  $\bar x_1$ is the composition $\Spec k\hook{} \Spec k[\epsilon]/\epsilon^2 \toup{x_1}\cX$. A morphism from $(x_1, \alpha_1)$ to $(x_2,\alpha_2)$ is a 2-morphism $\beta: x_1\to x_2$ such that the diagram commutes
$$
\begin{array}{ccc}
\bar x_1 &  \toup{\bar\beta} & \bar x_2\\
\uparrow\lefteqn{\scriptstyle \alpha_1} & \nearrow\lefteqn{\scriptstyle \alpha_2}\\
x
\end{array}
$$

The tangent groupoid to $_{D,P}\Bun_{G_k}$ at the $k$-point $p_D(\cF_{P_k})$ is isomorphic to the stack quotient of $\H^1(X,K)$ by the trivial action of $\H^0(X,K)$. The natural map $\gp_{\cF_{P_k}}\to \gg_{\cF_{P_k}}$ factors through $K\subset \gg_{\cF_{P_k}}$. We need to show that $\H^1(X, \gp_{\cF_{P_k}})\to \H^1(X, K)$ is surjective. We have an exact sequence $0\to K/(\gp_{\cF_{P_k}})\to (\gg/\gp)_{\cF_{P_k}}\to (\Omega\otimes \Sym^2 L^*)_D\to 0$. So, $K/(\gp_{\cF_{P_k}})\,\iso\, (\Omega\otimes \Sym^2 L^*)(-D)$, the desired surjectivity follows. 

 It is easy to deduce that $\tilde p_D: {^0_D\Bun_{P_k}}\to {_{D,P}\Bunt_{G_k}}$ is smooth. One checks that it is also surjective and has connected fibres. So, (\ref{iso_Pp_one_step_3}) descends to an isomorphism of (shifted) perverse sheaves on $_{D,P}\Bunt_{G_k}$
$$
\tilde\nu_D^*\Aut \,\iso\, 
(\tilde a_{\cH,D})_!(\tilde\nu_{\cH,D})^* \Aut[k(k+1)i+k^2 i]
$$

 Now from (\ref{iso_used_Pp_one_for_descent}) we get an isomorphism over $^{0,i}\Bun_{P_k}$
$$
S_{P,\psi}\,\iso\, \tilde p_D^* \tilde\nu_D^*\Aut [\dimrel(\nu_k)]
\eqno{\square}
$$

\medskip

 For the rest of the paper we fix the isomorphism (\ref{iso_gr_k}) over $\Bun_{P_k}$, some of our results will depend on this choice. 

\bigskip\bigskip
\centerline{\scshape 3. Theta-lifting for the pair $\GSp_{2n}, \GO_{2m}$} 

\bigskip\noindent
3.1  Let $n,m\in\NN$ and $\GG=G=\GSp_{2n}$. Pick a maximal torus and a Borel subgroup $\TT_{\GG}\subset \BB_{\GG}\subset \GG$. The stack $\Bun_G$ classifies $M\in\Bun_{2n}, \cA\in\Bun_1$ with symplectic form $\wedge^2 M\to\cA$.
We have a ($\ZZ/2\ZZ$-graded) line bundle $\cA_G$ on $\Bun_G$ with fibre $\det\RG(X,M)$ at $(M,\cA)$. 
 
 Let $\pi: \tilde X\to X$ be an \'etale degree 2 covering, $\sigma$ the nontrivial automorphism of $\tilde X$ over $X$ and $\Sigma=\{1,\sigma\}$. Let $\cE$ be the $\sigma$-anti-invariants in $\pi_*\cO$, it is equipped with a trivialization $\kappa: \cE^2\,\iso\,\cO$. Let $\cE_0$ denote the $\sigma$-anti-invariants in $\pi_*\Qlb$, it is equipped with $\cE_0^2\,\iso\,\Qlb$. Let $g$ (resp., $\tilde g$) denote the genus of $X$ (resp., of $\tilde X$). 
 
  Let $\HH:=\GO_{2m}^0$ be the connected component of unity of the split orthogonal similitude group $\GO_{2m}$ over $\Spec k$. Pick a maximal torus and a Borel subgroup $\TT_{\HH}\subset \BB_{\HH}\subset \HH$. Pick $\tilde\sigma\in\OO_{2m}(k)$ with $\tilde\sigma^2=1$ such that $\tilde\sigma\notin \SO_{2m}(k)$. We assume in addition that $\tilde\sigma$ preserves $\TT_{\HH}$ and $\BB_{\HH}$, so for $m\ge 2$ it induces the unique\footnote{except for $m=4$. The group $\GO_8$ also has trilitarian outer forms, we do not consider them.} nontrivial automorphism of the Dynkin diagram of $\HH$.  
For $m=1$ we identify $\HH\,\iso\,\Gm\times\Gm$ in such a way that $\tilde\sigma$ permutes the two copies of $\Gm$. 

 Realize $\HH$ as the subgroup of $\GL(k^{2m})$ preserving up to a multiple the symmetric form given by the matrix
$$
\left(
\begin{array}{cc}
0  & E_m\\
E_m & 0
\end{array}
\right),
$$ 
where $E_m\in\GL_m$ is the unity. Take $\TT_H$ to be the maximal torus of diagonal matrices, $\BB_{H}$ the Borel subgroup preserving for $i=1,\ldots,m$ the isotropic subspace generated by the first $i$ base vectors $\{e_1,\ldots,e_i\}$. Then one may take $\tilde\sigma$ interchanging $e_m$ and $e_{2m}$ and acting trivially on the orthogonal complement to $\{e_m, e_{2m}\}$. 
 
  Consider the corresponding $\Sigma$-action on $\HH$ by conjugation. Let $\tilde H$ be the group scheme on $X$, the twisting of $\HH$ by the $\Sigma$-torsor $\pi:\tilde X\to X$. 
  
  The stack $\Bun_{\tilde H}$ classifies: $V\in\Bun_{2m},  \, \cC\in\Bun_1$, a nondegenerate symmetric
form $\Sym^2 V\to\cC$, and a compatible trivialization $\gamma: \cC^{-m}\otimes\det V\,\iso\,\cE$. This means that the composition 
$$
\cC^{-2m}\otimes(\det V)^2\,\toup{\gamma^2}\,\cE^2\,\iso\,\cO
$$ 
is the isomorphism induced by $V\,\iso\, V^*\otimes\cC$. (Though $\det$ is involved, we view $\gamma$ as ungraded). 


 Let $\RCov^0$ denote the stack classifying a line bundle $\cU$ on $X$ together with a trivialization $\cU^{\otimes 2}\,\iso\,\cO$. Its connected components are indexed by 
$\H^1_{et}(X, \ZZ/2\ZZ)$, each connected component is isomorphic to the classifying stack $B(\mu_2)$. 
 
 Let $\Bun_H$ be the stack classifying $V\in\Bun_{2m}, \cC\in\Bun_1$ and
a symmetric form $\Sym^2 V\to \cC$ such that the corresponding trivialization $(\cC^{-m}\otimes\det V)^2\,\iso\,\cO$ lies in the component of $\RCov^0$ given by $(\cE,\kappa)$. Note that 
$$
\Bun_{\tilde H}\,\iso\, \Spec k\times_{\RCov^0} \Bun_H,
$$
where the map $\Spec k\to \RCov^0$ is given by $(\cE,\kappa)$. 
Write $\Bun_H^d\subset\Bun_H$ for the open substack given by $\deg\cC=d$, and similarly for $\Bun_{\tilde H}^d$. 

The projection $\rho_H:\Bun_{\tilde H}\to \Bun_H$ is a $\mu_2$-torsor. By extension of scalars $\mu_2\subset \Qlb^*$ it yields a rank one local system $\cN$ on $\Bun_H$, which we refer to as \select{the determinantal local system}. 

 Let $\cA_H$ be the ($\ZZ/2\ZZ$-graded) line bundle on $\Bun_H$ with fibre $\det\RG(X,V)$ at $(V,\cC)$. Set
$$
\Bun_{G,H}=\Bun_H\times_{\Pic X} \Bun_G, 
$$
where  the map $\Bun_H\to \Pic X$ sends $(V,\cC, \Sym^2 V\to\cC)$ to $\Omega\otimes\cC^{-1}$. The map $\Bun_G\to\Pic X$ sends $(M, \wedge^2 M\to\cA)$ to $\cA$. We have an isomorphism
$\cC\otimes \cA\,\iso\, \Omega$ for a point of $\Bun_{G,H}$. 
Let
$$
\tau: \Bun_{G,H}\to \Bun_{G_{2nm}}
$$ 
be the map sending a point as above to $V\otimes M$ with symplectic form $\wedge^2(V\otimes M)\to\Omega$. 

\begin{Pp} 
\label{Pp_line_bundle_description}
For a point of $\Bun_{G,H}$ as above we have a canonical $\ZZ/2\ZZ$-graded isomorphism
\begin{equation}
\label{iso_line_bundle_A}
\det\RG(X, V\otimes M)\,\iso\, \frac{\det\RG(X, V)^{2n}\otimes \det\RG(X, M)^{2m}}{\det\RG(X, \cO)^{2nm}\otimes\det\RG(X,\cA)^{2nm}}
\end{equation}
More precisely, we have a canonical $\ZZ/2\ZZ$-graded isomorphism of line bundles on $\Bun_{G,H}$
$$
\tau^*\cA_{G_{2nm}}\,\iso\, \cA_H^{2n}\otimes \cA_G^{2m}\otimes \cA_1^{-2nm}\otimes \det\RG(X,\cO)^{-2nm}
$$
\end{Pp}

\begin{Lm} 
\label{Lm_general_detRG}
i) For any $M\in\Bun_n$, $V\in\Bun_m$ there is a canonical $\ZZ/2\ZZ$-graded isomorphism
$$
\det\RG(X, M\otimes V)\,\iso\, \frac{\det\RG(X, M)^{\otimes m}\otimes\det\RG(X,V)^{\otimes n}}{\det\RG(X,\cA)\otimes\det\RG(X,\cB)} \otimes \frac{\det\RG(X, \cA\otimes\cB)}{\det\RG(X, \cO)^{\otimes nm-1}},
$$
where $\cA=\det M$, $\cB=\det V$. \\
ii) For any $\cA,\cB\in\Pic X$ there is a canonical isomorphism of $\ZZ/2\ZZ$-graded lines
$$
\frac{\det\RG(X, \cA\otimes\cB^m)\otimes \det\RG(X, \cA)^{m-1}}{\det\RG(X, \cA\otimes\cB)^m}\,\iso\, 
\frac{\det\RG(X, \cB^m)\otimes \det\RG(X, \cO)^{m-1}}{\det\RG(X, \cB)^m}
$$
\end{Lm}
\begin{Prf} 
i) Denote by $A(M,V)$ the $\ZZ/2\ZZ$-graded vector space
$$
\frac{\det\RG(X, M\otimes V)}{\det\RG(X, M)^{\otimes m}\otimes\det\RG(X,V)^{\otimes n}} \otimes
\frac{\det\RG(X,\cA)\otimes\det\RG(X,\cB)}{\det\RG(X, \cA\otimes\cB)},
$$
where $\cA=\det M$, $\cB=\det V$. 
View $A$ as a $\ZZ/2\ZZ$-graded line bundle on $\Bun_n\times\Bun_m$. Let us show that this line bundle is constant. 
 
  For an exact sequence of $\cO_X$-modules $0\to M\to M'\to M'/M\to 0$, where $M'/M$ is a torsion sheaf of length one at $x\in X$, we have a canonical $\ZZ/2\ZZ$-graded isomorphism $A(M',V)\,\iso\, A(M,V)$. 
Similarly, for an exact sequence of $\cO_X$-modules $0\to V\to V'\to V'/V\to 0$, where $V'/V$ is a torsion sheaf of length one at $x\in X$, we have
a canonical $\ZZ/2\ZZ$-graded isomorphism $A(M,V')\,\iso\, A(M,V)$.
 
  To conclude, note that $A(\cO^n, \cO^m)=\det\RG(X,\cO)^{\otimes 1-nm}$.
  
\smallskip\noindent
ii) The proof is similar. 
\end{Prf}

\medskip

\begin{Prf}\select{of Proposition~\ref{Pp_line_bundle_description}}

\smallskip\noindent
For a point $(M,\cA)$ of $\Bun_G$ the form $\wedge^2 M\to\cA$ induces an isomorphism $\det M\,\iso\, \cA^n$. So, by Lemma~\ref{Lm_general_detRG}, for a point of $\Bun_{G,H}$ as above we get a $\ZZ/2\ZZ$-graded isomorphism
\begin{equation}
\label{iso_for_proof_Pp2}
\det\RG(X, V\otimes M)\,\iso\, \frac{\det\RG(X, V)^{2n}\otimes\det\RG(X, M)^{2m}}{\det\RG(X, \cA^n)\otimes\det\RG(X,\det V)} \otimes
\frac{\det\RG(X, \cA^n\otimes\det V)}{\det\RG(X,\cO)^{4nm-1}}
\end{equation}
Applying it to $M=\oplus_{i=1}^n (\cO\oplus\cA)$ with the natural symplectic form $\wedge^2 M\to\cA$, we get
$$
\frac{\det\RG(X, \cO)^{2nm-1}}{\det\RG(X, \cA)^{2nm}}\,\iso\,
\frac{\det\RG(X, \cA^n\otimes\det V)}{\det\RG(X, \cA^n)\otimes\det\RG(X, \det V)}
$$
Combining the latter formula with (\ref{iso_for_proof_Pp2}), one concludes the proof.
\end{Prf}

\medskip
\noindent
3.2.1  By Proposition~\ref{Pp_line_bundle_description}, we get a map 
$
\tilde\tau: \Bun_{G,H}\to \Bunt_{G_{2nm}}
$ sending $(\wedge^2 M\to\cA, \Sym^2 V\to\cC, \cA\otimes\cC\,\iso\,\Omega)$ to $(\wedge^2(M\otimes V)\to\Omega, \cB)$. Here 
$$
\cB=\frac{\det\RG(X,V)^n\otimes\det\RG(X, M)^m}{\det\RG(X,\cO)^{nm}\otimes \det\RG(X,\cA)^{nm}},
$$
and $\cB^2$ is identified with $\det\RG(X, M\otimes V)$ via (\ref{iso_line_bundle_A}). 

\begin{Def} For the diagram of projections
$$
\Bun_H\getsup{\gq}  \Bun_{G,H}\toup{\gp} \Bun_G
$$
define $F_G: \D(\Bun_H)\to\D(\Bun_G)$ by
$$
F_G(K)=\gp_!(\tilde\tau^*\Aut\otimes \gq^*K)[\dimrel],
$$
where $\dimrel=\dim\Bun_{G_n}-\dim\Bun_{G_{2nm}}$. Define $F_H: \D(\Bun_G)\to\D(\Bun_H)$ by 
$$
F_H(K)=\gq_!(\tilde\tau^*\Aut\otimes \gp^*K)[\dimrel],
$$
where $\dimrel=\dim\Bun_{\SO_{2m}}-\dim\Bun_{G_{2nm}}$. Set also $F_{\tilde H}=\rho_H^*\comp F_H$. 
Replacing $\Aut$ by $\Aut_s$ (resp., by $\Aut_g$) in the above definitions, one defines the functors $F_{G,s}$, $F_{H,s}$, $F_{\tilde H,s}$ (resp., $F_{G,g}$, $F_{H,g}$, $F_{\tilde H,g}$). We write $F^G_H=F_H$ when we need to express the dependence of $F_H$ on $G$, and  similarly for $F^H_G=F_G$.
\end{Def}
 
 Let $\Bun_{G,\tilde H}$ be obtained from $\Bun_{G,H}$ by the base change $\Bun_{\tilde H}\to\Bun_H$. By abuse of notation, the restriction of $\tilde\tau:\Bun_{G,H}\to\Bunt_{G_{2nm}}$ to $\Bun_{G,\tilde H}$ is also denoted by  $\tilde\tau$. 
 
\medskip\noindent
3.2.2 Let $\Lambda_{\HH}$ (resp., $\check{\Lambda}_{\HH}$) denote the coweight (resp., weight) lattice for $\HH$. Write $\Lambda^+_{\HH}$ for the dominant coweights. The corresponding objects for $\GG$ are denoted $\Lambda_{\GG}, \check{\Lambda}_{\GG}$ and so on. 

For $m\ge 2$ let $\iota_m\in \Spin_{2m}$ be the central element of order 2 such that $\Spin_{2m}/\{\pm \iota_m\}\,\iso\,\SO_{2m}$. Here $\Spin_{2m}$ and $\SO_{2m}$ denote the corresponding split groups over $\Spec k$. For $m\ge 2$ denote by $\GSpin_{2m}$ the quotient of $\Gm\times\Spin_{2m}$ by the subgroup generated by $(-1,\iota)$. Let us convent that $\GSpin_2\,\iso\, \Gm\times\Gm$. The Langlands dual group is $\check{\HH}\,\iso\,\GSpin_{2m}$. 
We also have $\check{\GG}\,\iso\,\GSpin_{2n+1}$, where $\GSpin_{2n+1}$ is the quotient of $\Gm\times\Spin_{2n+1}$ by the diagonally embedded $\{\pm 1\}$.   
  
  Let $V_{\HH}$ (resp., $V_{\GG}$) denote the standard representation of $\SO_{2m}$ (resp., of $\SO_{2n+1}$). 
  
\smallskip\noindent
{\scshape CASE $m\le n$.} Pick an inclusion $V_{\HH}\hook{} V_{\GG}$ compatible with symmetric forms. It yields an inclusion $\check{\HH}\hook{}\check{\GG}$, which we assume compatible with the corresponding maximal tori. Pick an element $\sigma_{\GG}\in \SO(V_{\GG})\,\iso\, \check{\GG}_{ad}$ normalizing $\check{\TT}_{\GG}$ and preserving 
$V_{\HH}$ and $\check{\TT}_{\HH}\subset \check{\BB}_{\HH}$. Let $\sigma_{\HH}\in\OO(V_{\HH})$ be its restriction to $V_{\HH}$. We assume that $\sigma_{\HH}$ viewed as an automorphism of $(\check{\HH}, \check{\TT}_{\HH})$ extends the action of $\Sigma$ on the roots datum of $(\check{\HH},\check{\TT}_{\HH})$ defined in Section~3.1. 

 In concrete terms, one may take $V_{\GG}=k^{2n+1}$ with symmetric form given by the matrix
$$
\left(
\begin{array}{ccc}
0 & E_n & 0\\
E_n & 0 & 0\\
0 & 0 & 1
\end{array}
\right),
$$
where $E_n\in\GL_n$ is the unity. Take $\check{\TT}_{\GG}$ to be the maximal torus of diagonal matrices. Let $V_{\HH}\subset V_{\GG}$ be generated by $\{e_1,\ldots, e_m, e_{n+1},\ldots,e_{n+m}\}$. Let $\check{\TT}_{\HH}$ be the torus of diagonal matrices, and $\check{\BB}_{\HH}$ the Borel subgroup preserving for $i=1,\ldots,m$ the isotropic subspace generated by $\{e_1,\ldots,e_i\}$. Then one may take $\sigma_{\GG}$ permuting $e_m$ and $e_{n+m}$, sending $e_{2n+1}$ to $-e_{2n+1}$ and acting trivially on the other base vectors.  

 We let $\Sigma$ act on $\check{\HH}$ and $\check{\GG}$ via the elements $\sigma_{\HH}$, $\sigma_{\GG}$. So, the inclusion $\check{\HH}\hook{}\check{\GG}$ is $\Sigma$-equivariant and yields a morphism of the $L$-groups $\tilde H^L\to G^L$, where $\tilde H^L\,\iso\, \check{\HH}\ltimes\Sigma$ and $G^L\,\iso\,\check{\GG}\ltimes\Sigma$ (in the sense of B.2). 
 
\smallskip\noindent
{\scshape CASE $m>n$.} Pick an inclusion $V_{\GG}\hook{} V_{\HH}$ compatible with symmetric forms. It yields an inclusion $\check{\GG}\hook{}\check{\HH}$, which we assume compatible with the corresponding maximal tori. Let $\sigma_{\GG}$ be the identical automorphism of $V_{\GG}$. Extend it to an element $\sigma_{\HH}\in\OO(V_{\HH})$ by requiring that  $\sigma_{\HH}$ preserves $\check{\TT}_{\HH}\subset \check{\BB}_{\HH}$ and $\sigma_{\HH}\notin \SO(V_{\HH})$, $\sigma_{\HH}^2=\id$. 

  Let $\Sigma$ act on $\check{\HH}$ and $\check{\GG}$ via the elements $\sigma_{\HH}$, $\sigma_{\GG}$. The $\Sigma$-action on $(\check{\HH}, \check{\TT}_{\HH})$ extends the $\Sigma$-action (defined in Section~3.1) on the root datum of $(\check{\HH}, \check{\TT}_{\HH})$. Again, we get a morphism of the $L$-groups $G^L\to \tilde H^L$. Note that $\check{\GG}\times \Sigma\,\iso\,G^L$ is the direct product in this case. 
  
\medskip  
  
  As in B.2, in both cases the corresponding functoriality problem can be posed. As in B.1, for $\lambda\in\Lambda^+_{\HH}$ (resp., $\lambda\in\Lambda^+_G$) one defines the Hecke functors 
$$
\H^{\lambda}_{\tilde H}: \D(\Bun_{\tilde H})\to \D(\tilde X\times\Bun_{\tilde H})$$
and 
$$
\H^{\lambda}_G: \D(\Bun_G)\to \D(X\times\Bun_G)
$$
Write $V^{\lambda}_{\HH}$ (resp., $V^{\lambda}_{\GG}$) for the irreducible representation of $\check{\HH}$ (resp., of $\check{\GG}$) with highest weight $\lambda$. 

\begin{Con} 
\label{Con_functoriality}
i) Case $m=n$. For $\lambda\in\Lambda^+_{\GG}$ there is an isomorphism functorial in $K\in\D(\Bun_{\tilde H})$ 
$$
(\pi\times\id)^*\H^{\lambda}_G F_G(K)\,\iso\, \oplus_{\mu\in\Lambda^+_{\HH}}\,  (\id\boxtimes F_G)\H^{\mu}_{\tilde H}(K)\otimes \Hom_{\check{\HH}}(V^{\mu}_{\HH}, (V^{\lambda}_{\GG})^*)
$$
Here $\pi\times\id: \tilde X\times\Bun_G\to X\times\Bun_G$ and 
$\id\boxtimes F_G: \D(\tilde X\times \Bun_{\tilde H})\to \D(\tilde X\times \Bun_G)$ is the corresponding theta-lifting functor. 

\smallskip\noindent
ii) Case $m=n+1$. For $\mu\in\Lambda^+_{\HH}$ there is an isomorphism functorial in $K\in\D(\Bun_G)$
$$
\H^{\mu}_{\tilde H}F_{\tilde H}(K)\,\iso\, 
\oplus_{\lambda\in\Lambda^+_{\GG}}\, (\pi\times\id)^*
(\id\boxtimes F_{\tilde H})\H^{\lambda}_G(K)\otimes\Hom_{\check{\GG}}((V^{\lambda}_{\GG})^*, V^{\mu}_{\HH})
$$
Here $\pi\times\id: \tilde X\times\Bun_{\tilde H}\to X\times\Bun_{\tilde H}$ and $\id\boxtimes F_{\tilde H}: \D(X\times\Bun_G)\to \D(X\times\Bun_{\tilde H})$ is the corresponding theta-lifting functor. 

 In both cases these isomorphisms are compatible with the action of $\Sigma$ on both sides ($\Sigma$ acts on the Hecke operators for $\tilde H$ via (\ref{iso_action_Sigma_Hecke_sectionD})).
\end{Con}

\begin{Rem} For other pairs $(n,m)$ the relation between the theta-lifting functors and Hecke functors $F_G, F_H$ is expected to be essentially as in \cite{Ly4} involving the $\SL_2$ of Arthur.
\end{Rem}

\medskip\noindent 
3.2.3 Let $\act: \Pic X\times\Bun_G\to\Bun_G$ be the map sending $(\cL\in\Pic X, M,\cA)$ to $(M\otimes\cL, \cA\otimes\cL^2)$. Write also $\act:\Pic X\times\Bun_H\to\Bun_H$ for the map sending $(\cL\in\Pic X, V,\cC)$ to $(V\otimes\cL, \cC\otimes\cL^2)$. 

\begin{Def} 
\label{Def_aut_loc_sys}
For a local system $\cU$ on $X$ write $A\cU$ for the automorphic local system on $\Pic X$ corresponding to $\cU$. It is equipped with an isomorphism between the restriction of $A\cU$ under $X^{(d)}\to\Pic X, D\mapsto\cO(D)$ and $\cU^{(d)}$, this defines $A\cU$ up to a unique isomorphism.
\end{Def}
 
\begin{Def}  For a rank one local system $\cA$ on $X$ say that $K\in\D(\Bun_G)$ (resp., $K\in\D(\Bun_H)$)
has central character $\cA$ if $K$ is equipped with a $(\Pic X, A\cA)$-equivariant structure as in (as in \cite{Ly}, A.1, Definition~7). In particular, we have $\act^*K\,\iso\, A\cA\boxtimes K$. 
\end{Def}
   
\begin{Rem} 
\label{Rem_one_non_2_commutativity}
Using Lemma~\ref{Lm_general_detRG}, one checks that for $K\in\D(\Bun_{\tilde H})$ (resp., $K\in\D(\Bun_G)$) with central character $\chi$ the central character of $F_G(K)$ (resp., of $F_{\tilde H}(K)$) is $\chi^{-1}\otimes \cE_0^{\otimes n}$. 
The reason for that is the following. 
Let $\cX_{G,H}$ be the stack classifying $(M,\cA)\in\Bun_G$, $(V,\cC)\in\Bun_H$, $U\in\Pic X$ equipped with $\cA\otimes\cC\otimes U^2\,\iso\,\Omega$.  
We have a commutative diagram
$$
\begin{array}{ccc}
\Bun_{G,H} & \toup{\tau} & \Bun_{G_{2nm}}\\
\uparrow\lefteqn{\scriptstyle m_G} && \uparrow\lefteqn{\scriptstyle \tau}\\
\cX_{G,H} & \toup{m_H} & \Bun_{G,H},
\end{array}
$$
where $m_G$ (resp., $m_H$) sends the above collection to $(M\otimes U, \cA\otimes U^2)\in\Bun_G$, $(V,C)\in\Bun_H$ (resp., to $(M,\cA)\in\Bun_G$, $(V\otimes U, \cC\otimes U^2)\in\Bun_H$). Then the diagram  
$$
\begin{array}{ccc}
\Bun_{G,H} & \toup{\tilde\tau} & \Bunt_{G_{2nm}}\\
\uparrow\lefteqn{\scriptstyle m_G} && \uparrow\lefteqn{\scriptstyle \tilde\tau}\\
\cX_{G,H} & \toup{m_H} & \Bun_{G,H},
\end{array}
$$
is not 2-commutative in general. The key observation is as follows. Consider  a line bundle on $\Pic X$ whose fibre at $U\in\Pic X$ is
$$
\vartheta(U):=\frac{\det\RG(X, U\otimes\cE)\otimes\det\RG(X,\cO)}{\det\RG(X, U)\otimes\det\RG(X,\cE)}
$$
The tensor square of this line bundle is canonically trivialized, so defines a 2-sheeted covering of $\Pic X$. The corresponding local system of order 2 on $\Pic X$ is  $A\cE_0$.
\end{Rem} 
 
\medskip\noindent
3.3.1 Let $P\subset G$ be the Siegel parabolic, so $\Bun_P$ classifies: $L\in\Bun_n$, $\cA\in\Bun_1$, and an exact sequence of $\cO_X$-modules $0\to\Sym^2 L\to ?\to\cA\to 0$. Write $\nu_P:\Bun_P\to\Bun_G$ for the projection. Let $M_P$ be the Levi factor of $P$, so $\Bun_{M_P}\,\iso\,\Bun_n\times\Pic X$. 

 Set $\Bun_{P,H}=\Bun_P\times_{\Bun_G} \Bun_{G,H}$ and $\Bun_{P,\tilde H}=\Bun_P\times_{\Bun_G}\Bun_{G,\tilde H}$. We have a commutative diagram
\begin{equation}
\label{diag_tilde_Bun_PH_commutes}
\begin{array}{ccc}
\Bun_{G,H} & \toup{\tau} & \Bun_{G_{2nm}}\\ 
\uparrow && \uparrow\lefteqn{\scriptstyle \nu_{2nm}}\\
\Bun_{P,H}& \toup{\tau_P} & \Bun_{P_{2nm}},
\end{array}
\end{equation}
where $\tau_P$ sends $(V, \cC,\Sym^2 V\to\cC)$ and 
\begin{equation}
\label{ex_seq_P}
0\to \Sym^2 L\to ?\to\cA\to 0
\end{equation}
to the extension
\begin{equation}
\label{ex_seq_P_2nm}
0\to\Sym^2(V\otimes L)\to ?\to \Omega\to 0,
\end{equation}
which is the push-forward of (\ref{ex_seq_P})
under the composition
$$
\cA^{-1}\otimes\Sym^2 L\,\iso\, \cC\otimes\Omega^{-1}\otimes \Sym^2 L\to\Omega^{-1}\otimes\Sym^2 V\otimes\Sym^2 L\to \Omega^{-1}\otimes\Sym^2(V\otimes L)
$$  
Here the second map is induced by the form $\cC\to\Sym^2 V$.  We have used the fact that for $(V,\cC, \Sym^2 V\toup{h}\cC)$ of $\Bun_H$ the map $\Sym^2 V^*\,\iso\, (\Sym^2 V)\otimes \cC^{-2}\toup{h} \cC^{-1}$ induces a section $\cC\hook{} \Sym^2 V$ of $\Sym^2 V\toup{h} \cC$,  so $\cC$ is naturally a direct summand of $\Sym^2 V$. 

 Though (\ref{diag_tilde_Bun_PH_commutes}) commutes, the following diagram is not 2-commutative
$$
\begin{array}{ccc}
\Bun_{G,\tilde H} & \toup{\tilde\tau} & \Bunt_{G_{2nm}}\\ 
\uparrow\lefteqn{\scriptstyle \nu_P\times\id} && \uparrow\lefteqn{\scriptstyle \tilde\nu_{2nm}}\\
\Bun_{P,\tilde H}& \toup{\tau_P} & \Bun_{P_{2nm}},
\end{array}
$$
its non commutativity is measured by the following lemma. Write $\alpha_P:\Bun_P\to\Pic X$ for the map sending (\ref{ex_seq_P}) to $\det L$. 
Let $\alpha_{P,\tilde H}$ be the composition $\Bun_{P,\tilde H}\to\Bun_P\toup{\alpha_P}\Pic X$. 

\begin{Lm} 
\label{Lm_non_2_commutativity}
There is a canonical isomorphism over $\Bun_{P,\tilde H}$
$$
(\nu_P\times\id)^*\tilde\tau^*\Aut\,\iso\, \tau_P^*\tilde\nu_{2nm}^*\Aut\otimes \alpha_{P,\tilde H}^* A\cE_0
$$
\end{Lm}
\begin{Prf}
Write $\Bunt_{P,\tilde H}$ for the restriction of the gerbe $\Bunt_{G_{2nm}}\to\Bun_{G_{2nm}}$ under $\tau\comp(\nu_P\times\id): \Bun_{P,\tilde H}\to\Bun_{G_{2nm}}$. The map $\tilde\tau\comp(\nu_P\times\id)$ yields a trivialization $\Bunt_{P, \tilde H}\,\iso\, \Bun_{P,\tilde H}\times B(\mu_2)$ of this gerbe. So, $\tilde\nu_{2nm}\comp \tau_P$ gives rise to a map $\Bun_{P,\tilde H}\to B(\mu_2)$. The corresponding $\mu_2$-torsor over $\Bun_{P,\tilde H}$ is calculated using Lemma~\ref{Lm_general_detRG}. Namely, we have a line bundle on $\Bun_{P,\tilde H}$ whose fibre at $((\ref{ex_seq_P}), V,\cC, \Sym^2 V\to\cC,\gamma)$ is 
$$
\frac{\vartheta(\cC^m\otimes \det L)}{\vartheta(\cC^m)}
$$
The tensor square of this line bundle is canonically trivialized and gives rise to the corresponding $\mu_2$-torsor on $\Bun_{P,\tilde H}$. Our assertion follows by Remark~\ref{Rem_one_non_2_commutativity}. 
\end{Prf}
 
\begin{Def} Let $F_{P,\psi}: \D(\Bun_H)\to \D(\Bun_P)$ be the functor given by 
$$
F_{P,\psi}(K)=(\gp_P)_!(\gq_P^*K\otimes\tau_P^* S_{P,\psi})[\dimrel]
$$ 
for the diagram of projections
$$
\Bun_H\getsup{\gq_P}  \Bun_{P,H}\toup{\gp_P} \Bun_P,
$$
where $\dimrel=\dim\Bun_{P,H}-\dim\Bun_H-\dim\Bun_{P_{2nm}}$ is a function of a connected component of $\Bun_{P,H}$. Replacing $H$ by $\tilde H$ is the above diagram, one defines $F_{P,\psi}:\D(\Bun_{\tilde H})\to\D(\Bun_P)$ by the same formula. 
\end{Def}   

\begin{Cor} 
\label{Cor_relation_F_P_and_F_G}
The isomorphism (\ref{iso_gr_k}) yields 
an isomorphism $F_{P,\psi}\otimes \alpha_P^* A\cE_0\,\iso\,\nu_P^*F_G[\dimrel(\nu_P)]$ of functors  from $\D(\Bun_{\tilde H})$ to $\D(\Bun_P)$.
\QED
\end{Cor} 

\medskip
\noindent
3.3.2 Let $\cS_P$ be the stack classifying $L\in\Bun_n, \cA\in\Pic X$ and a section $\Sym^2 L\toup{s} \cA\otimes\Omega$. Then $\cS_P$ and $\Bun_P$ are dual (generalized) vector bundles over $\Bun_{M_P}\,\iso\,\Bun_n\times\Pic X$. Let $i_P:\Bun_{M_P}\hook{} \cS_P$ denote the zero section. 
 
  Let 
$$
\cV_{H,P}\to \Bun_n\times\Bun_H
$$ 
be the stack whose fibre over $(L\in\Bun_n, V,\cC, \Sym^2 V\to\cC)$ is $\Hom(V\otimes L, \Omega)$. We have a map $\gp_{\cV}: \cV_{H,P}\to \cS_P$ sending $(V,\cC)\in \Bun_H$, $L\in\Bun_n$, $L\toup{t} V^*\otimes\Omega$ to $(L,\cA,s)$, where $\cA=\Omega\otimes\cC^{-1}$ and $s$ is the composition
$$
\Sym^2 L\toup{t\otimes t} \Omega^2\otimes \Sym^2 V^*\to
\Omega^2\otimes\cC^{-1}\,\iso\, \cA\otimes\Omega
$$
\begin{Def} Let $F_{\cS}: \D(\Bun_H)\to\D(\cS_P)$ be given by
$F_{\cS}(K)=(\gp_{\cV})_! \gq_{\cV}^*K[\dimrel(\gq_{\cV})]$ 
for the diagram
\begin{equation}
\label{diag_cV}
\Bun_H \,\getsup{\gq_{\cV}} \,\cV_{H,P}\,\toup{\gp_{\cV}}\,\cS_P,
\end{equation}
where $\gq_{\cV}$ is the projection. 
\end{Def}

 The following is immediate from definitions.

\begin{Lm} There is a canonical isomorphism of functors $F_{P,\psi}\,\iso\, \Four_{\psi}\comp F_{\cS}$ from $\D(\Bun_H)$ to $\D(\Bun_P)$. \QED
\end{Lm}

\smallskip

 Let $\CT_P: \D(\Bun_G)\to\D(\Bun_{M_P})$ be the constant term functor given by $\CT_P(K)= \rho_P\nu_P^*K$ for the diagram $\Bun_G\getsup{\nu_P}\Bun_P\toup{\rho_P}\Bun_{M_P}$, where $\rho_P$ is the projection. Let $\alpha_{M_P}:\Bun_{M_P}\to\Pic X$ be the map sending $(L,\cA)$ to $\det L$. 
 
\begin{Cor} 
\label{Cor_CT_P}
The isomorphism (\ref{iso_gr_k}) induces 
an isomorphism $\CT_P\comp F_G\,\iso\, i_P^*F_{\cS}\otimes \alpha_{M_P}^* A\cE_0$ of functors (up to a shift) from $\D(\Bun_{\tilde H})$ to $\D(\Bun_{M_P})$. \QED
\end{Cor}

 Let $\cV_{\tilde H, P}$ be obtained from $\cV_{H,P}$ by the base change $\Bun_{\tilde H}\to\Bun_H$. Denote by
\begin{equation}
\label{map_xi_cV}
\xi_{\cV}: \cV_{H,P}\to \Bun_H\times_{\Pic X}\, \cS_P
\end{equation}
the map $(\gq_{\cV}, \gp_{\cV})$. The map $\cV_{\tilde H,P}\to \Bun_{\tilde H}\times_{\Pic X}\,\cS_P$ obtained from $\xi_{\cV}$ by the base change $\Bun_{\tilde H}\to\Bun_H$ is again denoted $\xi_{\cV}$ by abuse of notation. 

 Define the complex $^0\cT$ on $\Bun_{\tilde H}\times_{\Pic X}\, \cS_P$ by
$$
^0\cT=\xi_{\cV !}\Qlb[\dim\cV_{\tilde H,P}],
$$  
The group $S_2$ acts on $\xi_{\cV}$ changing the sign of $t:L\to V^*\otimes\Omega$, so $S_2$ acts also on $^0\cT$. Let $\Four_{\psi}: \D(\Bun_{\tilde H}\times_{\Pic X}\cS_P)\to \D(\Bun_{P,\tilde H})$ denote the Fourier transform. For the map $\tau_P:\Bun_{P,\tilde H}\to\Bun_{P_{2nm}}$ we have a $S_2$-equivariant isomorphism
$$
\tau_P^*S_{P,\psi}[\dimrel(\tau_P)]\,\iso\, \Four_{\psi}(^0\cT)
$$

\begin{Rem} If $G_1$ is a connected reductive group, which is not a torus, it is expected that for the projection $\pr_{G_1}: \Bun_{G_1}\to\Bun_{G_1/[G_1,G_1]}$ and a cuspidal complex $K\in\D(\Bun_{G_1})$ we have $(\pr_{G_1})_!K=0$. The reason to believe in this is that for $K$ cuspidal and automorphic the eigenvalues of Hecke operators acting on $K$ and on $\Qlb$ are different, so that $\Qlb$ and $K$ are `orthogonal'. 
 
 So, for $m>1$ it is expected that for the projection $\pr_{\tilde H}: \Bun_{\tilde H}\to \Pic X$ sending $(V,\cC)$ to $\cC$ and a cuspidal $K\in \D(\Bun_{\tilde H})$ we have $(\pr_{\tilde H})_!K=0$. If this is true then for such cupidal $K$ we have $\CT_P(F_G(K))=0$ (in particular, if $n=1$ then $F_G(K)$ is cuspidal). 
 
 Indeed, by Corolary~\ref{Cor_CT_P}, a fibre of $\CT_P(F_G(K))$ over $(L,\cA=\cC^{-1}\otimes\Omega)$ is an integral over $(V,\cC)\in\Bun_{\tilde H}$ equipped with a map $L\to V^*\otimes\Omega$, whose image is isotropic. 
But we can first fix the isotropic subbundle of $V^*\otimes\Omega$ generated by the image of $L$ and then integrate. The corresponding vanishing follows.
\end{Rem}

\medskip\noindent
3.4 {\scshape The case of split $H$.} In this subsection we assume the covering $\pi:\tilde X\to X$ split. Let $\tilde Q\subset \tilde H=\GO^0_{2m}$ be the Siegel parabolic, then $\Bun_{\tilde Q}$ is the stack classifying $L\in\Bun_m, \cC\in\Bun_1$ and an exact sequence $0\to\wedge^2 L\to ?\to\cC\to 0$ on $X$. The projection $\nu_{\tilde Q}:\Bun_{\tilde Q}\to\Bun_{\tilde H}$ sends the above point to $(V,\cC, \Sym^2V\to\cC,\gamma)$, where $V$ is included into an exact sequence $0\to L\to V\to L^*\otimes\cC\to 0$ and $\gamma: \det V\,\iso\,\cC^m$. 

 Let $M_{\tilde Q}$ be the Levi factor of $\tilde Q$, so $\Bun_{M_{\tilde Q}}\,\iso\, \Bun_m\times\Pic X$. Let $\Bun_{G,\tilde Q}$ be the stack obtained from $\Bun_{G,H}$ by the base change $\Bun_{\tilde Q}\to\Bun_H$.  Lemma~\ref{Lm_general_detRG} implies that the following diagram is 2-commutative
$$
\begin{array}{ccc}
\Bun_{G,H} & \toup{\tilde\tau} & \Bunt_{G_{2nm}}\\ 
\uparrow && \uparrow\lefteqn{\scriptstyle \tilde\nu_{2nm}}\\
\Bun_{G,\tilde Q}& \toup{\tau_{\tilde Q}} & \Bun_{P_{2nm}},
\end{array}
$$
where $\tau_{\tilde Q}$ sends $(M,\cA, \wedge^2 M\to\cA,\, 
0\to \wedge^2 L\to ?\to \cC\to 0)$ to the extension
$$
0\to \Sym^2(M\otimes L)\to ?\to\Omega\to 0,
$$
which is the push-forward of $0\to \cA\otimes\wedge^2 L\to ?\to \Omega\to 0$ under the composition $\cA\otimes\wedge^2 L\to \wedge^2 M\otimes\wedge^2 L\to \Sym^2(M\otimes L)$. Remind that we have a canonical direct sum decomposition $\Sym^2(M\otimes L)\,\iso\, (\Sym^2 M\otimes\Sym^2 L)\oplus (\wedge^2 M\otimes\wedge^2 L)$. 

\begin{Def} Let $F_{\tilde Q,\psi}:\D(\Bun_G)\to\D(\Bun_{\tilde Q})$ be the functor given by 
$$
F_{\tilde Q,\psi}(K)=\gp_{\tilde Q\, !}(\gq_{\tilde Q}^*K\otimes \tau_{\tilde Q}^*S_{P,\psi})[\dimrel]
$$
for the diagram of projections
$$
\Bun_G\getsup{\gq_{\tilde Q}}\Bun_{G,\tilde Q}\toup{\gp_{\tilde Q}}\Bun_{\tilde Q},
$$
where $\dimrel=\dim\Bun_{G,\tilde Q}-\dim\Bun_G-\dim\Bun_{P_{2nm}}$. 
\end{Def}

\begin{Cor} The isomorphism (\ref{iso_gr_k}) induces an isomorphism $F_{\tilde Q,\psi}\,\iso\, \nu_{\tilde Q}^*F_{\tilde H}[\dimrel(\nu_{\tilde Q})]$ of functors from $\D(\Bun_G)$ to $\D(\Bun_{\tilde Q})$.
\end{Cor}

 Let $\cS_{\tilde Q}$ be the stack classifying $L\in\Bun_m, \cC\in\Pic X$ and $\wedge^2 L\toup{s} \cC\otimes\Omega$. Then $\cS_{\tilde Q}$ and $\Bun_{\tilde Q}$ are dual (generalized) vector bundles over $\Bun_{M_{\tilde Q}}$. 

 Let $\cW_{G,\tilde Q}\to \Bun_m\times\Bun_G$ be the stack whose fibre over $L\in\Bun_m, (M,\cA)\in\Bun_G$ is $\Hom(M\otimes L,\Omega)$. 
We have a map $\gp_{\cW}: \cW_{G,\tilde Q}\to \cS_{\tilde Q}$ sending 
$L\in\Bun_m, (M,\cA)\in\Bun_G$ and $t:L\to M^*\otimes\Omega$ to $(L, \cC=\Omega\otimes\cA^{-1}, s)$, where $s$ is the composition
$$
\wedge^2 L\,\toup{\wedge^2 t}\,(\wedge^2 M^*)\otimes\Omega^2\to \cA^{-1}\otimes\Omega^2
$$
Define $F_{\cS_{\tilde Q}}:\D(\Bun_G)\to \D(\cS_{\tilde Q})$ by 
$$
F_{\cS_{\tilde Q}}(K)=(\gp_{\cW})_!\gq_{\cW}^*K[\dimrel(\gq_{\cW})]
$$ 
for the diagram
$$
\Bun_G\getsup{\gq_{\cW}} \cW_{G,\tilde Q}\toup{\gp_{\cW}}\cS_{\tilde Q}
$$
Here $\gq_{\cW}$ is the projection. 
 From definitions one gets
 
\begin{Lm} There is a canonical isomorphism of functors $F_{\tilde Q,\psi}\,\iso\,\Four_{\psi}\comp F_{\cS_{\tilde Q}}$ from $\D(\Bun_G)$ to 
$\D(\Bun_{\tilde Q})$. \QED
\end{Lm}

\medskip\noindent
3.5 {\scshape Weil representation and two-sheeted coverings} 

\smallskip\noindent
Write $\Bun_{G_n,\tilde X}$ for the stack classifying rank $2n$ vector bundles $W$ on $\tilde X$ with symplectic form $\wedge^2 W\to\Omega_{\tilde X}$. Let $\pi_n: \Bun_{G_n,\tilde X}\to\Bun_{G_{2n}}$ be the map sending the above point to $\pi_*W$ equipped with natural symplectic form $\wedge^2 (\pi_*W)\to\Omega$. 
 Denote by $\cA_{G_n,\tilde X}$ the line bundle on $\Bun_{G_n,\tilde X}$ with fibre $\det\RG(\tilde X, W)$ at $W$. Since $\pi_n^*\cA_{G_{2n}}\,\iso\, \cA_{G_n,\tilde X}$ canonically, $\pi_n$ lifts to a map
$$
\tilde\pi_n:\Bunt_{G_n,\tilde X}\to\Bunt_{G_{2n}}
$$ 
\begin{Pp} 
\label{Pp_Weil_and_coverings}
There is a canonical isomorphism  ${_{\tilde X}\Aut}\,\iso\, \tilde\pi_n^*\Aut[\dimrel(\tilde\pi_n)]$ preserving the generic and special parts.
\end{Pp}
\begin{Prf}
Let $_i\Bun_{G_k}\subset\Bun_{G_k}$ be the locally closed substack
given by $\dim\H^0(X, M)=i$ for $M\in\Bun_{G_k}$. Let $_i\Bunt_{G_k}$ be the restriction of the $\mu_2$-gerbe $\Bunt_{G_k}\to\Bun_{G_k}$ to $_i\Bun_{G_k}$. As in (\cite{Ly}, Remark~1), we have a cartesian square, where the vertical arrows are canonical sections of the corresponding $\mu_2$-gerbes
$$
\begin{array}{ccc}
{_i\Bun_{G_n,\tilde X}} & \toup{\pi_n} & {_i\Bun_{G_{2n}}} \\
\downarrow && \downarrow\\
{_i\Bunt_{G_n,\tilde X}} & \toup{\tilde\pi_n} &  {_i\Bunt_{G_{2n}}}
\end{array}
$$
This gives a canonical normalization of the sought-for isomorphism over $_i\Bunt_{G_n,\tilde X}$ for $i=0,1$. It remains to show its existence. 

 To do so, consider the commutative diagram
$$
\begin{array}{ccc}
\Bun_{P_n,\tilde X} & \toup{\tilde\nu_{n,\tilde X}} & \Bunt_{G_n,\tilde X}\\
\downarrow\lefteqn{\scriptstyle\pi_{n,P}} && \downarrow\lefteqn{\scriptstyle \tilde\pi_n}\\
\Bun_{P_{2n}} & \toup{\tilde\nu_{2n}} & \Bunt_{G_{2n}},
\end{array}
$$ 
where we denote by $\pi_{n,P}$ the following map. Given an exact sequence on $\tilde X$
\begin{equation}
\label{seq_point_of_Bun_Pn_tilde_X}
0\to \Sym^2 L\to ?\to \Omega_{\tilde X}\to 0,
\end{equation}
summate it with the one obtained by applying $\sigma^*$. The resulting exact sequence 
$$
0\to \Sym^2 (L\oplus \sigma^* L)\to ?\to \Omega_{\tilde X}\to 0
$$ 
is equipped with descent data for $\tilde X\toup{\pi} X$, so yields an exact sequence $0\to\Sym^2(\pi_*L)\to ?\to\Omega\to 0$, which is the image of (\ref{seq_point_of_Bun_Pn_tilde_X}) by $\pi_{n,P}$. 
 
 Let $_c\Bun_{n,\tilde X}\subset \Bun_{n,\tilde X}$ be the open substack classifying $W\in\Bun_{n,\tilde X}$ with $\H^0(\tilde X, W)=0$. Let $_c\cV_{n,\tilde X}\to {_c\Bun_{n,\tilde X}}$ be the vector bundle with fibre $\Hom(W, \Omega_{\tilde X})$ at $W$. Write $_c\Bun_{P_n,\tilde X}$ for the preimage of $_c\Bun_{n,\tilde X}$ under the projection $\Bun_{P_n,\tilde X}\to\Bun_{n,\tilde X}$. We get a commutative diagram
$$
\begin{array}{ccc}
\Sym^2 {_c\cV^*_{n,\tilde X}} & \gets & _c\Bun_{P_n,\tilde X}\\ 
\downarrow\lefteqn{\scriptstyle\pi_\cV} && \downarrow\lefteqn{\scriptstyle\pi_{n,P}}\\
\Sym^2 {_c\cV^*_{2n}} & \gets & _c\Bun_{P_{2n}}
\end{array}
$$ 
where the horizontal arrows are those of (\cite{Ly}, 5.2). Here $\pi_{\cV}$ is the map sending $L\in\Bun_{n,\tilde X}, b\in\Sym^2 H^1(\tilde X, L)$ to $\pi_*L\in\Bun_{2n}, b\in\Sym^2\H^1(X, \pi_*L)$.

 By (\select{loc.cit}, Proposition~7), it suffices to show that over $_c\Bun_{P_n,\tilde X}$ there is an isomorphism 
$$
\pi_{n,P}^*S_{P,\psi}[\dimrel(\pi_{n,P})]\,\iso\, S_{P,\psi}
$$ 
We have the sheaves $S_{\psi}$ on $\Sym^2 {_c\cV^*_{n,\tilde X}}$ and on $\Sym^2 {_c\cV^*_{2n}}$ defined in (\select{loc.cit.}, 4.3). Since $\pi_{\cV}^*S_{\psi}\,\iso\, S_{\psi}$ up to a shift, our assertion follows from (\select{loc.cit.}, 5.2).
\end{Prf}

\bigskip\noindent
3.6 {\scshape Whittaker type functors} 

\medskip\noindent
3.6.1 Write $\Bunt_P$ for the Drinfeld compactification of $\Bun_P$ introduced in (\cite{BG}, 1.3). So, $\Bunt_P$ classifies $(M,\cA)\in\Bun_G$ together with a Lagrangian subsheaf $L\subset M$, $L\in\Bun_n$. Then $\Bun_P\subset\Bunt_P$ is the open substack given by the condition that $L$ is a subbundle of $M$. 

  In the spirit of (\cite{Ly2}, Sect.~7), let $\cZ_1$ denote the stack obtained from $\Bun_{G,\tilde H}$ by the base change $\Bunt_P\to\Bun_G$. Let $\nu_{\cZ}: \cZ_1\to\Bun_{G,\tilde H}$ be the projection. 
  
  
  Denote by $\pi_{2,1}: \cZ_2\to\cZ_1$ the stack over $\cZ_1$ with fibre consisting of all maps $s:\Sym^2 L\to \cA\otimes\Omega$. A version of (Theorem~3, \select{loc.cit.}) holds. Namely, one defines a Whittaker category $\D^W(\cZ_2)$ as in (\select{loc.cit.}, 2.10). Here is its description on strata. 
  
   For $d\ge 0$ let $^{\le d}\cZ_1\subset\cZ_1$ be the closed substack given by the condition that $\wedge^n L\hook{} \wedge^n M$ has zeros of order $\le d$. Its open substack $^d\cZ_1\subset{^{\le d}\cZ_1}$ is given by: there is a subbundle $L'\subset M$ such that $L\subset L'$ is a subsheaf with $d=\deg(L'/L)$. Then $^{\le d}\cZ_1$ is stratified by $^i\cZ_1$ for $0\le i\le d$. 
   
   The stack $^d\cZ_1$ classifies collections: a modification of rank $n$ vector bundles $L\subset L'$ on $X$ with $\deg(L'/L)=d$, $\cA\in\Pic X$, and an exact sequence $0\to\Sym^2 L'\to ?\to \cA\to 0$ on $X$, $(V,\cC, \Sym^2 V\to\cC, \gamma)\in\Bun_{\tilde H}$ with $\cC\otimes\cA\,\iso\,\Omega$. 

 Set 
$$
^d\cZ_2=\cZ_2\times_{\cZ_1} {^d\cZ_1}\;\;\;\mbox{and}\;\;\; {^{\le d}\cZ_2}=\cZ_2\times_{\cZ_1} {^{\le d}\cZ_1}
$$ 
Let $^d\cZ'_2\hook{} {^d\cZ_2}$ be the closed substack given by the condition: $s$ factors as 
$$
\Sym^2 L\hook{} \Sym^2 L'\to \cA\otimes\Omega
$$ 
Let $^d\chi: {^d\cZ'_2}\to\A^1$ be the pairing of $s$ with the extension $0\to\Sym^2 L'\to ?\to\cA\to 0$. 

Let $^d\cP_2$ be the stack classifying: $(V,\cC, \Sym^2 V\to\cC, \gamma)\in\Bun_{\tilde H}$, a modification of rank $n$ vector bundles $L\subset L'$ on $X$ with $d=\deg(L'/L)$, and a section $s:\Sym^2 L'\to\cA\otimes\Omega$ with $\cC\otimes\cA\,\iso\,\Omega$. The projection
$$
\phi_2: {^d\cZ'_2}\to {^d\cP_2}
$$ 
is smooth. 

\begin{Lm} Any object of $\D^W(^d\cZ_2)$ is the extension by zero from $^d\cZ'_2$. The functor 
$$
^dJ(K)={^d\chi^*\cL_{\psi}}\otimes \phi^*_2 K[\dimrel]
$$ 
provides an equivalence of categories $^dJ: \D(^d\cP_2)\to \D^W(^d\cZ_2)$ and is exact for the perverse t-structures. Here $\dimrel$ is the relative dimension of $\phi_2$. 
\QED
\end{Lm}

\begin{Pp} 
\label{Pp_Whittaker_functors}
There is an equivalence of categories $W_{12}: \D(\cZ_1)\,\iso\, \D^W(\cZ_2)$, which is exact for the perverse t-structures, and $(\pi_{2,1})_!$ is quasi-inverse to it. Moreover, for any $K\in \D^W(\cZ_2)$ the natural map $(\pi_{2,1})_! K\to (\pi_{2,1})_* K$ is an isomorphism. \QED
\end{Pp}

\noindent
3.6.2  We have naturally $^0\cP_2\,\iso\, \Bun_{\tilde H}\times_{\Pic X}\cS_P$. Note that $^d\cZ_1$ is of codimension $d$ in $\cZ_1$. Define the complex $^d\cT\in \D(^d\cP_2)$ by
$$
^d\cT=\pr^*(^0\cT)[d+nd],
$$
where $\pr: {^d\cP_2}\to {^0\cP_2}$ is the projection forgetting $L$ (the   relative dimension of $\pr$ is $nd$). Set also $^dK_2={^dJ(^d\cT)}$. 
Let $S_2$ act on $^dK_2$ via its action on $^0\cT$. 
Using (\ref{iso_gr_k}), for the $*$-restriction we get 
$$
W_{12}(\nu_{\cZ}^*\tilde\tau^*\Aut)\mid_{^d\cZ_2}[\dimrel(\tilde\tau\comp\nu_{\cZ})]\,\iso\, {^dK_2}
$$ 

\noindent
{\bf Question.} Consider the object of the Grothendieck group of $\D(\cZ_2)$, whose $*$-restriction to $^d\cZ_2$ is $^dK_2$. Does it satisfy the Hecke property in the Grothendieck group of $\cZ_2$?


\bigskip\bigskip
\centerline{\scshape 4. The pair $\GL_2, \GO_2$}

\bigskip\noindent
4.1 Keep the notation of Sect.~3 assuming $n=1, m=1$. 
So, $G=\GG=\GL_2$. In Sect.~4 we assume $\tilde X$ connected.\footnote{Some of our results extend to the case of non connected $\tilde X$, but this case reduces to the study of renormalized geometric Eisenstein series from \cite{BG}.}
Identify $\HH=\GO_2^0$ with $\Gm\times\Gm$ in such a way that the automorphism $\sigma$ of $\HH$ permutes the two copies of $\Gm$, so $\tilde H=\pi_*\Gm$. 
We have canonically $\Lambda_{\HH}\,\iso\, \ZZ^2$, and $\sigma$ sends a coweight $\mu=(\mu_1,\mu_2)$ to $\sigma\mu=(\mu_2,\mu_1)$. 

We have a canonical isomorphism $\Pic\tilde X\,\iso\,\Bun_{\tilde H}$ sending $\cB\in\Pic\tilde X$ to $V=\pi_*\cB, \; \cC=N(\cB)$ equipped with natural symmetric form $\Sym^2 V\to\cC$ and isomorphism $\gamma: \cC\,\iso\,\cE\otimes\det V$. Here $N:\Pic\tilde X\to\Pic X$ is the norm map (cf. A.1). We also write $\sigma$ for the map $\Pic\tilde X\to\Pic\tilde X$ sending $\cB$ to $\sigma^*\cB$.
The following diagram is 2-commutative
$$
\begin{array}{ccc}
\Pic\tilde X & \toup{\sigma} & \Pic\tilde X\\
\downarrow\lefteqn{\scriptstyle\rho_H} & \swarrow\lefteqn{\scriptstyle\rho_H}\\ 
\Bun_H
\end{array}
$$

Write $X^{(d)}$ for the $d$-th symmetric power of $X$, we also view it as the scheme classifying effective divisors on $X$ of degree $d$. 

 Let $\Pic^d \tilde X$ be the connected component of $\Pic\tilde X$ classifying $\cL\in\Pic\tilde X$ with $\deg\cL=d$. Write $\Pic'^d \tilde X$ for the stack classifying $\cL\in\Pic^d\tilde X$ with a section $\cO\to\cL$. One defines $\Pic'^d X$ similarly. Let 
$$
\pi'_{ex}: \Pic'^d\tilde X\to \Pic'^d X\times_{\Pic X} \Pic\tilde X
$$ 
be the map sending $(\cL\in\Pic^d\tilde X, \cO\toup{t}\cL)$ to $(\cL, \cO\toup{N t} N(\cL))$. The group $S_2$ acts on $\pi'_{ex}$ sending $(\cL,t)$ to $(\cL, -t)$.
We have an open immersion $X^{(d)}\hook{} \Pic'^d X$ corresponding to nonzero sections. Restricting $\pi'_{ex}$ to the the corresponding open substacks, one gets a map  
$$
\pi': \tilde X^{(d)}\to X^{(d)}\times_{\Pic X} \Pic \tilde X
$$
 
 By abuse of notation, denote also by $\pi$ the direct image map $\pi: \tilde X^{(d)}\to X^{(d)}$. 

\begin{Lm}  
\label{Lm_evident}
i) For any local system $\tilde E$ on $\tilde X$ we have $\pi_* \tilde E^{(d)}\,\iso\, (\pi_*\tilde E)^{(d)}$. For a rank one local system $\chi$ on $X$ we have $\pi^*(\chi^{(d)})\,\iso\, (\pi^*\chi)^{(d)}$. \\
ii) For each $d\ge 0$ both the $S_2$-invariants and anti-invariants in $\pi'_!\Qlb[d]$ are irreducible perverse sheaves.
If $d>2\tilde g-2$ then the same holds for $(\pi'_{ex})_!\Qlb[d]$.
\end{Lm} 
\begin{Prf}
ii) The map $\pi'_{ex}$ is finite. Write $\uPic\tilde X$ for the Picard scheme of $\tilde X$ and similarly for $X$. We have a $\mu_2$-gerbe
$$
\und{\gr}: X^{(d)}\times_{\Pic X} \Pic \tilde X \to X^{(d)}\times_{\uPic X} \uPic \tilde X
$$
\Step 1  Let us show that the $S_2$-invariants in $\pi'_!\Qlb[d]$ is an irreducible perverse sheaf. It suffices to show that $\und{\gr}_!\pi'_!\Qlb[d]$ is an irreducible perverse sheaf.

 Let $Z$ denote the image of the (finite) map $\und{\gr}\comp \pi'$ (with reduced scheme structure). The projection $p_Z: Z\to X^{(d)}$ is a finite map. Take a rank one local system $\tilde E$ on $\tilde X$ that does not descend to $X$. Let $A\tilde E$ denote the corresponding automorphic local system on $\Pic\tilde X$ (cf. Definition~\ref{Def_aut_loc_sys}). Since $(p_Z)_!$ sends a nonzero perverse sheaf to a nonzero perverse sheaf, it suffices to show that 
$$
(p_Z)_!((\und{\gr}_!\pi'_!\Qlb[d])\otimes (\Qlb\boxtimes A\tilde E))
$$
is an irreducible perverse sheaf. But the latter identifies with $(\pi_!\tilde E)^{(d)}[d]$, so is irreducible.  

\smallskip
\Step 2 Let $\eta_Z$ (resp., $\eta$) denote the generic point of $Z$ (resp., of $\tilde X^{(d)}$). From Step 1 it follows that the map $\und{\gr}\comp \pi'$ yields an isomorphism $\eta\,\iso\,\eta_Z$. So, the restriction of the $\mu_2$-gerbe $\und{\gr}$ to $\eta_Z$ is trivial. For any map $\xi=(\id,\xi'): \eta\to \eta\times B(\mu_2)$, the $S_2$-anti-invariants in $\xi_!\Qlb$ is an irreducible local system. It follows that the $S_2$-anti-invariants in $\pi'_!\Qlb[d]$ is an irreducible perverse sheaf.

\smallskip
\Step 3  For $d>2\tilde g-2$ the stack $\Pic'^d\tilde X$ is smooth. Since $(\pi'_{ex})_!\Qlb[d]$ is the Goresky-MacPherson extension from $X^{(d)}\times_{\Pic X} \Pic \tilde X$, from Steps 1 and 2 we learn that both the $S_2$-invariants and anti-invariants in $(\pi'_{ex})_!\Qlb[d]$ are irreducible perverse sheaves.
\end{Prf}

\medskip

 For the map $\tilde\tau:\Bun_{G,\tilde H}\to\Bunt_{G_2}$ set 
$$
\Aut_{G,\tilde H,g}=\tilde\tau^*\Aut_g[\dimrel(\tilde\tau)]\;\;\;\; and\;\;\;\; \Aut_{G,\tilde H,s}=\tilde\tau^*\Aut_s[\dimrel(\tilde\tau)]
$$ 
and $\Aut_{G,\tilde H}=\Aut_{G,\tilde H,g}\oplus \Aut_{G,\tilde H,s}$. 

\begin{Pp} 
\label{Pp_irreducibility}
i) Both $\Aut_{G,\tilde H,g}$ and $\Aut_{G,\tilde H,s}$ are irreducible perverse sheaves, and we have $\DD(\Aut_{G,\tilde H})\,\iso\,\Aut_{G,\tilde H}$ canonically. 

\smallskip\noindent
ii) The sheaf $\Aut_{G,\tilde H}$ is ULA w.r.t. $\Bun_{G,\tilde H}\to\Bun_{\tilde H}$. 
\end{Pp}
\begin{Prf} i) Consider the map $\Bun_{P,\tilde H}\to\Bun_{G,\tilde H}$ obtained from $\Bun_P\to\Bun_G$ by base change $\Bun_{G,\tilde H}\to\Bun_G$. 
Let $^0\Bun_{P,\tilde H}$ be the open substack of $\Bun_{P,\tilde H}$ given by $2\deg L+\deg\cC<0$. We have a commutative diagram
$$
\begin{array}{ccc}
\Bun_{G,\tilde H} & \toup{\tau} & \Bun_{G_2}\\
\uparrow && \uparrow\lefteqn{\scriptstyle \nu_2}\\
^0\Bun_{P, \tilde H} & \toup{\tau_P} & ^0\Bun_{P_2},
\end{array}
$$
where the vertical arrows are smooth and surjective. By (\cite{Ly}, Proposition 7) and Lemma~\ref{Lm_non_2_commutativity}, it suffices to show that both $\tau_P^* S_{P,\psi,g}[\dimrel]$ and $\tau_P^* S_{P,\psi,s}[\dimrel]$ are irreducible perverse sheaves over each connected component of $^0\Bun_{P,\tilde H}$. 

 Remind the following notation introduced in 3.3.2. The stack
$\cS_P$ classifies $L\in\Pic X$, $\cA\in\Pic X$ and a map $L^{\otimes 2}\to\cA\otimes\Omega$. Let $^0\cS_P\subset \cS_P$ be the open substack given by $2\deg L-\deg\cA+\deg\Omega<0$. 
 
 The stack $\cV_{\tilde H,P}$ classifies $L\in\Pic X$, $\cB\in\Pic \tilde X$ and $\pi^*L\toup{t}\cB^{\star}$. (We may view $t$ as a map $L\otimes V\to\Omega$ for $V=\pi_*\cB$). The map $\xi_{\cV}: \cV_{\tilde H,P}\to \Bun_{\tilde H}\times_{\Pic X} \,\cS_P$ sends $(L,\cB,t)$ to $(L,\cB,s)$, where $s: L^2\to N(\cB^{\star})$ is the norm of $t$.  Set 
$$
^0\cV_{\tilde H,P}=\xi_{\cV}^{-1}(\Bun_{\tilde H}\times_{\Pic X} {^0\cS_P})
$$  
Note that $^0\cV_{\tilde H,P}$ is smooth (here the connectedness of $\tilde X$ is essentially used). 
  
  By definition, $\tau_P^* S_{P,\psi}[\dimrel]$ is the Fourier transform of $(\xi_{\cV})_!\Qlb[\dim \cV_{\tilde H, P}]$. From ii) of Lemma~\ref{Lm_evident} it follows that $(\xi_{\cV})_!\Qlb[\dim \cV_{\tilde H, P}]$ is a direct sum of two irreducible perverse sheaves over $\Bun_{\tilde H}\times_{\Pic X} {^0\cS_P}$. We are done. 
  
\smallskip\noindent
ii) 
We need the following general observation. If $f:Y\to S$ is a vector bundle over a smooth base $S$, and $K\in\D(Y)$ is ULA w.r.t. $f$ then $\Four_{\psi}(K)$ is ULA w.r.t. the projection $Y^{\check{}}\to S$.

 Apply this to the vector bundle $v: \Bun_{\tilde H}\times_{\Pic X} {^0\cS_P}\to {^0(\Bun_{\tilde H}\times_{\Pic X}\Bun_{M_P})}$, where the base classifies pairs $\cB\in\Pic\tilde X, L\in\Pic X$ with $2\deg L+\deg\cB<0$. 
 
 Since the projection $^0\cV_{\tilde H,P}\to {^0(\Bun_{\tilde H}\times_{\Pic X}\Bun_{M_P})}$ is smooth, $(\xi_{\cV})_!\Qlb$ is ULA w.r.t. $v$, so its Fourier transform is also ULA over $^0(\Bun_{\tilde H}\times_{\Pic X}\Bun_{M_P})$. 

 Since $^0\Bun_{P,\tilde H}\to \Bun_{G,\tilde H}$ is smooth and surjective, our assertion follows (the ULA property is local in the smooth topology on the source).
\end{Prf}

\medskip

 For our particular pair $(\tilde H,G)$ the projection $\Bun_{G,\tilde H}\to\Bun_G$ is proper (this phenomenon does not happen for $m>1$). So, Proposition~\ref{Pp_irreducibility} implies that $F_G$ commutes with the Verdier duality. 
 
\bigskip\noindent
4.2 {\scshape Hecke property}

\medskip\noindent
4.2.1.1 For a dominant coweight $\lambda$ of $G$ write $\ov{\cH}^{\lambda}_G$ for the Hecke stack classifying $x\in X, M, M'\in\Bun_G$ with an isomorphism $\beta: M\,\iso\, M'\mid_{X-x}$ such that $M'$ is in a position $\le\lambda$ w.r.t. $M$ in the sense of (\cite{BG}). We have a diagram
$$
X\times \Bun_G\;\getsup{\supp\times p_G}\;\ov{\cH}^{\lambda}_G\;\toup{p'_G}\;\Bun_G,
$$
where $p_G$ (resp., $p'_G$) sends $(x,M,M',\beta)$ to $M$ (resp., $M'$). 

  We fix an inclusion $\check{\HH}\hook{}\check{\GG}=\GL_2$ as the maximal torus of diagonal matrices. This yields isomorphisms $\Lambda_{\HH}\,\iso\, \Lambda_{\GG}\,\iso\,\ZZ$.  Given a coweight $\mu=(\mu_1,\mu_2)\in\Lambda_{\HH}$ such that $\lambda-\mu$ vanishes in $\pi_1(G)$ (that is, $\lambda_1+\lambda_2=\mu_1+\mu_2$), consider the diagram 
$$
\begin{array}{cccc}
\tilde X\times\Bun_{\tilde H,G} & \getsup{\supp\times p_{\tilde H,G}} & \ov{\cH}^{\mu, \lambda}_{\tilde H, G} & \toup{p'_{\tilde H, G}}\;\Bun_{\tilde H, G}\\
\downarrow && \downarrow\\
X\times\Bun_G & \getsup{\supp\times p_G} & \ov{\cH}^{\lambda}_G
\end{array} 
$$
where $\ov{\cH}^{\mu,\lambda}_{\tilde H,G}$ classifies collections: $(\cB, M)\in\Bun_{\tilde H,G}$, $\tilde x\in\tilde X$ for which we set $x=\pi(\tilde x)$, and $(x,M,M',\beta)\in\ov{\cH}^{\lambda}_G$. The map $p_{\tilde H,G}$ forgets $(\beta,M')$, so the left square is cartesian. The map $p'_{\tilde H,G}$ sends the above collection to $(\cB',M')\in\Bun_{\tilde H,G}$, where $\cB'=\cB(\mu_1\tilde x+\mu_2\sigma(\tilde x))$. 
 
  The Hecke functor $\H^{\mu,\lambda}: \D(\Bun_{\tilde H,G})\to \D(\tilde X\times \Bun_{\tilde H,G})$ is given by 
\begin{equation}
\label{Hecke_functor_GH}
\H^{\mu,\lambda}(K)=(\supp\times p_{\tilde H,G})_!(\IC_{\ov{\cH}^{\mu, \lambda}_{\tilde H, G}}\otimes (p'_{\tilde H,G})^*K)[-\dim\Bun_{\tilde H,G}]
\end{equation}
We have  
$$
\pr^*\IC_{\ov{\cH}^{\lambda}_G}[\dimrel]\,\iso\,\IC_{\ov{\cH}^{\mu, \lambda}_{\tilde H, G}},
$$ 
where $\dimrel=\dim\Bun_{G,\tilde H}-\dim\Bun_G$.
By B.1.3, $\H^{\mu,\lambda}$ commutes with the Verdier duality.
As in B.1.2, we have a canonical isomorphism 
\begin{equation}
\label{iso_Hecke_mu_lambda}
\delta_{\sigma}: (\sigma\times\id)^*\comp \H^{\mu,\lambda}\,\iso\, \H^{\sigma\mu, \lambda}
\end{equation}
(we used that $\sigma$ acts trivially on the dominant coweights of $G$). 

 Write $V^{\lambda}$ for the irreducible representation of $\check{\GG}$ with h.w. $\lambda$. For a $\check{\GG}$-representation $V$ and $\mu\in\Lambda_{\HH}$ denote by $V(\mu)$ the $\mu$-weight space of $\check{\HH}$ in $V$. 
 
\begin{Th} 
\label{Th_Hecke_property}
For $\lambda\in\Lambda^+_{\GG},\mu\in\Lambda_{\HH}$ such that $\lambda-\mu$ vanishes in $\pi_1(G)$ the sheaf $\H^{\mu,\lambda}(\Aut_{G,\tilde H})$ is perverse, and we have 
\begin{equation}
\label{iso_Th_1}
\H^{\mu,\lambda}(\Aut_{G,\tilde H})\,\iso\, \oplus_{\nu} \H^{\nu,0}(\Aut_{G,\tilde H})\otimes (V^{\lambda})^*(\nu-\mu)
\end{equation} 
In other words, the sum (without multiplicities) is over the coweights $\nu=(-a,a)$ such that $\lambda_1-\mu_1\ge a\ge \mu_2-\lambda_1$. If $\lambda_1-\lambda_2$ is even then this isomorphism preserves generic and special parts, otherwise it interchanges them. 
\end{Th} 
 
\begin{Rem} The isomorphism (\ref{iso_Th_1})  is compatible with the action of $\Sigma=\{1,\sigma\}$ on both sides, that is, the diagram commutes
$$
\begin{array}{ccc}
(\sigma\times\id)^*\H^{\mu,\lambda}(\Aut_{G,\tilde H}) & \iso & (\sigma\times\id)^* (\oplus_{\nu} \H^{\nu,0}(\Aut_{G,\tilde H})\otimes (V^{\lambda})^*(\nu-\mu))\\
\downarrow\lefteqn{\scriptstyle \delta_{\sigma}} && \downarrow\lefteqn{\scriptstyle \delta_{\sigma}}\\
\H^{\sigma\mu,\lambda}(\Aut_{G,\tilde H}) & \iso  & \oplus_{\nu} \H^{\sigma\nu,0}(\Aut_{G,\tilde H})\otimes (V^{\lambda})^*(\sigma\nu-\sigma\mu)
\end{array}
$$
\end{Rem}
 
\medskip\noindent 
4.2.1.2 Remind that $\Bun_{\tilde H,P}$ denotes the stack obtained from $\Bun_{\tilde H,G}$ by the base change $\Bun_P\toup{\nu_P}\Bun_G$. Remind the stack $^{\le d}\cZ_1$ (cf. Sect.~3.6.1).  
We have a commutative diagram
$$
\begin{array}{ccccc}
\tilde X\times \Bun_{\tilde H,P} & \getsup{\supp\times p_{\tilde H,P}} & \ov{\cH}^{\mu,\lambda}_{\tilde H,P} & \toup{p'_{\tilde H,P}} & ^{\le \lambda_1-\lambda_2}\cZ_1\\
\downarrow\lefteqn{\scriptstyle \id\times\nu_P} && \downarrow&& \downarrow\lefteqn{\scriptstyle \nu_{\cZ}}\\
\tilde X\times\Bun_{\tilde H,G} & \getsup{\supp\times p_{\tilde H,G}} & \ov{\cH}^{\mu, \lambda}_{\tilde H, G} & \toup{p'_{\tilde H,G}} & \Bun_{\tilde H,G},
\end{array}
$$
where the left square is cartesian, thus defining $\ov{\cH}^{\mu,\lambda}_{\tilde H,P}$, and the map $p'_{\tilde H,P}$ sends 
$(\tilde x, \cB, L\hook{} M, \beta: M\,\iso\, M'\mid_{X-x})$ to $(\cB', L(-\lambda_1 x)\hook{} M')$. Here $\cB'=\cB(\mu_1\tilde x+\mu_2\sigma(\tilde x))$. Write also $\cA'=\det M'\,\iso\, \cA(-(\lambda_1+\lambda_2)x)$ and $L'=L(-\lambda_1 x)$. 

 Define the Hecke functor
$$
\H^{\mu,\lambda}: \D(^{\le\lambda_1-\lambda_2}\cZ_1)\to \D(\tilde X\times 
\Bun_{\tilde H,P})
$$
by
$$
\H^{\mu,\lambda}(K)=(\supp\times p_{\tilde H,P})_! (\pr^*\IC_{\ov{\cH}^{\lambda}_G}\otimes (p'_{\tilde H,P})^*K)[(\lambda_2-\lambda_1)-\dim\Bun_G]
$$
We normalize it so that (in view of Theorem~\ref{Th_Hecke_property}) it should preserve perversity. The term $(\lambda_2-\lambda_1)$ appears, because the dimension of $\Bun_{\tilde H,P}$ depends on a connected component. So, 
$$
(\id\times\nu_P)^*\H^{\mu,\lambda}(K)[\dimrel(\nu_P)]\,\iso\, \H^{\mu,\lambda}(\nu_{\cZ}^*K[\dimrel(\nu_{\cZ})])
$$

 We have a commutative diagram
$$
\begin{array}{ccccc}
\cZ_2^{\le \lambda_1-\lambda_2} & \getsup{p_{\cZ}} & 
\ov{\cH}^{\mu,\lambda}_{\tilde H,P}\times_{\cZ_1} \cZ_2 & \toup{p'_{\cZ}} & ^{\le\lambda_1-\lambda_2}\cZ_2\\
\downarrow\lefteqn{\scriptstyle \pi_{2,1}} && \downarrow && \downarrow\lefteqn{\scriptstyle \pi_{2,1}}\\
\tilde X\times \Bun_{\tilde H,P} & \getsup{\supp\times p_{\tilde H,P}}  & \ov{\cH}^{\mu,\lambda}_{\tilde H,P} & \toup{p'_{\tilde H,P}} & ^{\le\lambda_1-\lambda_2}\cZ_1,
\end{array}
$$  
where the right square is cartesian. Here $\cZ_2^{\le i}$ is the stack over $\tilde X\times \Bun_{\tilde H,P}$ whose fibre over $(\tilde x, \cB, 0\to L\to M\to L^*\otimes\cA\to 0)$ is 
$\Hom(L^2, \cA\otimes\Omega( i x))$, where $x=\pi(\tilde x)$. By definition, $p_{\cZ}$ is the map that forgets $(\beta, M')$. 

 Define the Hecke functor
$
\H^{\mu,\lambda}: \D(^{\le\lambda_1-\lambda_2}\cZ_2)\to \D(\cZ_2^{\le\lambda_1-\lambda_2})
$
by
\begin{equation}
\label{Hecke_functor_cZ_2}
\H^{\mu,\lambda}(K)=p_{\cZ !} (\pr^*\IC_{\ov{\cH}^{\lambda}_G}\otimes (p'_{\cZ})^*K)[(\lambda_2-\lambda_1)-\dim\Bun_G],
\end{equation}
it commutes with the functor $(\pi_{2,1})_!$. 

\bigskip\noindent
4.2.2 {\scshape Grothendieck group calculation} Set $K_1=\nu_{\cZ}^*\Aut_{G,\tilde H}[\dimrel(\nu_{\cZ})]$ and $K_2=W_{12}(K_1)$. Remind that $\Bun_{\tilde H,P}$ classifies: $\cB\in\Pic\tilde X$, $L\in\Pic X$, and an exact sequence $0\to L^2\to ?\to \cA\to 0$ on $X$ with $\cC=N(\cB)$ and $\cA=\Omega\otimes\cC^{-1}$.
Let $_c\Bun_{\tilde H,P}\subset \Bun_{\tilde H,P}$ be the open substack given by
\begin{equation}
\label{condition_Sect422}
2\deg L+\deg\cC+ 2(\lambda_1-\lambda_2)<0
\end{equation} 
The projection $_c\Bun_{\tilde H,P}\to\Bun_{\tilde H,G}$ is smooth and surjective. We will derive Theorem~\ref{Th_Hecke_property} from a description of the complex 
$$
(\id\times\nu_P)^*\H^{\mu,\lambda}(\Aut_{G,\tilde H})[\dimrel(\nu_P)]\,\iso\, \H^{\mu,\lambda}(K_1)
$$
over $\tilde X\times{_c\Bun_{\tilde H,P}}$.  The latter will follow from the Hecke property of $K_2$ for (\ref{Hecke_functor_cZ_2}). 

 By equivariance, $\H^{\mu,\lambda}(K_2)$ is the extension by zero from the closed substack $\cZ_2^{\le 0}\hook{} \cZ_2^{\le\lambda_1-\lambda_2}$. 
  Write $\cZ_2^i\subset \cZ_2^{\le i}$ for the open substack given by the condition that $L^2\toup{s} \cA\otimes\Omega(ix)$ does not have a zero at $x$. 
Let $_c\cZ_2^{\le i}$ be the preimage of $\tilde X\times{_c\Bun_{\tilde H,P}}$ under $\pi_{2,1}: \cZ_2^{\le i}\to \tilde X\times\Bun_{\tilde H,P}$. Set 
$_c\cZ_2^i={_c\cZ_2^{\le i}}\cap {\cZ_2^i}$. 
Set also
$$
^dY^i=p_{\cZ}^{-1}(_c\cZ_2^i)\cap (p'_{\cZ})^{-1}(^d\cZ'_2)
$$
Let $\cK^i$ denote the $*$-restrictiction of $\H^{\mu,\lambda}(K_2)$ to $_c\cZ_2^i$. We can similarly define the category $\D^W(_c\cZ_2^i)$, then  $\cK^i\in \D^W(_c\cZ_2^i)$. 

\begin{Lm} 
\label{Lm_cutting_and_calculating}
The complex $\cK^0$ (resp., $\cK^i$ for $i<0$) is placed in non-positive (resp., strictly negative) perverse degrees. The $0$-th perverse cohomology of $\cK^0$ identifies with
\begin{equation}
\label{sum_for_Con_Hecke}
\oplus_{\nu} \H^{\nu,0}(^0K_2)\otimes (V^{\lambda})^*(\nu-\mu)
\mid_{_c\cZ_2^0}, 
\end{equation}
\end{Lm}
\begin{Prf}   
Denote by $^d\cK^i$ the $*$-restriction of 
\begin{equation}
\label{complex_integrant_Hecke}
(\pr^*\IC_{\ov{\cH}^{\lambda}_G}\otimes (p'_{\cZ})^*(^dK_2))[(\lambda_2-\lambda_1)-\dim\Bun_G]
\end{equation}
to $^dY^i$ followed by the direct image $p_{\cZ !}$. Let $S_2$ act on $^d\cK^i$ via its action on $^d K_2$ (cf.~3.6.2). We are reduced to the following lemma.
\end{Prf}

\begin{Lm} 
\label{Lm_in_422_S_2_equivariance_and_Hecke}
The complex $^d\cK^i$ is placed in perverse degrees $\le 0$. The inequality is strict for all terms except $^0\cK^0$. 
The $0$-th perverse cohomology of $^0\cK^0$ is $S_2$-equivariantly isomorphic to (\ref{sum_for_Con_Hecke}).
\end{Lm}
\begin{Prf}
Remind the diagram
$$
\begin{array}{cccc}
^d Y^i  & \toup{p'_{\cZ}} & {^d\cZ'_2} & \toup{^d\chi}\; \A^1\\
\downarrow\lefteqn{\scriptstyle p_{\cZ}}&& \downarrow\lefteqn{\scriptstyle \phi_2}\\
_c\cZ^i_2 && ^d\cP_2
\end{array}
$$
The scheme $^dY^i$ is empty unless $i\le (\lambda_1-\lambda_2)-2d$. 

 Assume $i\le (\lambda_1-\lambda_2)-2d$, then the map $p_{\cZ}: {^d Y^i}\to \cZ^i_2$ can be seen as a (twisted) projection
$$
_c\cZ^i_2\,\ttimes \,(\Grb^{\lambda}_G\cap S^{\lambda'})\to {_c\cZ^i_2}
$$ 
for $\lambda'=(\lambda_1-d,\lambda_2+d)$. 

 Remind that $\ov{\cH}^{\lambda}_G$ is a twisted product $(X\times\Bun_G)\,\ttimes\, \Grb^{\lambda}_G$, where the projection to $\Bun_G$ corresponds to $p_G$. We have $\IC_{\ov{\cH}^{\lambda}_G}\,\iso\, \IC_{X\times\Bun_G}\tilde\boxtimes \cA_{\lambda}$, where $\cA_{\lambda}\in\Sph(\Gr_G)$ is the spherical sheaf on $\Gr_G$ corresponding to $\lambda$. 

 View $_c\cZ^i_2$ as the stack classifying: $\tilde x\in\tilde X$, $\cB\in\Pic\tilde X$, an exact sequence $0\to L\to M\to L^*\otimes\cA\to 0$ on $X$, where $\cC\otimes\cA\,\iso\,\Omega$ and (\ref{condition_Sect422}) holds, and a section $s: L^2\to \cA\otimes\Omega(ix)$ that has no zero at $x$. Here $x=\pi(\tilde x)$ and $\cC=N(\cB)$.  
 
 View $^dY^i$ as the stack over $_c\cZ^i_2$, whose fibre over the above point is the scheme of pairs $(M',\beta)$ with $M'\in\Bun_G$ and $\beta: M\,\iso\, M'\mid_{X-x}$ such that $M'$ is in a position $\le\lambda$ w.r.t. $M$, and $\bar L\subset M'$ is a subbundle. Here $\bar L=L((d-\lambda_1)x)$. 
 
  View $^d\cZ'_2$ as the stack classifying: $\cB'\in\Pic\tilde X$, a modification $L'\subset \bar L$ of line bundles on $X$ with $\deg(\bar L/L')=d$, an exact sequence $0\to \bar L^2\to ?\to \cA'\to 0$, and a section $s: \bar L^2\to \cA'\otimes\Omega$. Here $N(\cB')\otimes\cA'\,\iso\,\Omega$. 
  
  The map $p'_{\cZ}$ sends the above collection to $\cB'=\cB(\mu_1\tilde x+\mu_2\sigma(\tilde x))$, $L(-\lambda_1 x)=L'\subset \bar L=L((d-\lambda_1)x)$, $s: \bar L^2\to \cA'\otimes\Omega$, and $0\to \bar L^2\to ?\to \cA'\to 0$. Note that $\cA'=\cA(-(\lambda_1+\lambda_2)x)$.  
  
  Now it is convenient to think of $^d\cP_2$ as the stack classifying: a modification $L'\subset \bar L$ of line bundles on $X$ with $\deg(\bar L/L')=d$, $\cB'\in\Pic\tilde X$, and a section $s: \bar L^2\to\cA'\otimes\Omega$, where $N(\cB')\otimes\cA'\,\iso\,\Omega$. 
  
  Denote by $p'_{\cP}: {_c\cZ^i_2}\to {^d\cP_2}$ the map sending the above collection to $\cB'=\cB(\mu_1\tilde x+\mu_2\sigma(\tilde x))$, $L(-\lambda_1 x)=L'\subset \bar L=L((d-\lambda_1)x)$, and $s: \bar L^2\to \cA'\otimes\Omega$. It fits into a commutative diagram
$$
\begin{array}{cccc}
^d Y^i  & \toup{p'_{\cZ}} & {^d\cZ'_2} & \toup{^d\chi}\; \A^1\\
\downarrow\lefteqn{\scriptstyle p_{\cZ}}&& \downarrow\lefteqn{\scriptstyle \phi_2}\\
_c\cZ^i_2 &\toup{p'_{\cP}} & ^d\cP_2
\end{array}
$$ 

 If $i\le 0$ then by (\cite{FGV}, 7.2.7(2)) the map 
$$
^d\chi\comp p'_{\cZ}: \cZ^i_2\,\ttimes \,(\Grb^{\lambda}_G\cap S^{\lambda'})\to\A^1
$$
identifies with 
$$
\cZ^i_2\,\ttimes \,(\Grb^{\lambda}_G\cap S^{\lambda'}) \toup{^0\chi\times \chi^{\lambda'}_{\nu}} \; \A^1\times \A^1\toup{sum}\A^1,
$$
where $\nu=(0,i)$, and $\chi^{\lambda'}_{\nu}$ is the notation from \select{loc.cit.} By Theorem~1 from \select{loc.cit.}, the complex
$$
\RG_c(\Grb^{\lambda}_G\cap S^{\lambda'}, \cA_{\lambda}\otimes (\chi^{\lambda'}_{\nu})^*\cL_{\psi})
$$
is placed in degree $\<2\lambda', \check{\rho}\>=\lambda_1-\lambda_2-2d$ and equals $\Hom_{\check{G}}(V^{\lambda}\otimes V^{\nu}, V^{\nu+\lambda'})$. 

A connected component of $^dY^i$ maps to a pair of connected components of $^d\cZ'_2$ and $_c\cZ^i_2$. For such pair of components we have
\begin{equation}
\label{dim_formula_used}
\dim{_c\cZ^i_2}-\dim {^d\cZ'_2}=1+3d+i+2(\lambda_2-\lambda_1)
\end{equation}

  The $*$-restriction of (\ref{complex_integrant_Hecke}) to $^dY^i$ identifies with 
$$
((p'_{\cP})^* (^d\cT)\otimes {^0\chi^*\cL_{\psi}})\tilde\boxtimes (\cA_{\lambda}\otimes (\chi^{\lambda'}_{\nu})^*\cL_{\psi}\mid_{\Grb^{\lambda}_G\cap S^{\lambda'}})[1+(\lambda_2-\lambda_1)+\dimrel(\phi_2)]
$$
The condition (\ref{condition_Sect422}) garantees that, over $p'_{\cZ}(^dY^i)$, the complex $^dK_2$ is placed in the usual cohomological degree $-d-\dim(^d\cZ'_2)$. 
Using (\ref{dim_formula_used}), we learn that $^d\cK^i$ is placed in usual cohomological degree $-\dim{_c\cZ^i_2}+i$. Since $i\le 0$, it is placed in perverse degree $\le 0$, and the inequality is strict unless $i=0$. 

  For $i=0$ we have $\Hom_{\check{G}}(V^{\lambda}, V^{\lambda'})=0$ unless $\lambda=\lambda'$. So, only $^0\cK^0$ contributes to the $0$-th perverse cohomology of $\cK^0$. 
  
  It remains to analize the $0$-th perverse cohomology of the $*$-restriction of $^0\cT$ under $p'_{\cP}:\cZ_2^0\to {^0\cP_2}$. We have to consider the space of sections $t: \pi^* L'\to (\cB')^{\star}$, that is, 
$$
t: \pi^*L\to \cB^*\otimes\Omega_{\tilde X}((\lambda_1-\mu_1)\tilde x+(\lambda_1-\mu_2)\sigma(\tilde x))
$$ 
such that $Nt: L^2\to\cA\otimes\Omega$ has no zero at $x$. This means that $t: \pi^*L\to\cB^*\otimes\Omega_{\tilde X}(a\tilde x-a\sigma(\tilde x))$ for some $a\in\ZZ$ such that
$$ 
a\tilde x-a\sigma(\tilde x)\le (\lambda_1-\mu_1)\tilde x+(\lambda_1-\mu_2)\sigma(\tilde x)
$$
as divisors on $\tilde X$. Our assertion follows.  
\end{Prf}

\medskip

\begin{Rem} In the above proof we have that $i-\dim{_c\cZ^i_2}=-\dim{_c\cZ^0_2}$ does not depend on $i$, so that $\H^{\mu,\lambda}(K_2)\mid_{_c\cZ_2^{\le 0}}$ is placed in the usual cohomological degree $-\dim{_c\cZ_2^0}$. 
\end{Rem}

\begin{Lm} 
\label{Lm_in422_for_Th_1}
The complex $\H^{\mu,\lambda}(K_1)$ over $\tilde X\times {_c\Bun_{\tilde H,P}}$ is placed in perverse degrees $\le 0$, and its $0$-th perverse cohomology identifies with 
$$
\oplus_{\nu} \; \H^{\nu,0}(K_1)\otimes (V^{\lambda})^*(\nu-\mu)
\mid_{\tilde X\times {_c\Bun_{\tilde H,P}}} 
$$
\end{Lm}
\begin{Prf} The intersection of a fibre of $_c\Bun_{P,\tilde H}\to \Bun_{G,\tilde H}$ with each connected component of $_c\Bun_{P,\tilde H}$ is
either connected or empty. So, by Proposition~\ref{Pp_irreducibility}, 
$K_1$ is an irreducible perverse sheaf over each connected component of $_c\Bun_{P,\tilde H}$. 

 Let $_c^0\cZ_2$ be the preimage of $_c\Bun_{P,\tilde H}$ under $\pi_{2,1}: {^0\cZ_2}\to \Bun_{P,\tilde H}$. 
By Proposition~\ref{Pp_Whittaker_functors}, $K_2$ is an irreducible perverse sheaf over each connected component of $_c^0\cZ_2$. So, if $\nu$ is a coweight of $\GO_2^0$ that vanishes in $\pi_1(G)$ then $\H^{\nu,0}(K_2)$ is an irreducible perverse sheaf over each connected component of $_c\cZ_2^{\le 0}$. 

The functor $(\pi_{2,1})_!: \D^W(\cZ_2)\to\D(\cZ_1)$ is exact for the perverse t-structures (and commutes with Hecke functors). Our assertion follows now from
Lemma~\ref{Lm_cutting_and_calculating}. 
\end{Prf}

\medskip
\begin{Prf}\select{of Theorem~\ref{Th_Hecke_property}}

\smallskip\noindent
Remind that the Hecke functor (\ref{Hecke_functor_GH}) commutes with Verdier duality. Since $\Aut_{G,\tilde H}$ is self-dual, our assertion follows from Lemma~\ref{Lm_in422_for_Th_1}. 

 The $S_2$-equivariance statement from Lemma~\ref{Lm_in_422_S_2_equivariance_and_Hecke} combined with (\cite{Ly}, Remark~3) imply the last assertion about generic and special parts. 
Indeed, for $L\in\Pic X, \cB\in\Pic\tilde X$ and $L'=L(-\lambda_1 x)$, $\cB'=\cB(\mu_1 \tilde x+\mu_2 \tilde x)$ we have $\chi(L'\otimes \pi_*\cB')-\chi(L\otimes\pi_*\cB)=\lambda_2-\lambda_1$. 
\end{Prf}

\bigskip\noindent
4.2.3 Let us derive from Theorem~\ref{Th_Hecke_property} that $F_G$ commutes with Hecke operators. 
As in B.1.3, for $\mu=(\mu_1,\mu_2)\in\Lambda_{\HH}$ we have a Hecke functor $\H^{\mu}_{\tilde H}: \D(\Bun_{\tilde H})\to \D(\tilde X\times\Bun_{\tilde H})$. It is given by 
$$
\H^{\mu}_{\tilde H}(K)=
(p'_{\tilde H})^*K[1]
$$ 
for the map $p'_{\tilde H}: \tilde X\times\Bun_{\tilde H}\to \Bun_{\tilde H}$ 
sending $(\tilde x,\cB)$ to $\cB(-\mu_1\tilde x-\mu_2\sigma(\tilde x))$, here $\cB\in\Pic\tilde X$. 

\begin{Cor}
\label{Cor_Hecke_property_GO_2}
 i) For the map $\pi\times\id: \tilde X\times \Bun_G\to X\times\Bun_G$ 
and a dominant coweight $\lambda$ of $G$ we have an isomoprhism of functors
\begin{equation}
\label{iso_Cor_Hecke_GO_2_GL_2}
(\pi\times\id)^*\comp \H^{\lambda}_G\comp F_G\,\iso\, \oplus_{\mu} (\id\boxtimes F_G)\comp  \H^{\mu}_{\tilde H}\otimes (V^{\lambda})^*(\mu)
\end{equation}
from $\D(\Bun_{\tilde H})$ to $\D(\tilde X\times\Bun_G)$. This isomorphism is compatible with the action of $\Sigma=\{1,\sigma\}$ on both sides. It is understood that $\Sigma$ acts on 
$
\oplus_{\mu} \, \H^{\mu}_{\tilde H}\otimes (V^{\lambda})^*(\mu)
$
via the isomorphisms (\ref{iso_action_Sigma_Hecke_sectionD}). So,
$$
\H^{\lambda}_G\comp F_G\,\iso\, \Hom_{\Sigma}(\triv, \oplus_{\mu} (\pi\times\id)_!\comp (\id\boxtimes F_G)\comp  \H^{\mu}_{\tilde H}\otimes (V^{\lambda})^*(\mu))
$$
Here $\id\boxtimes F_G$ is the corresponding functor $\D(\tilde X\times\Bun_{\tilde H})\to \D(\tilde X\times\Bun_G)$. If $\lambda_1-\lambda_2$ is even then (\ref{iso_Cor_Hecke_GO_2_GL_2}) preserves the generic and special parts of $F_G$, otherwise it interchanges them. 

\smallskip\noindent
ii) If $K$ is an automorphic sheaf on $\Bun_{\tilde H}$ corresponding to a rank one local system $\tilde E$ on $\tilde X$ then $F_G(K)\in\D(\Bun_G)$ is an automorphic sheaf corresponding to the local system $(\pi_*\tilde E)^*$. 
\end{Cor}
\begin{Prf}
i) Take $\tilde\mu=(\tilde\mu_1, 0)$ with $\tilde\mu_1=\lambda_1+\lambda_2$. Consider the diagram
$$
\begin{array}{cccccc}
\tilde X\times\Bun_{\tilde H,G} & \getsup{\supp\times p_{\tilde H,G}} & \ov{\cH}^{\tilde\mu,\lambda}_{\tilde H,G} & \toup{p'_{\tilde H,G}} & \Bun_{\tilde H,G}\ & \toup{\gq}\,  \Bun_{\tilde H}\\
\downarrow\lefteqn{\scriptstyle\id\times\gp} && \downarrow && \downarrow\lefteqn{\scriptstyle\gp}\\
\tilde X\times\Bun_G & \gets & \tilde X\times_X \ov{\cH}^{\lambda}_G & \to & \Bun_G,
\end{array}
$$
where both squares are cartesian. By Theorem~\ref{Th_Hecke_property}, for $K\in\D(\Bun_{\tilde H})$ we get an isomorphism
\begin{multline*}
(\pi\times\id)^*\H^{\lambda}_G F_G(K)\,\iso\, \\
\oplus_{\nu} (\id\times\gp)_! (\H^{\nu,0}(\Aut_{G,\tilde H})\otimes (V^{\lambda})^*(\nu-\tilde\mu)\otimes
\H^{-\tilde\mu}_{\tilde H}(K))[-1-\dim\Bun_{\tilde H}]\,\iso\\
\oplus _{\nu} (\id\boxtimes F_G) \H^{\nu-\tilde\mu}_{\tilde H}(K)\otimes (V^{\lambda})^*(\nu-\tilde\mu)
\end{multline*}
The assertion about generic and special parts also follows from Theorem~\ref{Th_Hecke_property}. 
\end{Prf}

\begin{Cor} 
For the map $\pi\times\id: \tilde X\times\Bun_{\tilde H}\to X\times \Bun_{\tilde H}$ and a dominant coweight $\lambda$ of $G$ we have 
\begin{equation}
\label{iso_Cor_5_Hecke}
(\pi\times\id)^*\comp (\id\boxtimes F_{\tilde H})\comp \H^{\lambda}_G\,\iso\, 
\oplus_{\mu} \H^{\mu}_{\tilde H} \comp F_{\tilde H}\otimes (V^{\lambda})(-\mu)
\end{equation}
This isomorphism is compatible with the action of $\Sigma=\{1,\sigma\}$ on both sides. So,
\begin{equation}
\label{iso_invariants_sigma_Cor5}
(\id\boxtimes F_{\tilde H})\comp \H^{\lambda}_G\,\iso\, \Hom_{\Sigma}(\triv, \;
\oplus_{\mu}\;  (\pi\times\id)_!\comp \H^{\mu}_{\tilde H} \comp F_{\tilde H}\otimes (V^{\lambda})(-\mu))
\end{equation}
Here $\id\boxtimes F_{\tilde H}: \D(X\times\Bun_G)\to \D(X\times \Bun_{\tilde H})$ is the corresponding functor. If $\lambda_1-\lambda_2$ is even then (\ref{iso_Cor_5_Hecke}) preserves the generic and special parts of $F_{\tilde H}$, otherwise it interchanges them.
\end{Cor}
\begin{Prf}
For a coweight $\tilde\mu$ such that $\lambda-\tilde\mu$ vanishes in $\pi_1(G)$
we have a commutative diagram
$$
\begin{array}{cccccccc}
\Bun_G \, \getsup{p'_G} & \ov{\cH}^{\lambda}_G & \gets & \ov{\cH}^{\tilde\mu, \lambda}_{\tilde H,G} & \toup{\supp\times p'_{\tilde H,G}} & \tilde X\times\Bun_{G,\tilde H} & \toup{\id\times\gq} & \tilde X\times\Bun_{\tilde H}\\
& \downarrow\lefteqn{\scriptstyle \supp\times p_G} && \downarrow\lefteqn{\scriptstyle \supp\times p_{\tilde H, G}} &&&
\swarrow\lefteqn{\scriptstyle \id\times p'_{\tilde H}}\\
& X\times\Bun_G & \getsup{\pi\times\gp} & \tilde X\times\Bun_{G,\tilde H} & \toup{\id\times\gq} & \tilde X\times\Bun_{\tilde H},
\end{array}
$$  
where the left square is cartesian. Here $\supp\times p'_{\tilde H,G}$ sends $(\tilde x, \cB, M\,\iso\, M'\mid_{X-x})$ to $(\tilde x, M', \cB')$ with $\cB'=\cB(\tilde\mu_1\tilde x+\tilde\mu_2\sigma(\tilde x))$. The map $p'_{\tilde H}$ sends $(\tilde x, \cB')$ to $\cB'(-\tilde\mu_1\tilde x-\tilde\mu_2\sigma(\tilde x))$. 

 In this notation we have
$$
\H^{-\tilde\mu, -w_0(\lambda)}(K)\,\iso\, (\supp\times p'_{\tilde H,G})_!(\IC_{\ov{\cH}^{\tilde\mu,\lambda}_{\tilde H,G}}\otimes p_{\tilde H,G}^*K)[ -\dim\Bun_{\tilde H,G}]
$$
So, for $K\in\D(\Bun_G)$ the above diagram yields an isomorphism
\begin{multline}
\label{iso_utility_in_Cor5}
(\pi\times\id)^* (\id\boxtimes F_{\tilde H})\H^{\lambda}_G(K)\,\iso\\
(\id\times p'_{\tilde H})_!(\id\times \gq)_!(\H^{-\tilde\mu, -w_0(\lambda)}(\Aut_{G,\tilde H})\otimes \gp^*K)[-\dim\Bun_G]
\end{multline}
By Theorem~\ref{Th_Hecke_property}, 
$$
\H^{-\tilde\mu, -w_0(\lambda)}(\Aut_{G,\tilde H})\,\iso\, \oplus_{\nu} \H^{\nu,0}(\Aut_{G,\tilde H})\otimes (V^{-w_0(\lambda)})^*(\nu+\tilde\mu)
$$
So, the RHS of (\ref{iso_utility_in_Cor5}) identifies with 
$
\oplus_{\nu} \H^{-\tilde\mu-\nu} F_{\tilde H}(K)\otimes V^{\lambda}(\nu+\tilde\mu)
$. 
\end{Prf}

\bigskip\noindent
4.3  Remind that $\cS_P$ classifies: $L\in\Bun_1, \cA\in\Bun_1$ and $L^2\toup{s}\cA\otimes\Omega$. We have open immersion $j_d: \Pic X\times X^{(d)}\hook{} \cS_P$ sending $(L,D)$ to $L, \cA=\Omega^{-1}\otimes L^2(D)$ with the canonical inclusion $L^2\hook{}\cA\otimes\Omega$.  
 
\begin{Def} Let $\tilde E$ be a rank one local system on $\tilde X$. Remind that $A\tilde E$ denotes the corresponding automorphic local system on $\Pic\tilde X$ (cf. Definition~\ref{Def_aut_loc_sys}). Set $E=\pi_*\tilde E$. Define the perverse sheaf $\tilde E_H\in P(\Bun_H)$ by
$$
\tilde E_H=\rho_{H\, !}A\tilde E[\dim\Bun_H]
$$ 
\end{Def}  
  
\begin{Lm} 
\label{Lm_strange}
For $d\ge 0$ we have canonically 
$$
j_d^*\, F_{\cS}(\tilde E_H)\,\iso\, AN(\tilde E^*)\boxtimes (E^*)^{(d)}\otimes AN(\tilde E)_{\Omega}[\dim\cS_P]
$$
\end{Lm}
\begin{Prf}
The stack $\cV_{\tilde H,P}$ classifies: $\cB\in\Pic\tilde X$, $L\in\Pic X$, and a map $t:L\otimes\pi_*\cB\to\Omega$. The datum of $t$ is equivalent to a datum of $t: \pi^*L\to \cB^{\star}$. We have $\cC=N(\cB)$. 

 Let $\tilde\gp_{\cV}$ be the composition $\cV_{\tilde H,P}\to\cV_{H,P}\toup{\gp_{\cV}}\cS_P$, it sends the above point to $(L, \cA=\Omega\otimes\cC^{-1}, s)$, where $s: L^2\to \cA\otimes\Omega$ equals the norm of $t$.  We have a cartesian square
$$
\begin{array}{ccc}
\Pic X\times \tilde X^{(d)} & \toup{j_{\cV,d}} & \cV_{\tilde H,P} \\
\downarrow\lefteqn{\scriptstyle \id\times\pi} && \downarrow\lefteqn{\scriptstyle \tilde\gp_{\cV}}\\
\Pic X\times X^{(d)} & \toup{j_d} & \cS_P,
\end{array}
$$
where $j_{\cV,d}$ sends $(L,\tilde D)$ to $L, \cB=(\pi^* L(\tilde D))^{\star}$ with the canonical inclusion $t: \pi^*L\hook{}\cB^{\star}$. We have canonically 
$$
j_{\cV,d}^*\, \gq_{\cV}^*A\tilde E\,\iso\, AN(\tilde E^*)\boxtimes (\tilde E^*)^{(d)}\otimes AN(\tilde E)_{\Omega}
$$
Our assertion follows by i) of Lemma~\ref{Lm_evident}.
\end{Prf}

\medskip

 Since $\Pic^d\tilde X$ is connected, the covering $\rho_H: \Pic^d\tilde X\to\Bun^d_H$ is nontrivial, and $\cN$ is a nontrivial local system on each $\Bun^d_H$. 

\begin{Lm} The following conditions are equivalent:
\begin{itemize}
\item $E$ is irreducible
\item $\tilde E$ does not descend to a rank one local system on $X$ 
\item $A\tilde E$ does not descend with respect to $\Pic\tilde X\toup{N} \Pic X$ 
\item the local system $\rho_{H\, !}A\tilde E$ on $\Bun_H$ is irreducible. \QED
\end{itemize}
\end{Lm}

\begin{Def} 
\label{Def_Aut_E}
For an irreducible rank $k$ local system $W$ on $X$ denote by $\Aut_W$ the corresponding automorphic sheaf on $\Bun_k$ normalized as in \cite{FGV}. By \select{loc.cit.}, if $\lambda$ is the dominant weight of the standard representation of $\check{\GL}_k$ then $\H^{\lambda}_{\GL_k}(\Aut_W)\,\iso\, W^*\boxtimes\Aut_W[1]$. 
\end{Def}

 Remind the normalization of $\Aut_W$ for $k=2$. The map $\nu_P:\Bun_P\to\Bun_G$ sends $0\to L^2\to ?\to \cA\to 0$ to $M$ (included into $0\to L\to M\to L^*\otimes\cA\to 0$).  
First, one considers the complex, say $\cK_W$, on $\Pic X\times X^{(d)}$ whose fibre at $L,D$ is 
$$
(A\det W)_{L\otimes\Omega^{-1}}\otimes W^{(d)}_D[\dim\cS_P]
$$ 
Then $\Four_{\psi}((j_d)_!\cK_W)$ identifies with $\nu_P^*\Aut_W[\dimrel(\nu_P)]$ over the components of 
$\Bun_P$ for which $\deg(\cA\otimes\Omega)\ge 2\deg L$. The sheaf $\Aut_W$ is perverse and irreducible on each connected component of $\Bun_G$. 

 To fix notation for Eisenstein series, denote by $\Bunb_P$ the stack classifying $M\in\Bun_2, L\in\Pic X$ and an inclusion of coherent sheaves $L\hook{} M$. Write $\Bunb^{d,d_1}_P$ for the connected component of $\Bunb_P$ given by $\deg L=d_1$ and $\deg M+\deg\Omega=2d_1+d$.
We have a diagram
$$
\Pic X\times\Pic X\getsup{\bar\gq_P}\Bunb_P\toup{\bar\gp_P} \Bun_G,
$$
where $\bar\gq_P$ sends $(L\subset M)$ to $(L, L^{-1}\otimes\det M)$, and $\bar\gp_P$ is the projection. For rank one local systems $E_1,E_2$ on $X$ set
$$
\Aut_{E_1\oplus E_2}=(AE_2)^{-1}_{\Omega}\otimes
(\bar\gp_P)_!\bar\gq_P^*(AE_1\boxtimes AE_2)[\dim\Bunb_P]
$$
This normalization is compatible with the above in the following sense. If $E_1$ and $E_2$ are not isomorphic then $\Four_{\psi}((j_d)_{!*}\cK_{E_1\oplus E_2})$ descends (over some open substack of $\Bun_P$) to $\Aut_{E_1\oplus E_2}$, and we have the functional equation $\Aut_{E_1\oplus E_2}\,\iso\, \Aut_{E_2\oplus E_1}$ (cf. \cite{BG}). Write 
\begin{equation}
\label{Eis_series_decomp}
\Aut_{E_1\oplus E_2}\,\iso\, \oplus_{(d,d_1)\in\ZZ^2} \Aut_{E_1\oplus E_2}^{d,d_1},
\end{equation}
where $\Aut_E^{d,d_1}$ is the contribution of $\Bunb_P^{d,d_1}$. 
 
 Remind the map $\alpha_P: \Bun_P\to\Pic X$ sending $0\to \Sym^2 L\to ?\to \cA\to 0$ to $L$. Let $^0\cS_P\subset\cS_P$ be the open substack classifying inclusions $L^2\hook{}\cA\otimes\Omega$ with $L, \cA\in\Pic X$. 

\begin{Pp} 
\label{Pp_lift_for_GL2}
1) If $E$ is irreducible then, over the connected components of $\Bun_P$ given by $\deg L<0$, there exists an isomorphism
\begin{equation}
\label{iso_used_Pp3_1}
F_{P,\psi}(\tilde E_H)\,\iso\, \nu_P^*\Aut_{E^*}\otimes (A\cE_0)_{\Omega}\otimes
\alpha_P^*A\cE_0[\dimrel(\nu_P)]
\end{equation}
So, (\ref{iso_gr_k}) gives rise to an isomorphism of perverse sheaves on $\Bun_G$
\begin{equation}
\label{iso_used_Pp3_2}
F_G(\tilde E_H)\,\iso\, \Aut_{E^*}\otimes (A\cE_0)_{\Omega}
\end{equation}
2) Assume $\tilde E=\Qlb$.  Then over the components of $\cS_P$ given by $\deg(\cA\otimes\Omega)-2\deg L>3g-3$ the sheaf $F_{\cS}(\tilde E_H)$ is perverse, the Goresky-MacPherson extension from $^0\cS_P$. Both (\ref{iso_used_Pp3_1}) (for $\deg L$ small enough) and (\ref{iso_used_Pp3_2}) remain valid, where now
$E\,\iso\, \Qlb\oplus \cE_0$. 
\end{Pp}
\begin{Prf}
1) Since $N_! (A\tilde E)=0$, our assertion follows from Lemma~\ref{Lm_strange} combined with Corolary~\ref{Cor_relation_F_P_and_F_G}. 

\smallskip\noindent
2) The components of $\cV_{\tilde H,P}$ given by $\deg\cB+2\deg L< 0$ are smooth (this is where the  connectedness of $\tilde X$ is essential!) The fibres of $N:\Pic\tilde X\to\Pic X$ are of dimension  
$g-1$, so over the corresponding components of $\cS_P$ the map $\tilde\gp_{\cV}: \cV_{\tilde H,P}\to \cS_P$ is small. The first assertion follows. The second one is obtained from Lemma~\ref{Lm_strange}. 
\end{Prf}

\medskip
\begin{Rem}
\label{Rem_not_refinement}
 i) Proposition~\ref{Pp_lift_for_GL2} implies that the constant term $\CT_P(\Aut_{\Qlb\oplus \cE_0})$ is essentially the cohomology of the Prym variety (cf. A.1).\\
ii) The formula (\ref{iso_used_Pp3_2}) for $\tilde E=\Qlb$ is a version of the classical theorem of Siegel (its proof given by A.~Weil can be found in \cite{We}, a version proved by Waldspurger is found in (\cite{Wa}, Sect.~I.5, Proposition 2)).\\
iii) If $\tilde E=\Qlb$ then $\Sigma=\{1,\sigma\}$ acts naturally on $\tilde E_H$ and, hence, on $F_G(\tilde E_H)$. Let $\Sigma$ acts on $\Aut_{\cE_0\oplus \Qlb}$ via (\ref{iso_used_Pp3_2}). The $\sigma$-invariants of $\Aut_{\cE_0\oplus \Qlb}$ are $\oplus_{d, d_1} \Aut_{\cE_0\oplus \Qlb}^{d,d_1}$, the sum over $(d,d_1)\in\ZZ^2$ with $d_1$ even.

\smallskip\noindent
iv) The stack $\Bun_{G,\tilde H}\times_{\Bun_G}\Bun_{G,\tilde H}$ splits as a disjoint union of the open substacks $\cU^0\sqcup \cU^1$,  where $\cU^a$ is given by the condition that $\cB_1\otimes\cB_2^{-1}\in\Bun_{U_{\pi}}^a$ for a point 
$(\cB_1,\cB_2\in\Pic\tilde X, M\in\Bun_2, N(\cB_1)\otimes\det M\,\iso\, N(\cB_2)\otimes\det M\,\iso\,\Omega)$ of  $\Bun_{G,\tilde H}\times_{\Bun_G}\Bun_{G,\tilde H}$ (cf. A.1). So, the restriction of $F_G(\tilde E_H)$ under $\Bun_{G,\tilde H}\to\Bun_G$ is naturally a direct sum $\cK^0\oplus \cK^1$, where $\cK^a$ is the contribution of $\cU^a$. If $\tilde E=\Qlb$ then (\ref{Eis_series_decomp}) is not a refinement of the decomposition $\cK^0\oplus \cK^1$.

 To see this, consider the line bundle $\cE^{(d)}$ on $X^{(d)}$, the $d$-th symmetric power of $\cE$. Its tensor square is canonically trivialized, so it defines a $\mu_2$-torsor $_{\pi} X^{(d)}\to X^{(d)}$. A fibre of the latter map over $D\in X^{(d)}$ can also be seen as the set of connected components of the stack of pairs $(\cB,\kappa)$, where $\cB\in\Pic\tilde X$, $\kappa: N(\cB)\,\iso\,\cO(D)$. The restriction of the covering $_{\pi}X^{(d)}\to X^{(d)}$ under $\pi: \tilde X^{(d)}\to X^{(d)}$ has a distinguised section, and 
$\pi_*(\Qlb^{(d)})\mid_{_{\pi} X^{(d)}}\,\iso\, K^0\oplus K^1$, where $K^0$ is the contribution of the distinguished section. 

 Remind that $(\Qlb\oplus\cE_0)^{(d)}\,\iso\, \oplus_{k=0}^d 
(\sym_{k,d-k})_!(\Qlb\boxtimes \cE_0^{(d-k)})$, where $\sym_{k, d-k}: X^{(k)}\times X^{(d-k)}\to X^{(d)}$ is the sum of divisors. So, 
$$
K^0\oplus K^1\,\iso\, (\oplus_{k=0}^d (\sym_{k,d-k})_!\Qlb\boxtimes \cE_0^{(d-k)})\mid_{_{\pi}X^{(d)}},
$$ 
but the RHS is \select{not} a refinement of the decomposition of the LHS. \end{Rem}

\medskip\noindent
4.4 {\scshape Local Rankin-Selberg type convolutions}

\medskip\noindent
4.4.1 Remind the following Laumon's construction for $\GL_2$. Let $\Bun'_2$ be the stack classifying $M\in\Bun_2$ with nonzero section $\Omega\hook{} M$. To a local system $E$ on $X$ one associates a complex $\Laum_E$ on $\Bun'_2$ defined as follows. 

 Let $\cQ$ be the stack classifying collections $(L_1\subset L_2\subset M)$ with $L_1\,\iso\,\Omega$, $L_2/L_1\,\iso\, \cO_X$, where $L_2\subset M$ is a modification of locally free $\cO_X$-modules of rank 2. Let $\ev_{\cQ}:\cQ\to\AA^1$ be the map sending the above point to the class of $0\to L_1\to L_2\to L_2/L_1\to 0$. Let $\gq_{\cQ}:\cQ\to\Sh_0$ be the map sending the above point to $M/L_2$, here $\Sh_0$ is the stack of torsion sheaves on $X$. Write $\cL_E$ for Laumon's sheaf corresponding to $E$ (\cite{FGV}). Let $\gp_{\cQ}:\cQ\to\Bun'_2$ be the map forgetting $L_2$. 
Set 
$$
\Laum_E=\gp_{\cQ\, !}(\gq_{\cQ}^*\cL_E\otimes\ev_{\cQ}^*\cL_{\psi})[\dim \cQ]
$$
Consider the map $\gq': \Bunb_P\to \Pic X\times \Bun'_2$ sending $(L\subset M)$ to $(L\otimes\Omega^{-1}, \Omega\subset M\otimes L^{-1}\otimes\Omega)$. For local systems $E, E_1$ on $X$, where $E_1$ is of rank one, set 
$$
\Laum_{E,E_1}=(A E_1)_{\Omega}\otimes(\gq')^*(A(E_1\otimes\det E)\boxtimes \Laum_E)[\dim\Pic X]
$$
Let $\Laum_{E,E_1}^{d,d_1}$ denote the restriction of $\Laum_{E,E_1}$ to the connected component $\Bunb^{d,d_1}_P$of $\Bunb_P$ given by $\deg M+\deg\Omega=2d_1+d$ and $\deg L=d_1$. Remind the projection $\bar\gp_P:\Bunb_P\to\Bun_G$. 
By (\cite{FGV}, 7.9), for $E$ irreducible and $d\ge 0$ we have 
$$
\Laum^{d,d_1}_{E,\Qlb}\,\iso\, \bar\gp_P^*\Aut_E
[\dimrel(\bar\gp_P)]
$$
over $\Bunb^{d, d_1}_P$, and
$$
\bar\gp_{P\, !}\Laum_{E,E_1}\,\iso\, \Aut_{E_1\oplus\Qlb}\otimes\Aut_E[-\dim\Bun_G]
$$ 

 Denote by 
$$
\mult^{d,d_1}:\Pic^{d_1} X\times \tilde X^{(d)}\to \Pic^{2d_1+d}\tilde X
$$ 
the map sending $(L, \tilde D)$ to $(\pi^*L)(\tilde D)$. Note that $\mult^{0, d_1}$ is a (representable) Galois $S_2$-covering over its image, the corresponding automorphism of $\Pic^{d_1} X$ sends $L$ to $L\otimes\cE$. 
Let $\epsilon: \Pic\tilde X\,\iso\,\Pic\tilde X$ be the involution sending $\cB$ to $\cB^{\star}$. The following is closely related to the main result of \cite{Ly3}. 

\begin{Th} 
\label{Th_Rankin_Selberg_GO_2_GL_2}
For any local systems $E,E_1$ on $X$ with $\rk(E_1)=1$, $\rk(E)=2$ there is an isomorphism
\begin{equation}
\label{iso_Th_2}
F_{\tilde H}(\bar\gp_{P\, !}\Laum_{E,E_1}^{d,d_1})\iso (A\det E)_{\Omega^{-1}}\otimes\epsilon_!\mult^{d,d_1}_!(A(\det E\otimes E_1\otimes\cE_0)\boxtimes (\pi^*E)^{(d)}) [d+g-1]
\end{equation}
depending only on a choice of (\ref{iso_gr_k}). If $d$ is even (resp.,  odd) then 
$F_{\tilde H,s}$ (resp., $F_{\tilde H,g}$) sends $\bar\gp_{P\, !}\Laum_{E,E_1}^{d,d_1}$ to zero. 
In particular, for $E$ irreducible we have 
\begin{multline*}
F_{\tilde H}(\Aut_{E_1\oplus\Qlb}\otimes\Aut_E)\,\iso\\
\oplus_{d\ge 0}  (A\det E)_{\Omega^{-1}}\otimes\epsilon_!\mult^{d, d_1}_!(A(\det E\otimes E_1\otimes\cE_0)\boxtimes (\pi^*E)^{(d)}) [d+g-1+\dim\Bun_G],
\end{multline*}
where $d_1$ is a function of a connected component of $\Bun_{\tilde H}$ given by $2d_1+d=\deg(\cB^{\star})$ for $\cB\in\Bun_{\tilde H}$. The sum is over $d\ge 0$ such that $d_1\in\ZZ$. 
\end{Th}

\begin{Rem} Write $s^a: X^{(a)}\times\tilde X^{(d-2a)}\to \tilde X^{(d)}$ for the map sending $(D,\tilde D)$ to $\pi^*D+\tilde D$. For any rank 2 local system $E$ on $X$ the sheaf $(\pi^*E)^{(d)}$ admits a filtration with succesive quotients being 
$$
s^a_!((\det E)^{(a)}\boxtimes \pi^*(E^{(d-2a)}))
$$ 
for $0\le a \le \frac{d}{2}$. This follows from the fact that for a 2-dimensional $\Qlb$-vector space $E$ we have
$$
\Sym^{a+b} E\otimes \Sym^a E\,\iso\, \oplus_{i=0}^a (\det E)^{\otimes i}\otimes \Sym^{2a+b-2i} E
$$
On the other hand, the complex $\bar\gp_{P\, !}\Laum_{E,E_1}^{d,d_1}$ has a filtration indexed by $a\ge 0$ coming from a stratification of $\Bunb_P$. The corresponding stratum of $\Bunb_P$ is given by the condition that there is a divisor $D\in X^{(a)}$ such that $L(D)\subset M$ is a subbundle. One may check that the corresponding graded of the left and the right hand side of (\ref{iso_Th_2}) coincide. 
\end{Rem}     
 
\bigskip\noindent
4.4.2 Following \cite{FGV}, for a local system $E$ on $X$ denote by $\Av^d_E: \D(\Bun_G)\to\D(\Bun_G)$ the following averaging functor. Let $\Mod_2^d$ be the stack classifying a modification $M\subset M'$ of rank 2 vector bundles on $X$ with $\deg(M'/M)=d$. Let $\gs: \Mod_2^d\to\Sh_0$ be the map sending this point to $M'/M$. We have a diagram
$$
\Bun_G \getsup{p} \Mod_2^d \toup{p'} \Bun_G,
$$
where $p$ (resp., $p'$) sends
$(M\subset M')$ to $M$ (resp., to $M'$). Set $\Av_E^d(K)=p'_!(p^*K\otimes \cL^d_E)[2d](d)$. 

By (\cite{FGV}, 9.5), for an irreducible rank 2 local system $W$ on $X$ we have 
$$
\Av^d_E(\Aut_W)\,\iso\, \Aut_W\otimes \RG(X^{(d)}, (E\otimes W^*)^{(d)})[d]
$$

 Write $S_d$ for the symmetric group on $d$ elements, set $\Sigma^d=\Aut_{X^d}(\tilde X^d)$. We have a semi-direct product
$\Sigma^d\rtimes S_d$ acting on $\tilde X^d$, it fits into an exact sequence $1\to \Sigma^d\to  \Sigma^d\rtimes S_d\to S_d\to 1$.
For a dominant coweight $\lambda$ of $G$ the functor 
$$
\oplus_{\mu_1,\ldots,\mu_d} \H^{\mu_1}_{\tilde H}\comp\ldots\comp\H^{\mu_d}_{\tilde H}\otimes (V^{\lambda})(-\mu_1)\otimes\ldots\otimes (V^{\lambda})(-\mu_d)
$$
is naturally a functor from $\D(\Bun_{\tilde H})$ to the equivariant derived category $\D^{\Sigma^d\rtimes S_d}(\tilde X^d\times\Bun_{\tilde H})$. So, we can introduce the functor  
$$
(\H^{\lambda}_{\tilde H, G})^{\boxtimes d}: \D(\Bun_{\tilde H})\to \D^{S_d}(X^d\times\Bun_{\tilde H})
$$ 
given by
$$
(\H^{\lambda}_{\tilde H, G})^{\boxtimes d}=\Hom_{\Sigma^d}(\triv, (\pi\times\id)_! \oplus_{\mu_1,\ldots,\mu_d} \H^{\mu_1}_{\tilde H}\comp \ldots\comp\H^{\mu_d}_{\tilde H} \otimes (V^{\lambda})(-\mu_1)\otimes\ldots\otimes (V^{\lambda})(-\mu_d)),
$$
where by abuse of notation we have written $\pi\times\id: \tilde X^d\times\Bun_{\tilde H}\to X^d\times\Bun_{\tilde H}$ for the projection. 

  For a local system $E$ on $X$ let $\Av^d_E: \D(\Bun_{\tilde H})\to\D(\Bun_{\tilde H})$ denote
the averaging functor given by
$$
\Av^d_E(K)=\Hom_{S_d}(\triv, (\pr_2)_! (\pr_1^* E^{\boxtimes d}\otimes (\H^{\lambda}_{\tilde H,G})^{\boxtimes d}(K)),
$$  
where $\lambda=(1,0)$, and 
$\pr_i$ are the two projections from
$X^d\times\Bun_{\tilde H}$ to $X^d$ and $\Bun_{\tilde H}$ respectively.
  
\begin{Pp} 
\label{Pp_avegaring_commutes}
For any local system $E$ on $X$ we have a canonical isomorphism of functors
$$
F_{\tilde H} \comp \Av^d_E\,\iso\, \Av^d_E\comp F_{\tilde H}
$$  
from $\D(\Bun_G)$ to $\D(\Bun_{\tilde H})$. If $d$ is even then this isomorphism preserves the generic and special parts of $F_{\tilde H}$, otherwise it interchanges them. 
\end{Pp}
\begin{Prf}
Take $\lambda=(1,0)$. By (\cite{G}, 1.8), the functor $(\H^{\lambda}_G)^{\boxtimes d}$ maps $\D(\Bun_G)$ to the equivariant derived category $\D^{S_d}(X^d\times\Bun_G)$. We have a canonical isomorphism of functors
from $\D(\Bun_G)$ to itself
$$
\Av^d_E\,\iso\, \Hom_{S_d}(\triv, \, (\pr_2)_!(\pr_1^*E^{\boxtimes d}\otimes (\H^{\lambda}_G)^{\boxtimes d}))
$$
where $\pr_i$ are the two projections from $X^d\times\Bun_G$ to $X^d$ and $\Bun_G$ respectively. Applying (\ref{iso_invariants_sigma_Cor5}) $d$ times
we get a $S_d$-equivariant isomorphism
$$
(\id\boxtimes F_{\tilde H})\comp (\H^{\lambda}_G)^{\boxtimes d} \,\iso\, (\H^{\lambda}_{\tilde H, G})^{\boxtimes d} \comp F_{\tilde H},
$$
where $\id\boxtimes F_{\tilde H}: \D(X^d\times\Bun_G)\to \D(X^d\times\Bun_{\tilde H})$ is the corresponding functor. If $d$ is even then this isomorphism preserves the generic and special parts of $F_{\tilde H}$, otherwise it interchanges them. 
Our assertion follows.  
\end{Prf}  
  
\medskip\noindent
\begin{Prf}\select{of Theorem~\ref{Th_Rankin_Selberg_GO_2_GL_2}}
\Step 1 Case $d=0$. Let $\cX_1$ be the stack classifying $L\in\Pic^{d_1} X$, $\cB\in\Pic\tilde X$, and an isomorphism $\xi: L^2\otimes\cC\,\iso\,\Omega^2$, where $\cC=N(\cB)$. We have a diagram of projections
$$
\Pic^{d_1} X \getsup{q_1} \cX_1 \toup{p_1} \Pic\tilde X,
$$
where $q_1$ (resp., $p_1$) sends the above point to $L$ (resp., to $\cB$). 

Let $p_{\cX}: \cX\to\cX_1$ denote the stack whose fibre over a point of $\cX_1$ is the stack of exact sequences of $\cO_X$-modules 
\begin{equation}
\label{seq_for_Th_2}
0\to L\to M\to L\otimes\Omega^{-1}\to 0
\end{equation}
Consider the diagram
$$
\Bun_{P_2}\getsup{\tau_P}\Bun_{P,\tilde H}\getsup{q_{\cX}} \cX\toup{\ev_{\cX}} \A^1,
$$
where $\ev_{\cX}$ is the map sending a point of $\cX$ to the class of (\ref{seq_for_Th_2}), and $q_{\cX}$ 
sends a point of $\cX$ to the collection $(L\subset M, \cB, \, \xi: \cC\otimes\det M\,\iso\, \Omega)$ with $\cC=N(\cB)$. Using (\ref{iso_gr_k}) and Corolary~\ref{Cor_relation_F_P_and_F_G} we get
$$
F_{\tilde H}(\bar\gp_{P\, !}\Laum_{E,E_1}^{0,d_1})\,\iso\, (A\det E)_{\Omega^{-1}}\otimes p_{1\, !}(p_{\cX\, !}(q_{\cX}^*\tau_P^*S_{P,\psi}\otimes\ev_{\cX}^*\cL_{\psi})\otimes q_1^*A(E_1\otimes\det E\otimes\cE_0))[b],
$$
where $b=\dim\Bun_{\SO_2}-\dim\Bun_{P_2}$, here $\dim\Bun_{P_2}$ is the dimension of (the unique) connected component of $\Bun_{P_2}$ containing $\tau_P(q_{\cX}(\cX))$. 

 Consider the map $s_1: \Pic^{d_1} X\to \cX_1$ sending $L$ to the collection $(L,\cB,\xi)$, where $\cB=(\pi^*L)^{\star}$ and $\xi:L^2\otimes\cC\,\iso\,\Omega^2$ is the natural isomorphism with $\cC=N(\cB)$. 
The diagram commutes
$$
\begin{array}{ccccc}
\Pic^{d_1} X & \getsup{q_1} & \cX_1 & \toup{p_1} & \Pic\tilde X\\
& \nwarrow\lefteqn{\scriptstyle\id} & \uparrow\lefteqn{\scriptstyle s_1} && \uparrow\lefteqn{\scriptstyle\epsilon}\\
&&\Pic^{d_1} X & \toup{\mult^{0,d_1}} & \Pic^{2d_1}\tilde X
\end{array}
$$
Remind that $S_2$ acts on $S_{P,\psi}$ (cf. 2.1). By definition of $S_{P,\psi}$, we get a $S_2$-equivariant isomorphism
$$
s_{1 \, !}\Qlb[g-1]\,\iso\,
p_{\cX\, !}(q_{\cX}^*\tau_P^*S_{P,\psi}\otimes\ev_{\cX}^*\cL_{\psi})[b]
$$  
 
 Note that $s_1$ is a (respresentable) $S_2$-covering over its image. We have a 2-automorphism $\eta$ of the identity functor $\id_{\cX_1}$ acting on $(L,\cB,\xi)$ as $-1$ on $L$ and trivially on $\cB$. Since $\eta$ acts as $-1$ on $\Hom_{S_2}(\sign, s_{1\, !}\Qlb)$ and trivially on $q_1^*A(E_1\otimes \det E\otimes\cE_0)$, it follows that
\begin{equation}
\label{complex_vanishing_in_Th_2}
p_{1\, !}(q_1^*A(E_1\otimes \det E\otimes\cE_0)\otimes \Hom_{S_2}(\sign, s_{1\, !}\Qlb))=0
\end{equation}
We have used that $\RG_c(B(\mu_2), W)=0$, where $W$ is the nontrivial rank one local system on $B(\mu_2)$ corresponding to the $S_2$-covering $\Spec k\to B(\mu_2)$. 
  
  Note that the genus of $\tilde X$ is odd. For a point of $\cX$ as above, we have $\chi(L\otimes\pi_*\cB)=0\mod 2$. By (\cite{Ly}, Remark~3), we get $(S_{P,\psi})^{S_2}\,\iso\, S_{P,\psi, g}$ over the connected component of $\Bun_{P_2}$ containing $\tau_P(q_{\cX}(\cX))$.  From (\ref{complex_vanishing_in_Th_2}) it follows that $F_{\tilde H,s}(\bar\gp_{P\, !}\Laum_{E,E_1}^{0,d_1})=0$. So, 
$$
F_{\tilde H,g}(\bar\gp_{P\, !}\Laum_{E,E_1}^{0,d_1})\,\iso\, (A\det E)_{\Omega^{-1}}\otimes\epsilon_!\mult^{0,d_1}_! A(E_1\otimes \det E\otimes\cE_0)
[g-1]
$$

\medskip  
\Step 2  For $d\ge 0$ we have $\Av^d_E(\bar\gp_!\Laum^{0,d_1}_{E,E_1})\,\iso\, \Laum^{d,d_1}_{E,E_1}$. By Step 1 and Proposition~\ref{Pp_avegaring_commutes}, we get
$$
F_{\tilde H}(\bar\gp_{P\, !}\Laum_{E,E_1}^{d,d_1})\,\iso\, 
(A\det E)_{\Omega^{-1}}\otimes \Av^d_E(\epsilon_!\mult^{0,d_1}_!A(E_1\otimes \det E\otimes\cE_0))
[g-1]
$$
It is easy to check that for any $K\in\D(\Pic^{d_1} X)$ we have 
$$
\Av^d_E(\epsilon_!\mult^{0,d_1}_! K)\,\iso\, \epsilon_!\mult^{d,d_1}_!(K\boxtimes (\pi^* E)^{(d)})[d]
$$
Our assertion follows.   
\end{Prf}

\bigskip\bigskip
\centerline{\scshape 5. The case $H=\GO_4$}

\bigskip\noindent
5.1 Keep the notation of Sect.~3 assuming $m=2$.

\begin{Rem} Given $k$-vector spaces $V_1,V_2$ of dimension 2, we have a canonical symmetric form $\Sym^2(V_1\otimes V_2)\to \det V_1\otimes\det V_2$. One may get a compatible isomorphism
$$
\gamma_{V_1,V_2}: \det(V_1\otimes V_2)\,\iso\, (\det V_1\otimes\det V_2)^2
$$
as follows. Denote by $\St$ (resp., $\det$) the standard (resp., the determinantal) representation of $\GL_2$. Fix an isomorphism $\gamma_{\St}: \det(\St\boxtimes\St)\,\iso\, (\det\boxtimes\det)^2$ of $\GL_2\times\GL_2$-representations compatible with the above symmetric form. It yields the desired isomorphism as follows. Given $V_i$ pick an isomorphism of vector spaces $b_i: V_i\,\iso\,\St$ and define $\gamma_{V_1, V_2}$ by the commutative diagram
$$
\begin{array}{ccc}
\det(V_1\otimes V_2) & \toup{\gamma_{V_1,V_2}} & (\det V_1\otimes\det V_2)^2\\
\downarrow\lefteqn{\scriptstyle b_1\otimes b_2} && \downarrow\lefteqn{\scriptstyle b_1\otimes b_2}\\
\det(\St\boxtimes\St) & \toup{\gamma_{\St}} & (\det\boxtimes\det)^2
\end{array}
$$
Then $\gamma_{V_1, V_2}$ does not depend on $b_i$. We have $\gamma_{V_2,V_1}=-\gamma_{V_1,V_2}$. 
\end{Rem} 

Denote by $\Bun_{k,\tilde X}$ the stack of rank $k$ vector bundles on $\tilde X$. Denote by $\rho: \Bun_{2,\tilde X}\to\Bun_{\tilde H}$ the map sending $W$ to $(V,\cC, \Sym^2 V\toup{h}\cC, \gamma)$, where $V\in\Bun_4$ is the descent of $W\otimes\sigma^*W$ equipped with natural descent data, $\cC=N(\det W)$, and $h$ is the descent of the canonical symmetric form
$\Sym^2(W\otimes\sigma^*W)\to \det W\otimes\sigma^*\det W$. The compatible trivialization 
$$
\det(W\otimes\sigma^*W)\,\iso\,(\det W\otimes\sigma^*\det W)^2
$$ 
descends to $\gamma: \cC^{-2}\otimes\det V\,\iso\,\cE$. The map $\rho$ is smooth and surjective. 

 Another way to spell the same construction is as follows. We have an exact sequence $1\to \Gm\to\GL_2\times\GL_2\to\GO_4^0\to 1$, where the first map sends $x\in\Gm$ to $(x,x^{-1})$. Then we can think of the automorphism $\tilde\sigma$ of $\GO_4^0$ chosen in 3.1 as an automorphism of this exact sequence permuting the two factors of $\GL_2\times\GL_2$. The corresponding twisting of this exact sequence by
the $\Sigma$-torsor $\pi:\tilde X\to X$ gives an exact sequence $1\to U_{\pi}\to \pi_*\GL_2\to \tilde H\to 1$.

 We have $\Bun_{\pi_*\GL_2}\,\iso\,\Bun_{2,\tilde X}$. The stack $\Bun_{U_{\pi}}$ classifies: $\cB\in\Pic\tilde X$ equipped with an isomorphism $N(\cB)\,\iso\,\cO$. The above map $\rho$ is the extension of scalars under $\pi_*\GL_2\to\tilde H$. Write also $\sigma$ for the automorphism of $\Bun_{\tilde H}$ sending $(V,\cC, \Sym^2 V\to \cC,\gamma)$ to $(V,\cC, \Sym^2 V\to \cC, -\gamma)$. Then the following diagram is 2-commutative
$$
\begin{array}{ccccc}
\Bun_{2,\tilde X} & \toup{\rho} & \Bun_{\tilde H}  & \toup{\rho_H} & \Bun_H\\
\downarrow\lefteqn{\scriptstyle \sigma^*} && \downarrow\lefteqn{\scriptstyle \sigma} & \nearrow\lefteqn{\scriptstyle \rho_H}\\
\Bun_{2,\tilde X} & \toup{\rho} & \Bun_{\tilde H}
\end{array}
$$

  Let $\Bun_{2,\tilde X}^d\subset\Bun_{2, \tilde X}$ 
be the substack given by $\deg W=d$. Remind that $\Bun_H^d$ is given by $\deg\cC=d$. For $\tilde X$ connected the irredicibility of $\Bun_{2,\tilde X}^d$ and surjectivity of $\rho$ implies that the stack $\Bun_{\tilde H}^d$ is irreducible, so $\cN$ is a nontrivial local system on each $\Bun_H^d$ in this case. 

 Let $\tilde E$ be an irreducible rank 2 local system on $\tilde X$. Let $\Aut_{\tilde E}$ be the corresponding automorphic sheaf on $\Bun_{2,\tilde X}$ normalized as in \cite{FGV} (cf. also Definition~\ref{Def_Aut_E}). We fix a rank one local system $\chi$ on $X$ and an isomorphism $\pi^*\chi\,\iso\, \det\tilde E$. This provides a descent data for $\Aut_{\tilde E}$ for the map $\Bun_{2,\tilde X}\to \Bun_{\tilde H}$, so we get a perverse sheaf, say $K_{\tilde E, \chi,\tilde H}$ on $\Bun_{\tilde H}$. 
 
  For $\tilde X$ connected the group stack $\Bun_{U_{\pi}}$ has two connected components (cf. A.1), write $\Bun^0_{U_{\pi}}$ for its connected component of unity. 
  
\begin{Def} The quotient of $\Bun_{2,\tilde X}$ by the action of $\Bun^0_{U_{\pi}}$ is a $\mu_2$-torsor over $\Bun_{\tilde H}$, we denote by $\cSN$ the corresponding local system (of order two) on $\Bun_{\tilde H}$. We refer to it as \select{the spinorial norm}. 
\end{Def}

We have $K_{\tilde E,\chi, \tilde H}\otimes\cSN\,\iso\, K_{\tilde E, \chi\otimes\cE_0, \tilde H}$. The central character of $K_{\tilde E,\chi,\tilde H}$ is $\chi$.  
   
 The local system $\pi_*\tilde E^*$ is equipped with a natural symplectic form $\wedge^2(\pi_*\tilde E^*)\to \chi^{-1}$, so gives rise to a $\check{G}$-local system $E_{\check{G}}$ on $X$, where $G=\GSp_4$ for $n=2$. 
 
 If $\tilde X$ splits then we fix a numbering of connected components of $\tilde X$. Then $\tilde E$ becomes a pair of irreducible rank 2 local systems $E_1,E_2$ on $X$. We get $\Bun_{2,\tilde X}\,\iso\,\Bun_2\times\Bun_2$ and $\Aut_{\tilde E}=\Aut_{E_1}\boxtimes\Aut_{E_2}$. The descent datum for $\det\tilde E$ becomes $\det E_1\,\iso\,\det E_2\,\iso\,\chi$. For $\tilde X$ split we have an exact sequence 
$0\to\ZZ/2\ZZ\to \pi_1(\tilde H)\toup{\lambda}\ZZ\to 0$, 
and $\Bun^d_{\tilde H}$ has two connected components $\Bun_{\tilde H}^{\theta}$ for $\theta\in\lambda^{-1}(d)$. 
For $d$ odd the stack $\Bun^d_H$ is connected and the covering $\rho_H:\Bun^d_{\tilde H}\to\Bun^d_H$ splits.  For $d$ even $\Bun^d_H$ has two connected components $\Bun_H^{\theta}$, $\theta\in\lambda^{-1}(d)$, and the covering $\rho_H: \Bun_{\tilde H}^{\theta}\to\Bun_H^{\theta}$ is nontrivial. 
 
  If $E$ is an irreducible rank 2 local system $E$ on $X$ such that $\pi^*E$ is irreducible then the perverse sheaf $K_{\pi^*E, \det E, \tilde H}$ has natural descent data with respect to $\rho_H:\Bun_{\tilde H}\to\Bun_H$, thus defining a perverse sheaf $K_{E, H}$ on $\Bun_H$. Remind the local system $\cE_0$ on $X$ (cf. 3.1). 
  
\begin{Con} 
\label{Con_cusp}
If $n=2$ then we have the following.\\
1) If $\tilde E$ does not descend with respect to $\tilde X\toup{\pi} X$ then $F_G(\rho_{H\, !}  K_{\tilde E, \chi, \tilde H})\in\D(\Bun_G)$ is a cuspidal automorphic sheaf on $\Bun_G$ for $E_{\check{G}}$. (For non connected $\tilde X$ our assumption says that $E_1, E_2$ are non isomorphic irreducible rank 2 local systems on $X$ equipped with $\det E_1\,\iso\,\det E_2$).

\smallskip\noindent
2)  If $E$ is an irreducible rank 2 local system on $X$ with $\pi^*E$ irreducible then we have two cases. If there is an isomorphism $E\,\iso\, E\otimes\cE_0$ on $X$ then $F_G(K_{E, H})\in\D(\Bun_G)$ is isomorphic to the geometric Eisenstein series (for the Siegel parabolic of $\check{G}$), otherwise it is a cuspidal automorphic sheaf on $\Bun_G$ for $E_{\check{G}}$. (In particular, for $\tilde X$ non connected we get an Eisenstein series this way). 
\end{Con}
  
\medskip\noindent
{\bf Question.} 
In case 2) of Conjecture~\ref{Con_cusp} for connected $\tilde X$ what about $F_G(K_{E,H}\otimes\cN)$? 

\medskip\noindent
5.2 In the rest of Sect.~5 we assume in addition $n=1$.
Let $E$ be an irreducible rank 2 local system on $X$. Assume that its restriction $\tilde E=\pi^*E$ is still irreducible. 

 The following is a geometric version of a theorem of Shimizu (\cite{Wa}, Theorem~1), it also an argument supporting Conjecture~\ref{Con_functoriality} in the case $n=1, m=2$. 

\begin{Pp} 
\label{Pp_Shimizu}
For $\tilde X$ split we have $F_{\tilde H}(\Aut_{E^*})\,\iso\, A(\det E)_{\Omega}\otimes K_{\tilde E, \cE_0\otimes\det E, \tilde H}$. This isomorphism depends on a choice of (\ref{iso_gr_k}). 
\end{Pp}
\begin{Prf}
\Step 1 Let $j_d: \Bun_2\times X^{(d)}\hook{} \cS_{\tilde Q}$ be the open immersion sending $(L,D)$ to $L, \cC=\wedge^2 L\otimes\Omega^{-1}(D)$ with canonical inclusion $\wedge^2 L\hook{} \cC\otimes\Omega$. Let $\cS_{\tilde Q}^d$ be the open substack of $\cS_{\tilde Q}$ given by $\deg(\cC\otimes\Omega)-\deg L=d$. We claim that 
$$
F_{\cS_{\tilde Q}}(\Aut_{E^*})\mid_{\cS_{\tilde Q}^d}\,\iso\, j_{d\, !} (\Aut_E\boxtimes E^{(d)})[d]
$$

 Remind the stack $\cW_{G,\tilde Q}$ classifying $M\in\Bun_2, L\in\Bun_2$ and $t:L\to M^*\otimes\Omega$. Let $^0\cW_{G,\tilde Q}\subset \cW_{G,\tilde Q}$ be the open substack given by the condition: $t:L\hook{} M^*\otimes\Omega$ is an inclusion. From cuspidality of $\Aut_{E^*}$ it follows that only $^0\cW_{G,\tilde Q}$ contributes to $F_{\cS_{\tilde Q}}(\Aut_{E^*})\mid_{\cS_{\tilde Q}^d}$, so the latter is extension by zero under $j_d$. 
 
 Let $\epsilon: \Bun_2\to\Bun_2$ be the involution sending $M$ to $M^*\otimes\Omega$. Then $\epsilon^*\Aut_E\,\iso\,\Aut_{E^*}$ canonically. Our assertion follows from Hecke property of $\Aut_E$.  

\Step 2 The map $\nu_{\tilde Q}:\Bun_{\tilde Q}\to\Bun_{\tilde H}$ sends $(L\in\Bun_2, 0\to\wedge^2 L\to V_1\to\cC\to 0)$ to $V=L^*\otimes V_1$ with symmetric form $\Sym^2 V\to\cC$. 
From Step 1 we get
$$
\nu^*_{\tilde Q} F_{\tilde H}(\Aut_{E^*})\,\iso\, A(\det E)_{\Omega}\otimes \nu_{\tilde Q}^*K_{\tilde E,\det E, \tilde H}
$$
 
 There is an open substack $^0\Bun_{\tilde Q}\subset \Bun_{\tilde Q}$ with the following properties. The projection $^0\Bun_{\tilde Q}\to\Bun_{\tilde H}$ is smooth and surjective with connected fibres, and $\nu^*_{\tilde Q}K_{\tilde E,\tilde H}[\dimrel(\nu_{\tilde Q})]$ is a perverse sheaf over $^0\Bun_{\tilde Q}$. Our assertion follows. 
\end{Prf}

\medskip

 Note that Proposition~\ref{Pp_Shimizu} (at least the corresponding non canonical isomorphism) would also follow from Conjecture~\ref{Con_functoriality}. We conjecture that Proposition~\ref{Pp_Shimizu} remains valid for $\tilde X$ nonsplit. 
  
\bigskip\bigskip
\centerline{\scshape 6. Bessel periods for $\GSp_4$}

\bigskip\noindent
6.1.1 Keep the notation of Sect. 5.1. In Sect. 6.1-6.2 we assume $n=2$. Remind the stack $\cS_P$ classifying $L\in\Bun_2,\cA\in\Bun_1$ and $\Sym^2 L\toup{s}\cA\otimes\Omega$. Denote by $\cS_P^r\subset \cS_P$ the open substack given by 
$$
2\deg(\cA\otimes\Omega)-2\deg L=r
$$  

 Let $^{rss}X^{(r)}\subset X^{(r)}$ be the open subscheme classifying divisors $x_1+\ldots+x_r$ on $X$ with $x_i$ pairwise distinct. Let $^{rss}\cS^r_P\subset \cS^r_P$ be the open substack given by the condition that $L\hook{s} L^*\otimes \cA\otimes\Omega$ and $\div(L^*\otimes\cA\otimes\Omega/L)\in {^{rss}X^{(r)}}$. 
Set 
$$
\RCov^r=\Pic X\times_{\Pic X} {^{rss}X^{(r)}},
$$ 
where the map $^{rss}X^{(r)}\to\Pic X$ sends $D$ to $\cO(-D)$, and $\Pic X\to\Pic X$ takes a line bundles to its tensor square. It is understood that $^{rss}X^{(0)}=\Spec k$. 

 We have a map $\gp_1: {^{rss}\cS^r_P}\to {\RCov^r}$ sending the above point to $\cE_{\phi}=(\cA\otimes\Omega)^{-1}\otimes\det L$ with the induced inclusion $\cE_{\phi}^2\hook{}\cO_X$. 
 
\begin{Lm} 
\label{Lm_citation_mine}
{\rm (\cite{Ly2}, 7.7.2)} For $r\ge 0$ the stack $\RCov^r$ classifies two-sheeted coverings $\phi:Y\to X$ ramified exactly at $D_X\in{^{rss}X^{(r)}}$ with $Y$ smooth. The stack $^{rss}\cS^r_P$ identifies with the one classifying collections: $D_X\in {^{rss}X^{(r)}}$, a two-sheeted covering $\phi: Y\to X$ ramified exactly at $D_X$, and $\cB\in\Pic Y$. \QED
\end{Lm} 
 
 The identification in Lemma~\ref{Lm_citation_mine} sends $\cB\in\Pic Y$ to $L=\phi_*\cB$ with symmetric form $\Sym^2 L\toup{s} N(\cB)\,\iso\,\cA\otimes\Omega$. Note that $s$ admits a canonical section $N(\cB)(-D_X)\hook{} \Sym^2 L$, which is a vector subbundle of $\Sym^2 L$. 
Let $\gq_1: {^{rss}\cS^r_P}\to\Bun_2$ be the map sending $(L,\cA,s)$ as above to $L$. 
 
\begin{Lm}
\label{Lm_smothness_curves_vary_1}
 For $r>4(g-1)$ the map $\gq_1: {^{rss}\cS^r_P}\to\Bun_2$ is smooth.
\end{Lm}
\begin{Prf}
Since the projection $^{rss}\cS^r_P\to \RCov^r$ is smooth, the stack $^{rss}\cS^r_P$ is smooth. Since $\Bun_2$ is also smooth, it suffices to show that the fibre of $\gq_1$ over a field-valued point $L\in\Bun_2$ is smooth. 

 Let us calculate the tangent space to the fibre of $\gq_1$ at a point $(L,\cA,s)$. Write for brevity $\cC=\cA\otimes\Omega$. Let $K$ denote the cokernel of $\cO\hook{s}\cC\otimes\Sym^2 L^*$. The sheaf $K$ is locally free. The tangent space in question identifies with $\H^0(X,K)$. We claim that $\H^1(X,K)=0$. 
 
  Indeed, suppose $(L,\cA,s)$ is given by a collection: a two-sheeted covering $\phi: Y\to X$ ramified at $D_X\in {^{rss}X^{(r)}}$ and a line bundle $\cB$ on $Y$.  So, $L\,\iso\,\phi_*\cB$ and $s: \Sym^2 L\to\cC\,\iso\, N(\cB)$ is the natural symmetric form. Let $D_Y\in{^{rss} Y^{(r)}}$ be the ramification divisor of $\phi$, so $D_X=\phi_* D_Y$. 
Then $\Sym^2 L$ is included into a cartesian square
$$
\begin{array}{ccc}
N(\cB)\oplus \phi_*(\cB^2) & \to & \cB^2\mid_{D_Y} \oplus \cB^2\mid_{D_Y} \\
\uparrow && \uparrow\lefteqn{\scriptstyle diag}\\
\Sym^2 L & \to & \cB^2_{D_Y}
\end{array}
$$ 
Let $\sigma_{\phi}$ be the nontrivial automorphism of $Y$ over $X$. We have $K^*\,\iso\, \phi_*(\cB\otimes\sigma_{\phi}^*\cB^{-1}(-D_Y))$. So, 
$\H^0(X, K^*\otimes\Omega)\,\iso\, \H^0(Y, \phi^*\Omega\otimes \cB\otimes\sigma_{\phi}^*\cB^{-1}(-D_Y))=0$, because the degree of the corresponding line bundle on $Y$ is $4(g-1)-r<0$. 

 As $\gq_1$ is separable, and the dimensions of the tangent spaces to the fibres are constant, $\gq_1$ is smooth. 
\end{Prf}

\medskip 
 
  Fix a two-sheeted covering $\phi:Y\to X$ ramified at $D_X\in {^{rss}X^{(r)}}$. Write $\sigma_{\phi}$ for the nontrivial automorphism of $Y$ over $X$, and $\cE_{\phi}$ for the $\sigma_{\phi}$-anti-invariants in $\phi_*\cO_Y$, it is equipped with $\cE_{\phi}^2\,\iso\, \cO(-D_X)$. Write $D_Y\in Y^{(r)}$ for the ramification divisor of $\phi$, so $D_X=\phi_*D_Y$.  Remind that $\phi$ is recovered from $(\cE_{\phi}, D_X)$ as $\Spec (\cO_X\oplus \cE_{\phi})$, where the structure of a $\cO_X$-algebra on $\cO_X\oplus \cE_{\phi}$ is given by $\cE_{\phi}^2\hook{}\cO_X$. 
  
  Remind the stack $\Bun_{U_{\phi}}$, its connected components $\Bun_{U_{\phi}}^a$ are indexed by $a\in\ZZ/2\ZZ$ (cf. A.1). Let $^0\Bun_{U_{\phi}}\subset \Bun_{U_{\phi}}$ be the open substack given by $\H^0(Y, \cV\otimes \phi^*\Omega)=0$ for $\cV\in\Bun_{U_{\phi}}$ equipped with $N(\cV)\,\iso\,\cO_X$. 
  
   Let $^0\Pic Y$ be the preimage of $^0\Bun_{U_{\phi}}$ under $e_{\phi}: \Pic Y\to\Bun_{U_{\phi}}$ (cf. A.1). Denote by $\phi_1: \Pic Y\to\Bun_2$ the map sending $\cB$ to $\phi_*\cB$. 

 For $g=0$ we have $^0\Bun_{U_{\phi}}=\Bun_{U_{\phi}}$. 
If $g=1$ then $^0\Bun_{U_{\phi}}\subset\Bun_{U_{\phi}}$ is given by the condition that $\cV$ is not isomorphic to $\cO_Y$.   

\begin{Lm}
i) If $r\ge 4g-4$ then $^0\Bun_{U_{\phi}}^a$ is nonempty for each $a\in\ZZ/2\ZZ$. So, the intersection of $^0\Pic Y$ with each connected component of $\Pic Y$ is nonempty. 

\smallskip\noindent
ii) The restriction of $\phi_1:\Pic Y\to\Bun_2$ to the open substack $^0\Pic Y\subset\Pic Y$ is smooth.
\end{Lm}  
\begin{Prf}
i) Write $\Ker \und{N}$ for the kernel of the norm map $\und{N}: \uPic Y\to \uPic X$ (cf. A.1). Let $^0\Ker\und{N}$ be the open subscheme given by $\H^0(Y, \cV\otimes \phi^*\Omega)=0$ for $\cV\in \Ker\und{N}$. Then $^0\Bun_{U_{\phi}}$ is the preimage of $^0\Ker\und{N}$ under the projection $\Bun_{U_{\phi}}\to \Ker\und{N}$.
 
 Let $Z$ denote the preimage of $\Omega^2$ under the map $X^{(4g-4)}\to \uPic X$ sending $D$ to $\cO(D)$. Here $\uPic X$ is the Picard scheme of $X$. Let $Z'$ be the preimage of $Z$ under $\phi: Y^{(4g-4)}\to X^{(4g-4)}$. We have $Z'=\emptyset$ for $g=0$, $Z'=\Spec k$ for $g=1$, and $\dim Z'=3g-4$ for $g>1$. Then $^0\Ker\und{N}$ is the complement to the image of the map $Z'\to \Ker\und{N}$ sending $D$ to $(\phi^*\Omega^{-1})(D)$. Since each connected component of $\Ker\und{N}$ is of dimension $g-1+\frac{r}{2}$, our assertion follows. 
 
\smallskip\noindent
ii) Since both $^0\Pic Y$ and $\Bun_2$ are smooth, it suffices to check that for $\cB\in {^0\Pic Y}$ the natural map $\H^1(Y,\cO)\to \H^1(X, \END(\phi_*\cB))\,\iso\, \H^1(Y, \cB\otimes\phi^*((\phi_*\cB)^*))$ is surjective. We have a cartesian square
$$
\begin{array}{ccc}
\cO(D_Y)\oplus \frac{\cB}{\sigma_{\phi}^*\cB}(D_Y) & \to & \cO(D_Y)/\cO\oplus \cO(D_Y)/\cO\\
\uparrow && \uparrow\\
\cB\otimes\phi^*((\phi_*\cB)^*) & \to & \cO(D_Y)/\cO,
\end{array}
$$
where the right vertical arrow is the diagonal map. This yields an exact sequence 
$$
0\to \cO_Y\to \cB\otimes\phi^*((\phi_*\cB)^*)\to  \frac{\cB}{\sigma_{\phi}^*\cB}(D_Y)\to 0
$$ 
We have $\H^1(Y, \cB(D_Y)\otimes\sigma_{\phi}^*\cB^{-1})=\H^0(Y, \phi^*\Omega\otimes \sigma_{\phi}^*\cB\otimes\cB^{-1})^*=
0$, because $\frac{\sigma_{\phi}^*\cB}{\cB}\in {^0\Bun_{U_{\phi}}}$. We are done.
\end{Prf}

\medskip  
  
 Our purpose is to study the $*$-restriction of $F_{\cS}(\rho_{H\, !}K_{\tilde E,\chi, \tilde H})$ under $\Pic Y\hook{} {^{rss}\cS^r_P}\subset \cS_P$.   
For $d\in\ZZ$ set $\bar d=2\deg\Omega-d$. The complex
$
F_{\cS}(\rho_{H\, !}K_{\tilde E,\chi, \tilde H}\mid_{\Bun_H^d})
$ 
lies in $\D(\Pic^{\bar d} Y)$. 

Define $\tilde Y$ by the cartesian square
\begin{equation}
\label{square_def_tilde_Y}
\begin{array}{ccc}
\tilde Y & \toup{\tilde\pi} & Y\\ 
\downarrow\lefteqn{\scriptstyle \tilde\phi} && \downarrow\lefteqn{\scriptstyle \phi}\\
\tilde X & \toup{\pi} & X
\end{array}
\end{equation}

 Consider the exact sequence $1\to\Gm\to\Gm\times\Gm\toup{N}\Gm\to 1$, where $N$ is the product map, and the first map sends $z$ to $(z,z^{-1})$. Twisting it by the $\ZZ/2\ZZ$-torsor $\tilde X\toup{\pi} X$ via the action permuting two factors of $\Gm\times\Gm$, we get an exact sequence $1\to U_{\pi}\to\pi_*\Gm\toup{N}\Gm\to 1$ of group schemes on $X$, here $N$ is the norm map. Consider the composition $U_{\pi}\to \pi_*\Gm\to\pi_*\tilde\phi_*\Gm\,\iso\, \phi_*\tilde\pi_*\Gm$. Define the group scheme $R_{\phi}$ on $X$ by the exact sequence 
\begin{equation}
\label{seq_R_phi}
1\to U_{\pi}\to \phi_*\tilde\pi_*\Gm\to R_{\phi}\to 1
\end{equation}
The corresponding map $\Pic\tilde Y\to\Bun_{R_{\phi}}$ is smooth and surjective. 

 Let $_{\phi}\GL_2$ be the group scheme of automorphisms of $\phi_*\cO_Y$, this is an inner form of $\GL_2$. We denote by the same symbol $_{\phi}\GL_2$ its restriction to $\tilde X$. We have a natural map $\tilde\phi_*\Gm\to {_\phi\GL_2}$ of group schemes on $\tilde X$. Let $_{\phi}\tilde H$ be the group scheme on $X$ included into 
a morphism of exact sequences 
$$
\begin{array}{ccccccc}
1\to  & U_{\pi} & \to & \pi_*\tilde\phi_*\Gm & \to & R_{\phi} & \to 1\\
 & \downarrow\lefteqn{\scriptstyle\id} && \downarrow && \downarrow\\
1\to & U_{\pi} & \to & \pi_*(_{\phi}\GL_2) & \to & _{\phi}\tilde H &\to 1
\end{array}
$$
Since $\Bun_{_{\phi}\tilde H}\,\iso\,\Bun_{\tilde H}$ canonically, we get a morphism $\gq_{R_{\phi}}: \Bun_{R_{\phi}}\to\Bun_{\tilde H}$.  

 For $\cL\in \Pic\tilde X$ we have $\phi^*N(\cL)\,\iso\, N(\tilde\phi^*\cL)$ canonically.  Consider the map $\phi_*\tilde\pi_*\Gm\toup{\phi_*N} \phi_*\Gm$ induced by the norm $\tilde\pi_*\Gm\toup{N}\Gm$. It is easy to check that $U_{\pi}\subset \Ker \phi_*N$, so we get a map $R_{\phi}\to \phi_*\Gm$. Let $\gp_{R_{\phi}}$ denote the composition of the extension of scalars $\Bun_{R_{\phi}}\to\Pic Y$ with the automorphism $\epsilon:\Pic Y\iso\Pic Y$ sending $\cB$ to $\cB^{\star}=\cB^*\otimes\Omega_Y$. So, the diagram commutes
\begin{equation}
\label{diag_phi}
\begin{array}{cccc}
\Pic \tilde Y & \to &\Bun_{R_{\phi}}& \toup{\gq_{R_{\phi}}} \Bun_{\tilde H}\\
\downarrow\lefteqn{\scriptstyle N} && \downarrow\lefteqn{\scriptstyle\gp_{R_{\phi}}}\\ 
\Pic Y & \toup{\epsilon} & \Pic Y
\end{array}
\end{equation}
  
  When $(\cE_{\phi}, D_X)$ run through $\RCov^r$, the group schemes $R_{\phi}$ are naturally organized into a group scheme $R$ over $\RCov^r\times X$. Let $\Bun_R$ denote the stack over $\RCov^r$ associating to a scheme $S$ the category: a map $S\to\RCov^r$, and a $R_S$-torsor on $S\times X$, where $R_S$ is the restriction of $R$ under $S\times X\to \RCov^r\times X$. 
 
 Diagrams (\ref{diag_phi}) naturally form a family 
$$
\Bun_{\tilde H} \;\getsup{\gq_R}\; \Bun_R\; \toup{\gp_R} \;{^{rss}\cS^r_P}
$$ 

\begin{Pp}
\label{Pp_R}
The restriction of $F_{\cS}(\rho_{H\, !}K_{\tilde E,\chi, \tilde H})$ to the open substack ${^{rss}\cS^r_P}\subset \cS_P$ is canonically isomorphic to 
$$
(\gp_R)_!\, \gq_R^*K_{\tilde E,\chi, \tilde H}[\dimrel(\gq_R)]
$$
In particular, for a $k$-point of $\RCov^r$ given by $\phi:Y\to X$, the $*$-restriction identifies canonically 
$$
F_{\cS}(\rho_{H\, !}K_{\tilde E,\chi, \tilde H})\mid_{\Pic Y}\,\iso\, 
(\gp_{R_{\phi}})_!\, \gq_{R_{\phi}}^*K_{\tilde E,\chi,\tilde H}[\dimrel(\gq_R)]
$$
\end{Pp}   
\begin{Prf}   
Define a map $\zeta_{\phi}:\Pic\tilde Y\to \cV_{\tilde H,P}$ as follows. Given $\cB\in\Pic \tilde Y$, let $W=\tilde\phi_*\cB$ and $(V,\cC,\Sym^2 V\to\cC,\gamma)\in\Bun_{\tilde H}$ be the corresponding $\tilde H$-torsor. We have $\sigma^*\tilde\phi_*\cB\,\iso\, \tilde\phi_*\sigma^*\cB$, so there is a natural map 
$$
\pi^*V\,\iso\, W\otimes\sigma^*W\to \tilde\phi_*\tilde\pi^*N_Y(\cB)\,\iso\, \pi^*\phi_*N_Y(\cB),
$$ 
where $N_Y:\Pic\tilde Y\to\Pic Y$ is the norm map. It  descends to a map $V\to  
\phi_*N_Y(\cB)$. So, for $L=\phi_*(N_Y(\cB)^{\star})$ we get a map $t: V\to L^*\otimes\Omega$. By definition, $\zeta_{\phi}$ sends $\cB$ to $(V,\cC,\Sym^2 V\to\cC,\gamma)\in\Bun_{\tilde H}$, $L\in\Bun_2$ and $t:V\to L^*\otimes\Omega$.
 
 Denote by $\tilde\phi_1:\Pic \tilde Y\to\Bun_{2,\tilde X}$ the map sending $\cB$ to $\tilde\phi_*\cB$. We have a commutative diagram
$$
\begin{array}{ccccc}
\Bun_{2,\tilde X} & \getsup{\tilde\phi_1} & \Pic\tilde Y & \toup{\epsilon\comp N} & \Pic Y \\
\downarrow && \downarrow\lefteqn{\scriptstyle \zeta_{\phi}} && \downarrow \\
 \Bun_{\tilde H} & \gets & \cV_{\tilde H,P} & \to  &\cS_P
\end{array}
$$
It extends naturally to a diagram
$$
\begin{array}{ccccc}
&& \Pic\tilde Y \\
& \swarrow\lefteqn{\scriptstyle\zeta_{\phi}} & \downarrow & \searrow\lefteqn{\scriptstyle \epsilon\comp N}\\
\cV_{\tilde H,P} & \gets & \Bun_{R_{\phi}} & \toup{\gp_{R_{\phi}}} &\Pic Y\\
& \searrow &&  \swarrow\\
&& \cS_P 
\end{array}
$$
As $\phi$ varies in $\RCov^r$ these diagrams form a family
\begin{equation}
\label{diag_I_realized}
\begin{array}{ccccc}
&& \Bun_R & \toup{\gp_R} & {^{rss}\cS^r_P}\\
& \swarrow\lefteqn{\scriptstyle \gq_R} &
\downarrow && \downarrow\\
\Bun_{\tilde H} & \gets & \cV_{\tilde H,P} &\to &\cS_P
\end{array}
\end{equation}
Our assertion is reduced to the following lemma.
\end{Prf} 
 
\begin{Lm} 
\label{Lm_cartesian_R}
The square in (\ref{diag_I_realized}) is cartesian.  
\end{Lm} 
 
Our proof of Lemma~\ref{Lm_cartesian_R} uses the following elementary observation.

\begin{Slm} 
\label{Slm_nice}
Let $K$ be a field of characteristic different from $2$. Let 
$V_i$ be 2-dimensional $K$-vector spaces. Let $\cB, \cB'$ be 1-dimensional $K$-vector spaces. Equip $\cB\oplus\cB'$ with the quadratic form 
$$
s:\Sym^2(\cB\oplus\cB')\,\iso\,\cB^2\oplus\cB'^2\oplus\cB\otimes\cB'\toup{(0,0,\id)} \cB\otimes\cB'
$$
Assume given a map $t: \cB\oplus\cB'\to V_1\otimes V_2$ such that there is a commutative diagram
$$
\begin{array}{ccc}
\Sym^2(\cB\oplus\cB') & \toup{t\otimes t} & \Sym^2(V_1\otimes V_2)\\
\downarrow\lefteqn{\scriptstyle s} && \downarrow\\
\cB\otimes\cB' & \iso & \det V_1\otimes\det V_2
\end{array}
$$
Then there exist  a unique decomposition into a direct sum of 1-dimensional subspaces $V_i\,\iso\, U_i\oplus U'_i$ and unique isomorphisms $\cB\,\iso\, U_1\otimes U_2$, $\;\cB'\,\iso\, U'_1\otimes U'_2$ under which $t$ identifies with a natural inclusion
$$
(U_1\otimes U_2)\oplus (U'_1\otimes U'_2)\hook{} V_1\otimes V_2
\eqno{\square}
$$
\end{Slm}  

\begin{Prf}\select{of Lemma~\ref{Lm_cartesian_R}} \\
Consider a point of $\cV_{\tilde H,P}$ given by $\cF_{\tilde H}=(V,\cC, \Sym^2 V\to\cC, \gamma)\in\Bun_{\tilde H}$, $L\in\Bun_2$ and $t: V\to L^*\otimes\Omega$.  
Assume that its image in $\cS_P$ is identified with a point 
$(\phi:Y\to X, \, \cB\in\Pic Y)$ of $^{rss}\cS^r_P$. So, we are given an isomorphism $L\,\iso\,\phi_*\cB$ and the diagram commutes
\begin{equation}
\label{diag_square_V_L}
\begin{array}{cccc}
\Sym^2 V & \toup{t\otimes t} & \Sym^2(L^*\otimes\Omega)\\ 
\uparrow && \uparrow\\
\cC & \iso & N(\cB)^{-1}\otimes\Omega^2 & \iso\; N(\cB^{\star})(-D_X)
\end{array}
\end{equation}

 Since $L^*\otimes\Omega\,\iso\,\phi_*\cB^{\star}$, we view the datum of $t$ as $t:\phi^*V\to\cB^{\star}$. We have a commutative diagram
\begin{equation}
\label{diag_t_sigmat}
\begin{array}{ccc}
\phi^* V & \toup{(t,\sigma_{\phi}^*t)} & \cB^{\star}\oplus\sigma_{\phi}^*\cB^{\star},\\  
\downarrow & \nearrow\\ 
\phi^*\phi_*\cB^{\star}
\end{array}
\end{equation}
where  
$$
\phi^*\phi_*\cB^{\star}=\{b\in \cB^{\star}\oplus \sigma_{\phi}^*\cB^{\star}\mid 
\mbox{the image of}\;
b\; \mbox{in}\; (\cB^{\star}\oplus \sigma_{\phi}^*\cB^{\star})\mid_{D_Y}\; \mbox{lies in the diagonal}\; \cB^{\star}\mid_{D_Y}\}
$$
Pick a lifting of $\cF_{\tilde H}\in\Bun_{\tilde H}$ to a point $W\in\Bun_{2,\tilde X}$. We get $\pi^*V\,\iso\, W\otimes\sigma^* W$. 
From (\ref{diag_t_sigmat}) we get a map
$$
\tilde\phi^*(W\otimes\sigma^* W)\; \toup{(t,\sigma_{\phi}^*t)}\;
\tilde\pi^*(\cB^{\star}\oplus\sigma_{\phi}^*\cB^{\star})
$$
whose tensor square fits into a commutative diagram
$$
\begin{array}{ccc}
\tilde\phi^*\Sym^2(W\otimes\sigma^* W) & \to & \tilde\pi^*(\cB^{\star 2} \oplus \sigma_{\phi}^*\cB^{\star 2}\oplus \cB^{\star}\otimes\sigma_{\phi}^*\cB^{\star})\\
\uparrow && \uparrow\lefteqn{\scriptstyle (0,0,1)}\\
\tilde\phi^*(\det W\otimes\sigma^*\det W) & \iso & \tilde\pi^*(\cB^{\star}\otimes\sigma_{\phi}^*\cB^{\star}(-2D_Y))
\end{array}
$$
By abuse of notation, we also write $\sigma$ and $\sigma_{\phi}$ for the involutions of $\tilde Y$ obtained by base change in the square (\ref{square_def_tilde_Y}). 

 Note that any surjection $\tilde\phi^*W\to \cL$, where $\cL$ is a line bundle on $\tilde Y$, gives rise to a map $\xi_{\cL}: V\to \phi_*N_Y(\cL)$. Indeed, the composition 
$$
\pi^*V\,\iso\, W\otimes\sigma^* W\to\tilde\phi_*\cL\otimes\sigma^*\tilde\phi_* \cL\to \tilde\phi_*(\cL\otimes\sigma^* \cL)\,\iso\, \tilde\phi_*\tilde\pi^*N_Y(\cL)
$$
descends to a map $\xi_{\cL}: V\to \phi_*N_Y(\cL)$. 

 By Sublemma~\ref{Slm_nice}, there is a unique rank 1 subbundle $W_1\subset \tilde\phi^*W$, for which we set $\cL=(\tilde\phi^*W)/W_1$, and a unique $\sigma$-invariant inclusion of coherent sheaves $\cL\otimes\sigma^*\cL\hook{} \tilde\pi^*\cB^{\star}$ with the following properties. The latter inclusion gives rise to an inclusion $N_Y(\cL)\hook{}\cB^{\star}$, and the composition
\begin{equation}
\label{eq_composition_in_Lm}
V\toup{\xi_{\cL}} \phi_*N_Y(\cL)\to \phi_*\cB^{\star}
\end{equation}
equals $t$. 

 Taking symmetric squares in (\ref{eq_composition_in_Lm}), we get a commutative diagram
\begin{equation}
\label{diag_one_more_C_maps}
\begin{array}{ccccccc}
\Sym^2 V & \to & \Sym^2(\phi_*N_Y(\cL)) & \hook{} & \Sym^2(\phi_*\cB^{\star}) & \to & N(\cB^{\star}),\\ 
\uparrow && \uparrow && \uparrow & \nearrow\\
 \cC && (N_XN_Y(\cL))(-D_X) & \hook{} & N_X(\cB^{\star})(-D_X)
\end{array}
\end{equation}
in which the middle square is cartesian (and all the three vertical arrows are subbundles). Using (\ref{diag_square_V_L}) we conclude that there is a unique isomorphism $\eta: \cC\,\iso\, (N_XN_Y(\cL))(-D_X)$ making (\ref{diag_one_more_C_maps}) commute, and the inclusion $N_Y(\cL)\hook{}\cB^{\star}$ is actually an isomorphism.  

 Let us show that the natural map $W\to\tilde\phi_*\cL$ is an isomorphism. 
We have inclusions
$$
\tilde\phi^*W\hook{} \tilde\phi^*\tilde\phi_*\cL\hook{} \cL\oplus \sigma_{\phi}^*\cL,
$$
whose determinants yield $\tilde\phi^*\det W\hook{} 
\det(\tilde\phi^*\tilde\phi_*\cL)\,\iso\, (\cL\otimes \sigma_{\phi}^*\cL)(-\tilde\pi^*D_Y)$. Symmetrizing with respect to the action of $\sigma$, one gets inclusions
\begin{multline*}
\tilde\phi^*\pi^*\cC\;\iso\; \tilde\phi^*(\det W\otimes\sigma^*\det W)\hook{} (\cL\otimes\sigma^* \cL)\otimes \sigma_{\phi}^*(\cL\otimes\sigma^* \cL)(-2\tilde\pi^*D_Y)\,\iso\\ (\tilde\pi^*\phi^*N_XN_Y(\cL))(-2\tilde\pi^*D_Y)
\end{multline*}
whose composition is an isomorphism (equal to restriction of $\eta$). So, $W\,\iso\, \tilde\phi_*\cL$ is an isomorphism. 

 Viewing $\cL$ as a $\phi_*\tilde\pi_*\Gm$-torsor on $X$, let $\cF_{R_{\phi}}$ be the $R_{\phi}$-torsor on $X$ obtained from it by extension of scalars (\ref{seq_R_phi}). Then $\gq_{R_{\phi}}(\cF_{R_{\phi}})\,\iso\, \cF_{\tilde H}$ equips $\cF_{\tilde H}$ with a $R_{\phi}$-structure that does not depend on a choice of a lifting of $\cF_{\tilde H}$ to a $\pi_*\GL_2$-torsor. We are done. 
\end{Prf} 
 
\medskip\noindent
\begin{Rem}
\label{Rem_sum_two_components}
Consider the case of $\tilde X\toup{\pi} X$ split. 
We have an exact sequence $0\to \ZZ/2\ZZ\to \pi_1(\tilde H)\to\ZZ\to 0$. The $*$-restriction
$$
F_{\cS}(\rho_{H\, !}K_{\tilde E, \chi, \tilde H})\mid_{\Pic^{\bar d} Y}
$$
is naturally a direct sum of two complexes indexed by $\theta\in\pi_1(\tilde H)$ whose image in $\ZZ$ is $d$. 
If $Y$ is connected then $\gq_{R_{\phi}}: \Bun_{R_{\pi}}\to\Bun_{\tilde H}$ induces a bijection at the level of connected components $\pi_0(\Bun_{R_{\phi}})\,\iso\, \pi_0(\Bun_{\tilde H})\,\iso\,\pi_1(\tilde H)$.  
\end{Rem} 
  
\medskip\noindent
6.1.2 Remind the diagram
$$
\Bun_2 \getsup{\gq_1} {^{rss}\cS^r_P}\toup{\gp_1}\RCov^r
$$
introduced in Sect.~6.1.1. In this subsection we prove the following acyclicity result.

\begin{Th} 
\label{Th_ULA}
Let $E$ be a rank 2 local system on $X$, $K\in\D(\Bun_2)$ be a Hecke eigensheaf with eigenvalue $E$. Then $\gq_1^*K$ is ULA with respect to $\gp_1: {^{rss}\cS^r_P}\toup{\gp_1}\RCov^r$.
\end{Th}
\begin{Prf}
\Step 1 The difficulty comes from the fact that $\gq_1\times\gp_1: {^{rss}\cS^r_P}\to \Bun_2\times\RCov^r$ is not smooth (for $g\ge 1$), we come around it using the Hecke property of $K$. Namely, for $d\ge 0$ consider the diagram 
$$
X^{(d)}\times\Bun_2\getsup{\supp\times p'}\Mod_2^d\toup{p} \Bun_2,
$$
where $\Mod_2^d$ is the stack classifying a lower modification $(L\subset L')$ of rank 2 vector bundles on $X$ with $\deg(L'/L)=d$, the map $p$ (resp., $p'$) sends $(L\subset L')$ to $L$ (resp., $L'$). The map $\supp$ sends this point to $\div(L'/L)$. As in (\cite{FGV}, Sect.~9.5) one shows that 
$$
(\supp\times p')_!p^*K[d]\,\iso\, (E^*)^{(d)}\boxtimes K,
$$
this is the only property of $K$ that we actually use. 

 Define temporary the stack $\cX$ and the maps $p_{\cX}, p'_{\cX}$ by the diagram, where the square is cartesian
$$
\begin{array}{ccc}
&& \Bun_2\times\RCov^r\\
&\nearrow\lefteqn{\scriptstyle p_{\cX}} & \uparrow\lefteqn{\scriptstyle p\times\id}\\
\cX & \to & \Mod_2^d\times\RCov^r\\ 
\downarrow\lefteqn{\scriptstyle p'_{\cX}} &&\downarrow\lefteqn{\scriptstyle \supp\times p'\times\id}\\
X^{(d)}\times{^{rss}\cS^r_P} & \toup{\id\times \gq_1\times\gp_1} &
X^{(d)}\times \Bun_2\times\RCov^r
\end{array}
$$
The above property of $K$ yields an isomorphism
\begin{equation}
\label{iso_for_acyclicity_Pp} 
(p'_{\cX})_!p_{\cX}^*(K\boxtimes\Qlb)[d]\,\iso\, (E^*)^{(d)}\boxtimes \gq_1^*K
\end{equation}
 
\Step 2  Let us show that for $d>4g-4$ the map $p_{\cX}:\cX\to \Bun_2\times\RCov^r$ is smooth. 
The projection $\cX\to {^{rss}\cS^r_P}$ is smooth of relative dimension $2d$. Since $^{rss}\cS^r_P$ is smooth, $\cX$ is also smooth. 

 Fix a $k$-point of $\RCov^r$ given by a two-sheeted covering $\phi:Y\to X$. The corresponding objects $D_Y$, $\cE_{\phi}$ and $\sigma_{\phi}$ are as in 6.1.1.
Let $p_{\cX,\phi}:\cX_{\phi}\to\Bun_2$ be obtained from $p_{\cX}$ by the base change $\Bun_2\times\Spec k\to \Bun_2\times\RCov^r$. The stack $\cX_{\phi}$ classifies: $L\in\Bun_2$, $\cB\in\Pic Y$ with an inclusion of coherent sheaves $L\subset \phi_*\cB$ such that $\div((\phi_*\cB)/L)$ is of degree $d$. It suffices to show that $p_{\cX,\phi}$ is smooth.
 
 Write $\Bunb_{B,Y}$ for the stack classifying: $V_1\in\Pic Y$, $V\in\Bun_{2,Y}$ and an inclusion of coherent sheaves $V_1\subset V$. Let $^0\Bunb_{B,Y}$ be the open substack given by $\H^1(Y, V_1^*\otimes(V/V_1))=0$. One checks that the projection $^0\Bunb_{B,Y}\to\Bun_{2,Y}$ is smooth. Set 
$$
\cY=\Bun_2\times_{\Bun_{2,Y}} {^0\Bunb_{B,Y}},
$$ 
where the map $\Bun_2\to \Bun_{2,Y}$ sends $L$ to $\phi^*L^*$. So, the projection $\cY\to\Bun_2$ is smooth. 
 
  We have an open immersion $j: \cX_{\phi}\hook{} \Bun_2\times_{\Bun_{2,Y}}\Bunb_{B,Y}$ sending $(L\subset \phi_*\cB)$ to $L\in\Bun_2$, $V_1=\cB^{-1}$ with the induced inclusion $\cB^{-1}\hook{}\phi^*L^*$. It suffices to show that the image of $j$ is contained in $\cY$. 
  
 Let $(L\subset \phi_*\cB)$ be a $k$-point of $\cX_{\phi}$ with $D=\div((\phi_*\cB)/L)$. Note that $(\det L)(D)\,\iso\, \cE_{\phi}\otimes N(\cB)$. Define an effective divisor $D'$ on $Y$ and $L_1\in\Pic Y$ by the exact sequence 
$$
0\to L_1\to \phi^*L\to \cB(-D')\to 0
$$ 
Then $L\subset \phi_*(\cB(-D'))$, and taking the determinants we get
$\det L\subset \cE_{\phi}\otimes N(\cB)(-\phi_* D')$, 
so $D\ge \phi_*D'$. We must show that $\H^1(Y, \cB\otimes L_1^*)=0$. 

 We have $L_1\otimes \cB(-D')\,\iso\, \phi^*\det L\,\iso\, \cB\otimes\sigma_{\phi}^*\cB(-D_Y-\phi^*D)$, because $\phi^*\cE_{\phi}\,\iso\, \cO(-D_Y)$. Our assertion follows from the fact that
$$
\cB^*\otimes L_1\otimes\Omega_Y\,\iso\, 
\frac{\sigma_{\phi}^*\cB}{\cB}\otimes\phi^*\Omega(D'-\phi^*D)
$$ 
is of degree $4g-4-2d+\deg D'\le 4g-4-d<0$. 

\medskip

\Step 3 Assume $d>4g-4$. Since $K\boxtimes\Qlb$ is ULA with respect to the projection $\Bun_2\times\RCov^r\to\RCov^r$, it follows that $p_{\cX}^*(K\boxtimes\Qlb)$ is ULA over $\RCov^r$.
Since $p'_{\cX}$ is proper, 
$$
(p'_{\cX})_!p_{\cX}^*(K\boxtimes\Qlb)
$$ 
is ULA over $\RCov^r$ (cf. \cite{BG}, Sect.~5.1.2). 
 Using (\ref{iso_for_acyclicity_Pp}), we learn that $(E^*)^{(d)}\boxtimes \gq_1^*K$ is ULA over $\RCov^r$. 
So, the restriction of $(E^*)^{(d)}\boxtimes \gq_1^*K$ to
$^{rss}X^{(d)}\times{^{rss}\cS^r_P}$
is also ULA over $\RCov^r$. Since $E^{(d)}$ is a local system over $^{rss}X^{(d)}$, and the ULA property is local in the smooth topology of the source, our assertion follows.
\end{Prf}  

\medskip
  
 Consider a $k$-point of $\RCov^r$ given by 
$\phi:Y\to X$ as in 6.1.1. Remind that $\phi_1:\Pic Y\to\Bun_2$ sends $\cB$ to $\phi_*\cB$. 
Using (\cite{BG}, Property 4 of Sect.~5.1.2) and Lemma~\ref{Lm_smothness_curves_vary_1} we obtain  
 
\begin{Cor} Let $E$ be a rank 2 local system on $X$, $K$ be a $E$-Hecke eigensheaf on $\Bun_2$. Then 
$$
\DD(\phi_1^*K[\dimrel(\phi_1)])\,\iso\, \phi_1^*\DD(K)[\dimrel(\phi_1)]
$$ 
Besides, if $K$ is perverse then for $r>4g-4$ the sheaf $\phi_1^*K[\dimrel(\phi_1)]$ is perverse.
\QED
\end{Cor}
 
\begin{Rem} 
\label{Rem_ULA_and_C_W}
i) Let $E$ be an irreducible rank 2 local system on $X$, $\Aut_E$ be the corresponding automorphic sheaf on $\Bun_2$ normalized as in \cite{FGV}. Then for any $r$ the complex $\phi_1^*\Aut_E$ is a direct sum of (possibly shifted) perverse sheaves. Indeed, take $d>4g-4$ and apply the decomposition theorem for the (shifted) perverse sheaf $p_{\cX,\phi}^*\Aut_E$ and 
the proper map $p'_{\cX,\phi}: \cX_{\phi}\to X^{(d)}\times\Pic Y$ as in the proof of Theorem~\ref{Th_ULA}. 

\smallskip\noindent
ii) The map $\phi_1:\Pic Y\to\Bun_2$ is not flat, because its fibres have different dimensions (this is related to the fact that the dimension of the scheme of automorphisms of $L\in\Bun_2$ varies). 
\end{Rem} 
 
\smallskip  
\medskip\noindent
6.2.1 We may view $\sigma$ and $\sigma_{\phi}$ as (commuting) automorphisms of $\tilde Y$. Let $Z$ be the quotient of $\tilde Y$ by the involution $\sigma\comp\sigma_{\phi}$, so we get two-sheeted coverings $\tilde Y\toup{\alpha} Z\toup{\beta} X$. Note that $Z$ is smooth, and $\alpha$ is unramified. Let $D_Z$ be the ramification divisor of $\beta$, then $\beta_*D_Z=D_X$. Let $\sigma_{\beta}$ be the nontrivial automorphism of $Z$ over $X$.

 Another way is to say that we let $\cE_{\beta}=\cE\otimes\cE_{\phi}$, it is equipped with $\cE_{\beta}^2\,\iso\, \cO(-D_X)$. Then $Z\,\iso\,\Spec (\cO_X\oplus\cE_{\beta})$, the structure of an $\cO_X$-algebra on $(\cO_X\oplus\cE_{\beta})$ is given by $\cE_{\beta}^2\hook{} \cO_X$. 
Let $\cE_{0,\phi}$ be the $\sigma_{\phi}$-anti-invariants in $\phi_*\Qlb$. 
Then we have $\beta_*\Qlb\,\iso\, \Qlb\oplus \cE_{0,\beta}$ with 
$$
\cE_{0,\beta}\,\iso\, \cE_{0,\phi}\otimes\cE_0
$$    
 
 Let $\Pic (Y,Z)$ be the stack classifying $\cB_1\in\Pic Y, \cB_2\in\Pic Z$, an isomorphism of line bundles $N(\cB_1)\,\iso\, N(\cB_2)$ on $X$, and its refinement $\gamma_{12}: \cB_1\mid_{D_Y}\iso\; \cB_2\mid_{D_Z}$ over $D_Y\,\iso\, D_Z$. This means that $\gamma_{12}^2$ coincides with 
$$
N(\cB_1)\mid_{D_X}\,\iso\, N(\cB_2)\mid_{D_X}
$$ 
We have used the fact that $\beta$ and $\phi$ yield isomorphisms of (reduced) schemes $D_Z\,\iso\, D_X\,\iso\, D_Y$.  
 
\begin{Lm} 
\label{Lm_structure_of_R_phi}
The map $\Pic\tilde Y\to \Pic (Y,Z)$ sending $\cB$ to $(N_Y(\cB), N_Z(\cB))$ yields  an isomorphism $\Bun_{R_{\phi}}\,\iso\, \Pic (Y,Z)$.
\end{Lm}
\begin{Prf}
Denote by $\tilde R_{\phi}$ the preimage of the diagonal $\Gm\hook{}\Gm\times\Gm$ under the homomorphism 
$$
\phi_*\Gm\times\beta_*\Gm\toup{N\times N} \Gm\times\Gm
$$ 
of group schemes on $X$. The product of norms yields a homomorphism $\phi_*\tilde\pi_*\Gm\to \tilde R_{\phi}$ of group schemes on $X$, and $U_{\pi}$ lies in its kernel. The induced map $R_{\phi}\to\tilde R_{\phi}$ is an isomorphism over $X-D_X$, but not everywhere (if $\phi$ is ramified). The group scheme $\tilde R_{\phi}\mid_{D_X}$ over $D_X$ has several connected components, and $R_{\phi}\mid_{D_X}$ is its component of unity. Our assertion follows.
\end{Prf} 
 
\bigskip

  The map $\gp_{R_{\phi}}:\Bun_{R_{\phi}}\to\Pic Y$ sends $(\cB_1,\cB_2, \gamma_{12}, N(\cB_1)\,\iso\, N(\cB_2))$ to $\cB_1^{\star}$. The square is cartesian
$$
\begin{array}{ccc}
\tilde Y & \toup{\tilde\phi} & \tilde X\\ 
\downarrow\lefteqn{\scriptstyle \alpha} && \downarrow\lefteqn{\scriptstyle \pi}\\
Z & \toup{\beta} & X
\end{array}
$$
The map $\rho_H\comp \gq_{R_{\phi}}: \Pic(Y,Z)\to\Bun_H$ sends $(\cB_1,\cB_2, \gamma_{12}, N(\cB_1)\,\iso\, N(\cB_2))$ as above to the collection $(V,\cC, \Sym^2 V\to\cC)$, where $V\subset \phi_* \cB_1\oplus \beta_*\cB_2$ is the lower modification defined by the cartesian square
\begin{equation}
\label{diag_def_of _lower_V}
\begin{array}{ccc}
\phi_* \cB_1\oplus \beta_*\cB_2 & \to (\phi_* \cB_1\oplus \beta_*\cB_2)\mid_{D_X}\to & \cB_1\mid_{D_Y}\oplus \cB_2\mid_{D_Z}\\
\uparrow && \uparrow\lefteqn{\scriptstyle \id+\gamma_{12}}\\
V & \to & \cB_1\mid_{D_Y},
\end{array}
\end{equation}
$\cC=N(\cB_1)(-D_X)$, and the quadratic form $\Sym^2 V\to\cC$ is the restriction of the difference of forms on $\cB_i$
$$
\Sym^2(\phi_* \cB_1\oplus \beta_*\cB_2)\to N(\cB_1)
$$ 

 Denote by $ZY$ the gluing of $Z$ and $Y$ along the isomorphism $D_Z\,\iso\, D_Y$. If $D_X$ is empty then $ZY$ is the disjoint union of $Z$ and $Y$. The projection $\varrho: ZY\to X$ is a 4-sheeted covering. 
 
  For a point $(\cB_1,\cB_2,\gamma_{12}, N(\cB_1)\,\iso\, N(\cB_2))\in \Pic(Y, Z)$ denote by $\cB_{12}$ the line bundle on $ZY$ obtained by gluing $\cB_1$ and $\cB_2$ along $\gamma_{12}$. Then $V$ from diagram (\ref{diag_def_of _lower_V}) is nothing but $V\,\iso\, \varrho_*\cB_{12}$.  Let $K_{ZY}$ denote the dualizing complex on $ZY$ then $K_{ZY}[-1]$ is the gluing of $\Omega_Y(D_Y)$ and $\Omega_Z(D_Z)$ via 
$$
\Omega_Y(D_Y)\mid_{D_Y}\,\iso\, \cO_{D_Y}\,\iso \,\cO_{D_Z}\, \iso\, \Omega_Z(D_Z)\mid_{D_Z}
$$
Set $(\cB_{12})^{\star}=\HOM(\cB_{12}, K_{ZY})[-1]$, this is
the gluing of $\cB_1^*\otimes\Omega_Y(D_Y)$ and $\cB_2^*\otimes\Omega_Z(D_Z)$ along the isomorphism 
$$
\gamma_{12}^*: \cB_2^*\otimes\Omega_Z(D_Z)\mid_{D_Z}\,\iso\, 
\cB_1^*\otimes\Omega_Y(D_Y)\mid_{D_Y}
$$
Then $(\varrho_*\cB_{12})^{\star}\,\iso\, \varrho_*(\cB_{12}^{\star})$ canonically. If $\phi$ is unramified then $ZY$ is smooth, and the definition of $(\cB_{12})^{\star}$ concides with that of 1.1. 
 
\medskip\noindent
6.3.1 In Sections 6.3.1-6.3.4 we assume $n=1$. 

Fix a two-sheeted covering $\phi:Y\to X$ as in 6.1.  For a given rank one local system $\cJ$ on $Y$ we want to calculate 
\begin{equation}
\label{complex_BP_square}
\RG(\Pic Y, A\cJ\otimes  (\gp_{R_{\phi}})_!\, \gq_{R_{\phi}}^*K_{\tilde E,\chi, \tilde H})
\end{equation}
Assume given an isomorphism $\chi\,\iso\, NJ$ of local systems on $X$. Then the complex that we integrate will descend under $e_{\phi}: \Pic Y\to\Bun_{U_{\phi}}$, we will actually integrate over $\Bun_{U_{\phi}}$. 

First, consider the situation when $E$ is an irreducible rank 2 local system on $X$ and $\tilde E=\pi^*E$. Let $\Aut_{E^*}$ denote the corresponding automorphic sheaf on $\Bun_2$ normamized as in \cite{FGV}. The following result is a calculation of
\begin{equation}
\label{complex_special_case_BP}
(\gp_{R_{\phi}})_!\, \gq_{R_{\phi}}^*F_{\tilde H}(\Aut_{E^*})[\dimrel(\gq_{R_{\phi}})],
\end{equation} 
under these assumptions. 

\begin{Th} 
\label{Th_special_case_BP}
Assume that $Y\toup{\phi} X$ is unramified and nonsplit, and the coverings $\tilde X\toup{\pi} X$ and $Y\toup{\phi} X$ are not isomorphic over $X$. 
Then the complex (\ref{complex_special_case_BP}) is isomorphic to
$$
\oplus_{d\ge 0} \, A(\cE_{0,\beta}\otimes \det E)_{\Omega}\otimes\mult^{d, d_1}_!(A(\cE_0\otimes\det E^*)\boxtimes (\phi^*E^*)^{(d)}) [d+g-1],
$$
where $d_1$ is a function of a connected component of $\Pic Y$ given by $2d_1+d=\deg\cB$ for $\cB\in\Pic Y$ (the sum is over $d\ge 0$ such that $d_1\in\ZZ$). Here 
$$
\mult^{d,d_1}: \Pic^{d_1} X\times Y^{(d)}\to\Pic^{2d_1+d} Y
$$ 
is the map introduced in Sect.~4.4.1.
\end{Th}

\begin{Rem} One may extend Theorem~\ref{Th_special_case_BP} to the case of $Y\toup{\phi} X$ split. To do so, first extend Theorem~\ref{Th_Rankin_Selberg_GO_2_GL_2} to the case of split $\tilde X$, then the argument of 6.3.2 will go through. We leave this to an interested reader. 
\end{Rem}

\medskip\noindent
6.3.2  Define $\Bun_{G,R_{\phi}}$ by the commutative diagram, where the square is cartesian
$$
\begin{array}{ccccc}
& \Bun_{G,R_{\phi}} & \to & \Bun_{G,\tilde H} & \toup{\gp} \Bun_G\\
& \downarrow && \downarrow\lefteqn{\scriptstyle \gq} \\
\Pic Y \getsup{\gp_{R_{\phi}}} & \Bun_{R_{\phi}} & \toup{\gq_{R_{\phi}}} & \Bun_{\tilde H}
\end{array}
$$ 
Then $\Bun_{G, R_{\phi}}$ classifies: $M\in\Bun_2$ for which we set $\cA=\det M$ and $\cC=\Omega\otimes\cA^{-1}$, $\cB_1\in\Pic Y$, $\cB_2\in\Pic Z$,  isomorphisms $N(\cB_1)\,\iso\, N(\cB_2)\,\iso\, \cC(D_X)$ and its refinement $\gamma_{12}: \cB_1\mid_{D_Y}\,\iso\, \cB_2\mid_{D_Z}$. 
 
\medskip 
 
\begin{Prf}\select{of Theorem~\ref{Th_special_case_BP}}

\smallskip
\noindent
Let $\tilde H_{\beta}$ (resp., $\tilde H_{\phi}$) denote the group scheme on $X$ obtained as a twisting of $\GO_2^0$ by the $\ZZ/2\ZZ$-torsor $\beta: Z\to X$ (resp., $\phi: Y\to X$) as in 3.1.  

 We have a commutative diagram
$$
\begin{array}{ccc}
\Bun_{G, R_{\phi}} & \toup{f}  \Bun_{G, \tilde H_{\beta}}\times_{\Bun_G} \Bun_{G,\tilde H_{\phi}}  \toup{\tilde\tau\times\tilde\tau} & \Bunt_{G_2}\times\Bunt_{G_2}\\
\downarrow && \downarrow\lefteqn{\scriptstyle \tilde\pi_2}\\
\Bun_{G,\tilde H} & \toup{\tilde\tau} & \Bunt_{G_4},
\end{array}
$$ 
where $\tilde\pi_2$ is the map defined in 3.5 for $X\sqcup X$. Here $f$ is the isomorphism sending the above point of $\Bun_{G, R_{\phi}}$ to $(\cB_2, M, N(\cB_2)\otimes\cA\,\iso\,\Omega)\in\Bun_{G, \tilde H_{\beta}}$, $(\cB_1, M, N(\cB_1)\otimes\cA\,\iso\,\Omega)\in\Bun_{G,\tilde H_{\phi}}$.  

 By Proposition~\ref{Pp_Weil_and_coverings}, 
$$
\tilde\pi_2^*\Aut[\dimrel]\,\iso\, \Aut\boxtimes\Aut,
$$ 
where $\dimrel=2\dim\Bun_{G_2}-\dim\Bun_{G_4}$.  We must calculate the direct image under the composition of projections
$$
\Bun_{G, \tilde H_{\beta}}\times_{\Bun_G} \Bun_{G,\tilde H_{\phi}} \to \Bun_{G,\tilde H_{\phi}}\to \Bun_{\tilde H_{\phi}}
$$
To calculate the direct image with respect to the first map we use
Proposition~\ref{Pp_lift_for_GL2} applied to the functor $F_G^{\tilde H_{\beta}}: \D(\Bun_{\tilde H_{\beta}})\to \D(\Bun_G)$. It yields an isomorphism
$$
F_G^{\tilde H_{\beta}}(\Qlb)[\dim\Bun_{\tilde H_{\beta}}]\,\iso\, \Aut_{\cE_{0,\beta}\oplus\Qlb}\otimes (A\cE_{0,\beta})_{\Omega}
$$
So, (\ref{complex_special_case_BP}) identifies with 
$$
(A\cE_{0,\beta})_{\Omega}\otimes \epsilon_!F_{\tilde H_{\phi}}^G(\Aut_{E^*}\otimes\Aut_{\cE_{0,\beta}\oplus\Qlb})[-\dim\Bun_G]
$$
Applying Theorem~\ref{Th_Rankin_Selberg_GO_2_GL_2} for $F_{\tilde H_{\phi}}$, one identifies (\ref{complex_special_case_BP}) with the direct sum
$$
\oplus_{d\ge 0} \, A(\cE_{0,\beta}\otimes \det E)_{\Omega}\otimes\mult^{d, d_1}_!(A(\cE_0\otimes\det E^*)\boxtimes (\pi^*E^*)^{(d)}) [d+g-1],
$$
where $d_1$ is a function of a connected component of $\Pic Y$ given by $2d_1+d=\deg\cB$ for $\cB\in\Pic Y$. The sum is over $d\ge 0$ such that $d_1\in\ZZ$. 
\end{Prf}

\medskip\noindent
6.3.3 {\scshape Geometric Waldspurger periods}  \  
In this subsection we assume $\tilde X$ split, $\phi: Y\to X$ nonsplit. From Theorem~\ref{Th_special_case_BP} combined with Proposition~\ref{Pp_Shimizu} one derives the following.

\begin{Cor} 
\label{Cor_Cor6}
Assume $\phi: Y\to X$ nonramified. For an irreducible rank 2 local system $E$ on $X$ we have
\begin{multline}
\label{complex_Cor6}
(\gp_{R_{\phi}})_!\gq_{R_{\phi}}^* K_{\pi^*E, \det E, \tilde H}[\dimrel(\gq_{R_{\phi}})] \,\iso\,\\
\oplus_{d\ge 0} \, (A\cE_{0,\phi})_{\Omega}\otimes \mult^{d,d_1}_!
(A(\det E^*)\boxtimes (\phi^* E^*)^{(d)})[d+g-1],
\end{multline}
where $d_1$ is a function of a connected component of $\Pic Y$ given by $2d_1+d=\deg\cB$ for $\cB\in\Pic Y$ (the sum is over $d\ge 0$ such that $d_1\in\ZZ$). 
\QED
\end{Cor}

\begin{Rem} By Remark~\ref{Rem_sum_two_components}, (\ref{complex_Cor6}) is naturally a direct sum of 2 complexes indexed by those $\theta\in\pi_1(\tilde H)$ whose image in $\ZZ$ is $2\deg\Omega-\deg\cB$, $\cB\in\Pic Y$. However, the RHS of (\ref{complex_Cor6}) seems not to be the refinement of this decomposition (cf. Remark~\ref{Rem_not_refinement}).
\end{Rem}

Remind the exact sequence 
$1\to \Gm\to\phi_*\Gm\to U_{\phi}\to 1$ on $X$ (cf. A.1). The corresponding extension of scalars map $e_{\phi}: \Pic Y\to \Bun_{U_{\phi}}$ sends $\cB$ to $\cB^{-1}\otimes\sigma_{\phi}^*\cB$. If $\phi: Y\to X$ is unramified then $U_{\phi}$ is also the kernel of the norm map $\phi_*\Gm\to\Gm$. 
 
 For $a\in\ZZ/2\ZZ$ we write $\Bun^a_{U_{\phi}}$ for the connected component of $\Bun_{U_{\phi}}$ corresponding to $a$, so $\Bun^0_{U_{\phi}}$ is the connected component of unity (cf. A.1). Let $\phi_1:\Pic Y\to\Bun_2$ be the map sending $\cB$ to $\phi_*B$. 
 
\begin{Def} 
\label{Def_WP}
Let $\cJ$ be a rank one local system on $Y$. Let $K\in\D(\Bun_2)$ be a complex with central character $N(\cJ)$. Then the complex $A\cJ^{-1}\otimes \phi_1^* K$ is equipped with natural descent data for the map $e_{\phi}:\Pic Y\to\Bun_{U_{\phi}}$. Assume that \footnote{though each perverse cohomology of the latter complex descends with respect to $e_{\phi}$, we ignore if the same is true for the complex itself, as the fibres of $e_{\phi}$ are not contractible.} 
\begin{itemize}
\item[($C_W$)] $\cK_K$ is a complex on $\Bun_{U_{\phi}}$ equipped with 
$$
e_{\phi}^*\cK_K[\dimrel(e_{\phi})]\,\iso\, A\cJ^{-1}\otimes \phi_1^* K[\dimrel(\phi_1)]
$$
\end{itemize}
For $a\in\ZZ/2\ZZ$ the \select{Waldspurger period} of $K$ is
$$
\WP^a(K,\cJ)=\RG_c(\Bun^a_{U_{\phi}}, \cK_K)
$$
\end{Def}
 
  Let $m_{\phi,d}: Y^{(d)}\to \Bun_{U_{\phi}}$ be the map sending $D$ to $\cO(D-\sigma_{\phi}^*D)$ with natural trivializations $N(\cO(D-\sigma_{\phi}^*D))\,\iso\,\cO$ and $\cO(D-\sigma_{\phi}^*D)\mid_{D_Y}\,\iso\,\cO_{D_Y}$. The map $m_{\phi,d}$ is proper. 
 
 Let $\mult_{\phi}: \Bun_{U_{\phi}}\times\Bun_{U_{\phi}} \to \Bun_{U_{\phi}}$ denote the multiplication map ($\Bun_{U_{\phi}}$ has a natural structure of a group stack). If $\phi$ is ramified then $\mult_{\phi}$ is proper. 
 
\begin{Th} 
\label{Th_WP}
Assume $\phi:Y\to X$ unramified. Let $E$ be an irreducible rank 2 local system on $X$, $\cJ$ be a rank one local system on $Y$ equipped with $\det E\,\iso\, N(\cJ)$. The condition ($C_W$) is satisfied for $\Aut_E$ giving rise to $\cK_E:=\cK_{\Aut_E}$. The complex $\cK_E$ is a direct sum of (possibly shifted) perverse sheaves. We have
$$
\mult_{\phi !}(\cK_E\boxtimes \cK_E)\,\iso\, \oplus_{d\ge 0} (A\cE_{0,\phi}\otimes A(N\cJ)^*)_{\Omega}\otimes
(m_{\phi,d})_!(\cJ\otimes\phi^* E^*)^{(d)}[d]
$$
In particular, for $a\in\ZZ/2\ZZ$ there are isomorphisms
\begin{multline*}
\mathop{\oplus}_{
\begin{array}{c}
\scriptstyle  a_1+a_2=a, \\
\scriptstyle  a_i\in\ZZ/2\ZZ
\end{array}
} \!\WP^{a_1}(\Aut_E,\cJ)\otimes\WP^{a_2}(\Aut_E,\cJ)\,
\iso\\ 
(A\cE_{0,\phi}\otimes A(N\cJ)^*)_{\Omega}\otimes
(\!\!\!\!\!\mathop{\oplus}_{
\begin{array}{c}
\scriptstyle d\ge 0,\\ 
\scriptstyle a=d\!\!\!\!\mod \!\! 2
\end{array}
}\!\!\!\!\! \RG(Y^{(d)}, (\cJ\otimes\phi^* E^*)^{(d)})[d])
\end{multline*}
If $\phi^*E$ is irreducible then the latter complex is a vector space (placed in cohomological degree zero). 
\end{Th}
\begin{Prf} Under our assumptions the sequence (\ref{seq_R_phi}) fits as the low row in the diagram
$$
\begin{array}{ccccccc}
1\to & \Gm\times\Gm & \to & \phi_*\Gm\times\phi_*\Gm & \to & U_{\phi}\times U_{\phi}& \to 1\\
 & \uparrow && \uparrow\lefteqn{\scriptstyle \id} && \uparrow \\ 
1\to &  \Gm & \to & \phi_*\Gm\times\phi_*\Gm & \to  & R_{\phi} & \to 1,
\end{array}
$$ 
where the left vertical arrow sends $z$ to $(z, z^{-1})$. Let $\kappa_{\phi}: \Bun_{R_{\phi}}\to \Bun_{U_{\phi}}\times\Bun_{U_{\phi}}$ denote the corresponding extension of scalars map. The following diagram is cartesian
$$
\begin{array}{ccc}
\Bun_{R_{\phi}} & \toup{\gp_{R_{\phi}}} & \Pic Y \\
\downarrow\lefteqn{\scriptstyle \kappa_{\phi}} && \downarrow\lefteqn{\scriptstyle e_{\phi}\comp\epsilon}\\
\Bun_{U_{\phi}}\times\Bun_{U_{\phi}} & \toup{\mult_{\phi}} & \Bun_{U_{\phi}},
\end{array}
$$ 
 
 By Remark~\ref{Rem_ULA_and_C_W}, the condition ($C_W$) is satisfied, and we get
$$
(A\cJ)^{-1}_{\Omega_Y}\otimes \gp_{R_{\phi}}^*A\cJ\otimes \gq_{R_{\phi}}^*K_{\pi^*E, \det E, \tilde H}[\dimrel(\gq_{R_{\phi}})]\,\iso\, \kappa_{\phi}^*(\cK_E\boxtimes \cK_E)[\dimrel(\kappa_{\phi})]
$$
By Corolary~\ref{Cor_Cor6}, we get
\begin{multline*}
(e_{\phi}\comp\epsilon)^*\mult_{\phi !}(\cK_E\boxtimes \cK_E)[\dimrel(\kappa_{\phi})]\,\iso\\
(A\cE_{0,\phi}\otimes A(N\cJ)^*)_{\Omega}\otimes 
(\oplus_{d\ge 0} \; A\cJ\otimes\mult^{d,d_1}_!
(A(\det E^*)\boxtimes (\phi^* E^*)^{(d)})[d+g-1])\,\iso\\
(A\cE_{0,\phi}\otimes A(N\cJ)^*)_{\Omega}\otimes 
(\oplus_{d\ge 0} \; \mult^{d,d_1}_! (\Qlb\boxtimes (\cJ\otimes\phi^* E^*)^{(d)})[d+g-1]),
\end{multline*}
where $d_1$ is a function of a connected component of $\Pic Y$ given by $2d_1+d=\deg\cB$ for $\cB\in\Pic Y$ (the sum is over $d\ge 0$ such that $d_1\in\ZZ$). 

  The square is cartesian 
$$
\begin{array}{ccc}
\Pic^{d_1} X\times Y^{(d)} & \toup{\mult^{d,d_1}} & \Pic^{2d_1+d} Y\\
\downarrow && \downarrow\lefteqn{\scriptstyle e_{\phi}\comp\epsilon}\\
Y^{(d)} & \toup{m_{\phi,d}} & \Bun_{U_{\phi}},
\end{array}
$$
where the left vertical arrow is the projection. 
 Since $\dimrel(\kappa_{\phi})=g-1$, we get an isomorphism
$$
(e_{\phi}\comp\epsilon)^*\mult_{\phi !}(\cK_E\boxtimes \cK_E)\,\iso\, 
(A\cE_{0,\phi}\otimes A(N\cJ)^*)_{\Omega}\otimes
(e_{\phi}\comp\epsilon)^*(\oplus_{d\ge 0} \; (m_{\phi,d})_!(\cJ\otimes\phi^* E^*)^{(d)}[d])
$$
compatible with the descent data for $e_{\phi}\comp\epsilon$. So, 
$$
\mult_{\phi !}(\cK_E\boxtimes \cK_E)\,\iso\,
(A\cE_{0,\phi}\otimes A(N\cJ)^*)_{\Omega}\otimes
(\oplus_{d\ge 0} \; (m_{\phi,d})_!(\cJ\otimes\phi^* E^*)^{(d)}[d]),
$$
the sum over all $d\ge 0$. For $d$ even (resp., odd) $m_{\phi,d}$ maps $Y^{(d)}$ to $\Bun_{U_{\phi}}^0$ (resp., to $\Bun_{U_{\phi}}^1$). 

 If $\phi^*E$ is irreducible then $\RG(Y^{(d)}, (\cJ\otimes \phi^* E^*)^{(d)})[d]\,\iso\wedge^d V$ with $V=\H^1(Y, \cJ\otimes \phi^* E^*)$. The last statement follows. 
\end{Prf}

\medskip\noindent
6.3.4 In this subsection $\phi:Y\to X$ is allowed to be ramified. Let us calculate the geometric Waldspurger periods of Eisenstein series on $\Bun_2$.
Let $E_1, E_2$ be rank one local systems on $X$, $\cJ$ be a rank one local system on $Y$ equipped with $N(\cJ)\,\iso\, E_1\otimes E_2$. Remind the complex $\Aut_{E_1\oplus E_2}$ on $\Bun_2$ (cf. 4.3). Let 
$$
\tilde m_{\phi,d}: Y^{(d)}\to\Bun_{U_{\phi}}
$$ 
be the map sending $D$ to $\cO(\sigma_{\phi}^*D-D)$ with canonical trivialization $N(\cO(\sigma_{\phi}^*D-D))\,\iso\,\cO_X$. 

\begin{Pp} 
\label{Pp_Eis_WP}
The condition ($C_W$) is satisfied for $\Aut_{E_1\oplus E_2}$. For the corresponding complex $\cK_{E_1\oplus E_2}:=\cK_{\Aut_{E_1\oplus E_2}}$ we have
\begin{equation}
\label{complex_K_Eis}
\cK_{E_1\oplus E_2}\,\iso\, \oplus_{d\ge 0} (AE_2)_{\cE_{\phi}\otimes\Omega^{-1}}\otimes (\tilde m_{\phi,d})_!(\cJ^{-1}\otimes\phi^* E_2)^{(d)}[d]
\end{equation}
In particular, 
$$
\mult_{\phi\, !}(\cK_{E_1\oplus E_2}\boxtimes \cK_{E_1\oplus E_2})\,\iso\,  \oplus_{d\ge 0} A(E_1\otimes E_2)_{\cE_{\phi}\otimes\Omega^{-1}}\otimes (\tilde m_{\phi,d})_!(\cJ^{-1}\otimes\phi^*(E_1\oplus E_2))^{(d)}[d]
$$
\end{Pp}
\begin{Prf}
We have a cartesian square
$$
\begin{array}{ccc}
\mathop{\sqcup}\limits_{d\ge 0} \Pic X\times Y^{(d)} & \toup{mult} & \Pic Y\\
\downarrow\lefteqn{\scriptstyle \phi_{1,P}} && \downarrow\lefteqn{\scriptstyle \phi_1}\\
\Bunb_P & \toup{\bar\gp_P} & \Bun_2,
\end{array}
$$
where $\phi_{1,P}$ sends $(L,D)$ to $(L\subset M)$, $M=L\otimes\phi_*\cO(D)$. The map $\mult^{d,d_1}:\Pic^{d_1} X\times Y^{(d)}\to \Pic^{2d_1+d}$ sends $(L,D)$ to $(\phi^*L)(D)$. So, we have
$$
\phi_1^*\Aut_{E_1\oplus E_2}\otimes A\cJ^{-1}[\dimrel(\phi_1)]\,\iso\, \mathop{\oplus}\limits_{d\ge 0} (AE_2)_{\cE_{\phi}\otimes\Omega^{-1}}\otimes \mult^{d,d_1}_!(\Qlb\boxtimes (\cJ^{-1}\otimes\phi^* E_2)^{(d)})[d+g-1],
$$
where $d_1$ is a function of a connected component of $\Pic Y$ given by $2d_1+d=\deg\cB$, $\cB\in\Pic Y$ (and the sum is over $d\ge 0$ such that $d_1\in\ZZ$). 
We have used that $\dimrel(\phi_1)=2(1-g)+\frac{1}{2}\deg D_X$. 

 The square is cartesian 
$$
\begin{array}{ccc}
\Pic^{d_1} X\times Y^{(d)} & \toup{\mult^{d,d_1}} & \Pic^{2d_1+d} Y\\
\downarrow && \downarrow\lefteqn{\scriptstyle e_{\phi}}\\
Y^{(d)} & \toup{\tilde m_{\phi,d}} & \Bun_{U_{\phi}},
\end{array}
$$
where the left vertical arrow is the projection. This yields an isomorphism
$$
\phi_1^*\Aut_{E_1\oplus E_2}\otimes A\cJ^{-1}[\dimrel(\phi_1)]\,\iso\, \mathop{\oplus}\limits_{d\ge 0} (AE_2)_{\cE_{\phi}\otimes\Omega^{-1}}\otimes e_{\phi}^*
(\tilde m_{\phi,d})_!(\cJ^{-1}\otimes\phi^* E_2)^{(d)})[d+g-1]
$$
Since $\dimrel(e_{\phi})=g-1$, the first assertion follows. 

  To get the second one, remind that for rank one local systems $V_i$ on $Y$ we have 
$$
(V_1\oplus V_2)^{(d)}\,\iso\, \oplus_{0=k}^d (\sym_{k,d-k})_!(V_1^{(k)}\boxtimes V_2^{(d-k)}),
$$
where $\sym_{k,d-k}: Y^{(k)}\times Y^{(d-k)}\to Y^{(d)}$ is the sum of divisors. 
\end{Prf}

\medskip
\begin{Rem} i) Remind that the Eisenstein series $\Aut_{E_1\oplus E_2}$ above is called regular if $E_1$ and $E_2$ are not isomorphic (cf. \cite{BG}, Sect.~2.1.7). 
Under this assumptions $(\tilde m_{\phi,d})_!(\cJ^{-1}\otimes\phi^* E_2)^{(d)}=0$ for $d> 2g_Y-2$ (so, the sum in (\ref{complex_K_Eis}) is actually by $0\le d\le 2g_Y-2$). 

 Indeed, $\tilde m_{\phi,d}$ decomposes as $Y^{(d)}\toup{s} \Pic Y\toup{e_{\phi}}\Bun_{U_{\phi}}$, where $s$ sends $D$ to $\cO(D)$. If $d> 2g_Y-2$ then $s$ is a vector bundle of rank $d+1-g_Y$ over $\Pic Y$ with zero section removed. For a rank one local system $\tilde E$ on $Y$ we have $s_! \tilde E^{(d)}\,\iso\, A\tilde E\otimes s_!\Qlb$. Further, $e_{\phi}: \Pic Y\to\Bun_{U_{\phi}}$ is a homomorphism of group stacks, each fibre of $e_{\phi}$ identifies with $\Pic X$. So, if $N(\tilde E)$ is nontrivial then $(e_{\phi})_!A\tilde E=0$.

\smallskip\noindent
ii) If $E$ is a rank 2 local system on $X$, $\cJ$ is a rank one local system on $Y$ equipped with $\det E\,\iso\, N(\cJ)$ then we have $\cJ^{-1}\otimes\phi^*E\,\iso\, \sigma_{\phi}^*\cJ\otimes \phi^* E^*$ canonically. Write also $\sigma_{\phi}: Y^{(d)}\to Y^{(d)}$ for the map sending $D$ to $\sigma_{\phi}^*D$. Then the composition $Y^{(d)}\toup{\sigma_{\phi}} Y^{(d)}\toup{m_{\phi,d}} \Bun_{U_{\phi}}$ equals $\tilde m_{\phi,d}$. So, 
$$
(\tilde m_{\phi,d})_!(\cJ^{-1}\otimes\phi^* E)^{(d)}\,\iso\, 
(m_{\phi,d})_!(\cJ\otimes \phi^* E^*)^{(d)}
$$
canonically. Thus, Theorem~\ref{Th_WP} and Proposition~\ref{Pp_Eis_WP} are consistent. 
\end{Rem}

\medskip

 Theorem~\ref{Th_WP} and Proposition~\ref{Pp_Eis_WP} suggest the following conjecture (it is a theorem if one of the following holds
\begin{itemize} 
\item  $E$ is irreducible and $\phi$ is nonramified, 
\item $E$ is a direct sum of two rank one local systems. 
\end{itemize}
 
\begin{Con}[Waldspurger periods]
\label{Con_WP} Let $E$ be a rank 2 local system on $X$. Let $K\in\D(\Bun_2)$ be an automorphic sheaf with eigenvalue $E$. Let $\phi: Y\to X$ be a (possibly ramified) degree 2 covering, $\cJ$ be a rank one local system on $Y$. Assume the condition ($C_W$) satisfied for $\cJ$ and $K$ giving rise to $\cK_K\in\D(\Bun_{U_{\phi}})$. Then for a suitable normalization of $K$ there exists an isomorphism
$$
\mult_{\phi !}(\cK_K\boxtimes \cK_K)\,\iso\, \oplus_{d\ge 0} 
(m_{\phi,d})_!(\cJ\otimes\phi^* E^*)^{(d)}[d]
$$
\end{Con}

\bigskip\noindent
6.3.5 {\scshape Geometric Bessel periods}

\medskip\noindent
In this subsection we assume $\phi:Y\to X$ nonsplit. Assume $n=2$, so $G=\GSp_4$. 

\begin{Def} 
\label{def_BP}
Let $K\in\D(\Bun_G)$ be a complex with central character $\chi^{-1}$. Let $\cJ$ be a rank one local system on $Y$ equipped with $N(\cJ)\,\iso\, \chi$. For the inclusion $\Pic Y\hook{} {^{rss}\cS_P}\subset \cS_P$ the $*$-restriction
$A\cJ\otimes
\Four_{\psi}(\nu_P^*K)\mid_{\Pic Y}$
is equipped with natural descent data for $e_{\phi}:\Pic Y\to \Bun_{U_{\phi}}$. Assume that
\begin{itemize}
\item[($C_B$)] $\cK_K$ is a complex on $\Bun_{U_{\phi}}$ equipped with 
$$
e_{\phi}^*\cK_K[\dimrel(e_{\phi})]\,\iso\, 
A\cJ\otimes
\Four_{\psi}(\nu_P^*K)[\dimrel(\nu_P)]\mid_{\Pic Y}
$$
\end{itemize}
For $a\in\ZZ/2\ZZ$ the \select{Bessel period} of $K$ is
$$
\BP^a(K,\cJ)=\RG_c(\Bun_{U_{\phi}}^a, \cK_K)
$$
\end{Def} 

 Assume $\tilde X$ connected. Let $\tilde E$ be an irreducible rank 2 local system on $\tilde X$, $\chi$ be a rank one local system on $X$ equipped with $\pi^*\chi\,\iso\, \det\tilde E$. Remind the complex $K_{\tilde E,\chi,\tilde H}$ defined in 5.1. 
Remind the map $\tilde m_{\phi,d}: Y^{(d)}\to \Bun_{U_{\phi}}$ introduced in 6.3.4. 
  
\begin{Th} 
\label{Th_BP}
Assume $\phi: Y\to X$ nonramified. Let $\cJ$ be a rank one local system on $Y$ equipped with $N(\cJ)\,\iso\,\chi$. The condition ($C_B$) is satisfied for $\cJ$ and $K:=F_G(\rho_{H\, !}K_{\tilde E,\chi,\tilde H})$ giving rise to $\cK_K\in\D(\Bun_{U_{\phi}})$. We have
$$
(\mult_{\phi})_!(\cK_K\boxtimes\cK_K)\,\iso\, 
\oplus_{d\ge 0} (\tilde m_{\phi,d})_!(\cJ\otimes\phi^*(\pi_* \tilde E^*))^{(d)}[d]
$$
In particular, for $a\in\ZZ/2\ZZ$ there are isomorphisms
$$
\mathop{\oplus}_{
\begin{array}{c}
\scriptstyle  a_1+a_2=a, \\
\scriptstyle  a_i\in\ZZ/2\ZZ
\end{array}
} \!\BP^{a_1}(K,\cJ)\otimes\BP^{a_2}(K,\cJ)\,
\iso\\ 
\!\!\!\!\!\mathop{\oplus}_{
\begin{array}{c}
\scriptstyle d\ge 0,\\ 
\scriptstyle a=d\!\!\!\!\mod \!\! 2
\end{array}
}\!\!\!\!\! \RG(Y^{(d)}, (\cJ\otimes\phi^*(\pi_*\tilde E^*))^{(d)})[d]
$$
\end{Th}
\begin{Prf}  According to Proposition~\ref{Pp_R} and Corollary~\ref{Cor_relation_F_P_and_F_G}, 
we must find a complex $\cK_K\in\D(\Bun_{U_{\phi}})$ together with an isomorphism
\begin{equation}
\label{iso_to_find_Th5}
e_{\phi}^*\cK_K[\dimrel(e_{\phi})]\,\iso\,
A\cJ\otimes (\gp_{R_{\phi}})_!\gq_{R_{\phi}}^*K_{\tilde E,\chi,\tilde H}[\dimrel(\gq_{R_{\phi}})]
\end{equation}

 Let us show that $R_{\phi}$ fits into a cartesian square
\begin{equation}
\label{square_for_Th5}
\begin{array}{ccc}
R_{\phi} & \to & \pi_* U_{\tilde\phi}\\
\downarrow\lefteqn{\scriptstyle N_Y} &&  \downarrow\lefteqn{\scriptstyle N_Y}\\
\phi_*\Gm & \to & U_{\phi}
\end{array}
\end{equation}
Indeed, (\ref{seq_R_phi}) fits into a commutative diagram
$$
\begin{array}{ccccccc}
1\to &\Gm & \to & \phi_*\Gm & \to & U_{\phi} & \to 1\\
& \uparrow\lefteqn{\scriptstyle N} && \uparrow\lefteqn{\scriptstyle N_Y} && \uparrow\lefteqn{\scriptstyle N_Y}\\
1\to  & \pi_*\Gm & \to &\pi_*\tilde\phi_*\Gm & \to &\pi_* U_{\tilde\phi}& \to 1\\
& \uparrow && \uparrow\lefteqn{\scriptstyle \id} && \uparrow\\
1\to &U_{\pi} & \to &\pi_*\tilde\phi_* \Gm & \to &R_{\phi}& \to 1\\
 & \downarrow && \downarrow\lefteqn{\scriptstyle N_Y} && \downarrow\lefteqn{\scriptstyle N_Y}\\
& 1 & \to & \phi_*\Gm & \toup{\id} & \phi_*\Gm &\to 1
\end{array}
$$
The latter diagram together with the exact sequence $1\to \Gm\to \phi_*\Gm\to U_{\phi}\to 1$ yield (\ref{square_for_Th5}). 

 Let $\kappa_{\phi}:\Bun_{R_{\phi}}\to \Bun_{U_{\tilde\phi}}$ be the extension of scalars map given by the upper row in (\ref{square_for_Th5}). 
The composition $\Pic \tilde Y\to\Bun_{R_{\phi}}\toup{\kappa_{\phi}}\Bun_{U_{\tilde\phi}}$ is the map $e_{\tilde\phi}$ sending $\cB$ to $\sigma_{\phi}^*\cB\otimes \cB^{-1}$. We get a cartesian square
$$
\begin{array}{ccc}
\Bun_{R_{\phi}} & \toup{\kappa_{\phi}} & \Bun_{U_{\tilde\phi}}\\
\downarrow\lefteqn{\scriptstyle \gp_{R_{\phi}}} && \downarrow\lefteqn{\scriptstyle \tilde N_Y}\\
\Pic Y & \toup{e_{\phi}} & \Bun_{U_{\phi}},
\end{array}
$$
where $\tilde N_Y(\cB):=N_Y(\cB)^{-1}$ for $\cB\in\Pic\tilde Y$. 

 Consider the commutative diagram
$$
\begin{array}{cccc}
\Pic \tilde Y & \to & \Bun_{R_{\phi}} & \toup{\kappa_{\phi}} \Bun_{U_{\tilde\phi}}\\
\downarrow\lefteqn{\scriptstyle \tilde\phi_1} && \downarrow\lefteqn{\scriptstyle \gq_{R_{\phi}}}\\
\Bun_{2,\tilde X} & \to & \Bun_{\tilde H},
\end{array}
$$
where $\tilde\phi_1$ sends $\cB$ to $\tilde\phi_*\cB$. 

 By Remark~\ref{Rem_ULA_and_C_W}, the condition ($C_W$) is satisfied for $\Aut_{\tilde E}$, the covering $\tilde\phi: \tilde Y\to\tilde X$ and the local system  $\tilde\pi^*\cJ$. So, there is a complex $\cK_{\tilde E}\in\D(\Bun_{U_{\tilde\phi}})$ equipped with
$$
e_{\tilde\phi}^*\cK_{\tilde E}[\dimrel(e_{\tilde\phi})]\,\iso\, 
A(\pi^*\cJ)^{-1}\otimes\tilde\phi_1^*\Aut_{\tilde E}[\dimrel(\tilde\phi_1)]
$$
Set $\tilde\cK=\cK_{\tilde E}\otimes (A\cJ)_{\Omega_Y}$, it is equipped with an isomorphism
$$
\gp_{R_{\phi}}^*A\cJ\otimes \gq_{R_{\phi}}^* K_{\tilde E,\chi, \tilde H}[\dimrel(\gq_{R_{\phi}})]\,\iso\, \kappa_{\phi}^*\tilde\cK[\dimrel(\kappa_{\phi})]
$$
over $\Bun_{R_{\phi}}$. Set $\cK_K=\tilde N_{Y !}(\tilde \cK)$, it is equipped with an isomorphism (\ref{iso_to_find_Th5}). 
We have a commutative diagram
$$
\begin{array}{ccc}
\Bun_{U_{\tilde\phi}}\times\Bun_{U_{\tilde\phi}} & \toup{\mult_{\tilde\phi}} & \Bun_{U_{\tilde\phi}}\\
\downarrow\lefteqn{\scriptstyle \tilde N_Y\times \tilde N_Y} && \downarrow\lefteqn{\scriptstyle \tilde N_Y}\\
\Bun_{U_{\phi}}\times\Bun_{U_{\phi}} & \toup{\mult_{\phi}} & \Bun_{U_{\phi}}
\end{array}
$$
By Theorem~\ref{Th_WP}, we have
$$
(\mult_{\tilde\phi})_!(\tilde\cK\boxtimes\tilde\cK)\,\iso\, 
\oplus_{d\ge 0} (A\cE_{0,\tilde\phi})_{\Omega_{\tilde X}}\otimes
(m_{\tilde\phi,d})_!(\tilde\pi^*\cJ\otimes\tilde\phi^* \tilde E^*)^{(d)}[d],
$$
where $m_{\tilde\phi,d}: \tilde Y^{(d)}\to\Bun_{U_{\tilde\phi}}$ sends $D$ to $\cO(D-\sigma_{\phi}^*D)$ with natural trivialization 
$$
N_{\tilde X}(\cO(D-\sigma_{\phi}^*D))\,\iso\, \cO_{\tilde X}
$$
Since $\cE_{0,\tilde\phi}\,\iso\, \pi^*\cE_{0,\phi}$, we get $(A\cE_{0,\tilde\phi})_{\Omega_{\tilde X}}\,\iso\, A(N\cE_{0,\tilde\phi})_{\Omega}\,\iso\, \Qlb$. So, 
$$
(\mult_{\phi})_!(\cK_K\boxtimes\cK_K)\,\iso\, (\tilde N_Y)_!(\mult_{\tilde\phi})_!(\tilde\cK\boxtimes\tilde\cK)\,\iso\,
\oplus_{d\ge 0} (\tilde N_Y)_!(m_{\tilde\phi,d})_!(\tilde\pi^*\cJ\otimes\tilde\phi^* \tilde E^*)^{(d)}[d]
$$
The diagram commutes
$$
\begin{array}{ccc}
\tilde Y^{(d)} & \toup{m_{\tilde\phi,d}} & \Bun_{U_{\tilde\phi}}\\
\downarrow\lefteqn{\scriptstyle \tilde\pi} && \downarrow\lefteqn{\scriptstyle \tilde N_Y}\\
Y^{(d)} & \toup{\tilde m_{\phi,d}} & \Bun_{U_{\phi}}
\end{array}
$$
Our assertion follows. 
\end{Prf}

\medskip

 Theorem~\ref{Th_BP} combined with Conjecture~\ref{Con_cusp} suggest the following.
 
\begin{Con}[Bessel periods]
\label{Con_BP}
For $G=\GSp_4$ let $E_{\check{G}}$ be a $\check{G}$-local system on $X$ viewed as a pair $(E,\chi^{-1})$, where $E$ (resp., $\chi$) is a rank 4 (resp., rank 1) local system on $X$ equipped with a symplectic form $\wedge^2 E\to \chi^{-1}$. Let $K$ be an automorphic sheaf on $\Bun_G$ with eigenvalue $E_{\check{G}}$ (in particular, the central character of $K$ is $\chi^{-1}$). 

 Let $\phi: Y\to X$ be a (possibly ramified) degree 2 covering. Let $\cJ$ be a rank one local system on $Y$ equipped with $N(\cJ)\,\iso\, \chi$. 
Assume the condition ($C_B$) satisfied for $\cJ$ and $K$.  If the corresponding complex $\cK_K\in\D(\Bun_{U_{\phi}})$ is nonzero then 
$$
(\mult_{\phi})_!(\cK_K\boxtimes\cK_K)\,\iso\, 
\oplus_{d\ge 0} (\tilde m_{\phi,d})_!(\cJ\otimes\phi^* E)^{(d)}[d]
$$
for a suitable normalization of $K$. 
In particular, for $a\in\ZZ/2\ZZ$ there are isomorphisms
\begin{equation}
\label{complex_Con_square_BP}
\mathop{\oplus}_{
\begin{array}{c}
\scriptstyle  a_1+a_2=a, \\
\scriptstyle  a_i\in\ZZ/2\ZZ
\end{array}
} \!\BP^{a_1}(K,\cJ)\otimes\BP^{a_2}(K,\cJ)\,
\iso\\ 
\!\!\!\!\!\mathop{\oplus}_{
\begin{array}{c}
\scriptstyle d\ge 0,\\ 
\scriptstyle a=d\!\!\!\!\mod \!\! 2
\end{array}
}\!\!\!\!\! \RG(Y^{(d)}, (\cJ\otimes\phi^* E)^{(d)})[d]
\end{equation}
\end{Con}

\begin{Rem}
\label{Rem_sym_form_on_V}
Under the assumptions of Conjecture~\ref{Con_BP} we have $\H^2(Y, \cJ\otimes\phi^*E)^*\,\iso\, \H^0(Y, \cJ\otimes\phi^*E)$. Consider the case $\H^i(Y, \cJ\otimes\phi^*E)=0$ for $i=0,2$. Then $\BP(K,\cJ)\otimes\BP(K,\cJ)$
identifies with the vector space (placed in degree zero)
\begin{equation}
\label{complex_Clifford}
\mathop{\oplus}_{i\ge 0} \wedge^i V,
\end{equation}
where $V=\H^1(Y, \cJ\otimes\phi^*E)$. The symplectic form on $E$ induces a map 
\begin{equation}
\label{map_sympl_form_H2}
\H^2(Y, \cJ\otimes\sigma_{\phi}^*\cJ\otimes\phi^*(E\otimes E))\to
\H^2(Y, \Qlb)\,\iso\,\Qlb
\end{equation}
Since the cup-product 
$$
V\otimes V\,\iso\, \H^1(Y, \cJ\otimes\phi^*E)\otimes \H^1(Y, \sigma_{\phi}^*\cJ\otimes\phi^*E)
\to \H^2(Y, \cJ\otimes\sigma_{\phi}^*\cJ\otimes\phi^*(E\otimes E))
$$
is anti-symmetric, composing it with (\ref{map_sympl_form_H2}) one gets a nondegenerate symmetric form $\Sym^2 V\to\Qlb$ on $V$. We have $\dim V=8(g_Y-1)$, where $g_Y$ is the genus of $Y$. 

 Let $\Spin(V)$ denote the simply-connected covering of $\SO(V)$. Let $\Gamma_{\alpha}$ and $\Gamma_{\beta}$ be the half-spin representations of $\Spin(V)$, here $\alpha$ and $\beta$ are the corresponding fundamental weights of $\Spin(V)$ (cf. \cite{FH}, 19.2, p. 287). Then
$$
\Gamma_{\alpha}\otimes\Gamma_{\alpha}\oplus
\Gamma_{\beta}\otimes\Gamma_{\beta}\,\iso\, \wedge^0 V\oplus\wedge^2 V\oplus\wedge^4 V\ldots
$$
and
$$
\Gamma_{\alpha}\otimes\Gamma_{\beta} \oplus
\Gamma_{\beta}\otimes\Gamma_{\alpha}\,\iso\,
\wedge^1 V\oplus\wedge^3 V\oplus \wedge^5 V\oplus\ldots
$$
\end{Rem}

\smallskip

\begin{Con}[Bessel periods refined] 
\label{Con_BP_refined}
Under the assumptions of Conjecture~\ref{Con_BP} consider the case $\H^i(Y, \cJ\otimes\phi^*E)=0$ for $i=0,2$. Set $V=\H^1(Y, \cJ\otimes\phi^*E)$. Then there is a numbering $\alpha_a$ ($a\in\ZZ/2\ZZ$) of the half-spin fundamental weights of $\Spin(V)$ and isomorphisms for $a\in\ZZ/2\ZZ$
$$
\BP^a(K,\cJ)\,\iso\, \Gamma_{\alpha_a},
$$
where $\Gamma_{\alpha_a}$ is the irreducible (half-spin) representation of $\Spin(V)$ with highest weight $\alpha_a$.
\end{Con}

\bigskip\medskip

\centerline{\scshape 7. The case $H=\GO_6$}

\bigskip\noindent
7.1 In this section we assume $m=3$ and $\tilde X$ split, so $\HH=\tilde H=\GO^0_6$. We have an exact sequence $1\to \mu_2\to \GL_4\to \GO_6^0\to 1$ of group schemes over $\Spec k$. By abuse of notaion, we write $\rho: \Bun_4\to\Bun_{\tilde H}$ for the corresponding extension of scalars map. It sends $W\in\Bun_4$ to 
$$
(V=\wedge^2 W, \cC=\det W, \Sym^2 V\toup{h}\cC, \gamma)
$$ 
Here $h$ is the symmetric form induced by the exteiour product $\wedge^2 W\otimes\wedge^2 W\to \det W$, and $\gamma: \det V\,\iso\, \cC^3$ is a compatible trivialization. 
  The connected components of $\Bun_{\tilde H}$ are indexed by $\pi_1(\tilde H)$. We have an exact sequence $0\to \pi_1(\GL_4)\to \pi_1(\tilde H)\to \ZZ/2\ZZ\to 0$, and the image of $\rho$ is $_0\Bun_{\tilde H}:=\sqcup_{a\in\pi_1(\GL_4)} \Bun_{\tilde H}^a$. 
  
  Remind that 
$$
\check{\HH}\,\iso\,\GSpin_6\,\iso\,\{(c,b)\in \Gm\times\GL_4\mid \det b=c^2\}
$$ 
Consider a $\check{\HH}$-local system on $X$ given by a collection: local systems $E$ and $\chi$ on $X$ of ranks 4 and 1 respectively, an isomorphism $\det E\,\iso\, \chi^2$.   
Assume $E$ irreducible on $X$. Let $\Aut_E$ denote the corresponding automorphic sheaf on $\Bun_4$ (cf. Definition~\ref{Def_Aut_E}). Then $\Aut_E$ is equipped with natural descent data with respect to $\rho:\Bun_4\to\Bun_{\tilde H}$, so gives rise to a perverse sheaf $K_{E,\chi,\tilde H}$ on $_0\Bun_{\tilde H}$. 

\begin{Lm} The sheaf $K_{E,\chi,\tilde H}$ extends naturally to a perverse sheaf (still denoted by the same symbol) over $\Bun_{\tilde H}$ with central character $\chi$. 
\end{Lm}
\begin{Prf} Let $\act:\Pic X\times\Bun_{\tilde H}\to\Bun_{\tilde H}$ be the action map sending $\cL\in\Pic X$ and $(V,\cC, \Sym^2 \to\cC)\in\Bun_{\tilde H}$ to $(V\otimes\cL, \, \cC\otimes\cL)$. Let $_1\Bun_{\tilde H}\subset \Bun_{\tilde H}$ denote the complement to $_0\Bun_{\tilde H}$. Then $\act$ sends  $\Pic^k X\times {_a\Bun_{\tilde H}}$ to $_b\Bun_{\tilde H}$, where $b=a+k\mod 2$. 
The perverse sheaf 
$$
A\chi\boxtimes K_{E,\chi,\tilde H}[g-1]
$$ 
on $\Pic^1 X\times {_0\Bun_{\tilde H}}$ is equipped with natural descent data for $\act: \Pic^1 X\times {_0\Bun_{\tilde H}}\to {_1\Bun_{\tilde H}}$. This yields a perverse sheaf $K_{E,\chi,\tilde H}$ on the whole of $\Bun_{\tilde H}$ equipped with
$\act^* K_{E,\chi,\tilde H}\,\iso\, A\chi\boxtimes K_{E,\chi,\tilde H}$. Here $A\chi$ is the automorphic local system corresponding to $\chi$. 
\end{Prf}

\medskip

 Assume $n=2$, so $G=\GSp_4$. View a $\check{G}$-local system $E_{\check{G}}$ on $X$ as a pair $(E,\chi)$, where $E$ (resp., $\chi$) is a rank 4 (resp., rank 1) local system on $X$ with symplectic form $\wedge^2 E\to \chi$. The symplectic form induces the isomorphism $\det E\,\iso\,\chi^2$, and $(E,\chi)$ identifies with the $\check{\HH}$-local system $E_{\check{\HH}}$ obtained from $E_{\check{G}}$ via the extension of scalars $\check{G}\hook{}\check{\HH}$. 
    
\begin{Con}
\label{Con_autom_GSp4_from_GO6}
i) Let $E_{\check{G}}$ be a $\check{G}$-local system on $X$ and $E_{\check{\HH}}$ be the induced $\check{\HH}$-local system on $X$ given by $(E,\chi)$. 
There exists $K\in \D(\Bun_G)$ which is a $E_{\check{G}}$-Hecke eigensheaf satisfying 
$F_{\tilde H}(K)\,\iso\, K_{E^*,\chi^*,\tilde H}$. 

\smallskip\noindent
ii) Assume in addition that $E$ is irreducible (as a local system of rank 4). Then $K=F_G(K_{E^*,\chi^*,\tilde H})$ satisfies the properties of i). 
\end{Con}

\medskip\noindent
7.2 Remind the stack $\RCov^r$ introduced in 6.1. Denote by $\Bun_{k, r}$ the following stack. For a scheme $S$, an $S$-point of $\Bun_{k,r}$ is a collection consisting of a map $S\to\RCov^r$ giving rise to a two-sheeted covering $_SY\to S\times X$, and a rank $k$ vector bundle on $_SY$. Let us precise that a map $S\to\RCov^r$ is given by a collection $(\cE,\kappa, D)$, where $D\hook{} S\times X$ is the preimage of the incidence divisor on $^{rss} X^{(r)}\times X$ under $S\times X\to {^{rss} X^{(r)}}\times X$, and $\cE$ is a line bundle on $S\times X$ equipped with $\kappa: \cE^2\,\iso\, \cO_{S\times X}(-D)$. Then $\cO_{S\times X}\oplus \cE$ is a $\cO_{S\times X}$-algebra, and $_SY=\Spec (\cO_{S\times X}\oplus \cE)$

  We simply think of $\Bun_{k,r}$ as the stack classifying $D_X\in{^{rss}X^{(r)}}$, a two-sheeted covering $\phi: Y\to X$ ramified exactly at $D_X$ (with $Y$ smooth), and a rank $k$ vector bundle $U$ on $Y$. 
  
  Remind that we assume $n=2$, so $G=\GSp_4$. Note that $^{rss}\cS^r_P=\Bun_{1,r}$. We have a diagram
$$
\Bun_4 \getsup{\gq_2} \Bun_{2,r} \toup{\gp_2}  {^{rss}\cS^r_P},
$$
where $\gq_2$ (resp., $\gp_2$) is the map sending $(\phi: Y\to X, U)$ as above to $W=\phi_*U$ (resp., to the point $(\phi: Y\to X, (\det U)^{\star})$ of ${^{rss}\cS^r_P}$). 
Extend it to a commutative diagram
\begin{equation}
\label{diag_GO_6_GSp_4}
\begin{array}{ccccc}
\Bun_4 &\getsup{\gq_2} &\Bun_{2,r}\\
\downarrow\lefteqn{\scriptstyle \rho} && \downarrow & \searrow\lefteqn{\scriptstyle \gp_2}\\
\Bun_{\tilde H} & \getsup{\gq_R} & \Bun_R & \toup{\gp_R} & 
{^{rss}\cS^r_P}
\end{array}
\end{equation}
defined as follows. 

 For a point of $\RCov^r$ given by $\phi:Y\to X$ let $R_{\phi}$ now denote the group scheme on $X$ included into an exact sequence $1\to \mu_2\to \phi_*\GL_2\to R_{\phi}\to 1$. Let $_{\phi}\GL_4$ be the group scheme of automorphisms of $\phi_*(\cO^2_Y)$. Define $_{\phi}\tilde H$ by the commutative diagram, where the rows are exact sequences
$$
\begin{array}{ccccccc}
1\to & \mu_2 & \to & \phi_*\GL_2 & \to & R_{\phi} &\to 1\\
       & \downarrow\lefteqn{\scriptstyle\id} && \downarrow&&
\downarrow\\
1\to & \mu_2 & \to & _{\phi}\GL_4 & \to & _{\phi}\tilde H & \to 1
\end{array}
$$
Since $_{\phi}\tilde H$ is an inner form of $\tilde H$, $\Bun_{\tilde H}\,\iso\, \Bun_{_{\phi}\tilde H}$ canonically. Let $\gq_{R_{\phi}}: \Bun_{R_{\phi}}\to \Bun_{\tilde H}$ denote the corresponding extension of scalars map. Since $\mu_2$ lies in the kernel of the determinant map $\phi_*\GL_2\toup{\det}\phi_*\Gm$, it yields
a map $R_{\phi}\to \phi_*\Gm$. Let $\gp_{R_{\phi}}: \Bun_{R_{\phi}}\to\Pic Y$ denote the composition of the corresponding extension of scalars map $\Bun_{R_{\phi}}\to\Pic Y$ with the automorphism $\epsilon:\Pic Y\,\iso\,\Pic Y$ sending $\cB$ to $\cB^{\star}$. 
As $\phi$ runs through $\RCov^r$, the group schemes $R_{\phi}$ organize into a group scheme $R$ over $X\times\RCov^r$, and the diagrams
$$
\Bun_{\tilde H} \getsup{\gq_{R_{\phi}}} \Bun_{R_{\phi}}\toup{\gp_{R_{\phi}}} \Pic Y
$$
form a family giving rise to (\ref{diag_GO_6_GSp_4}). Remind the functor $F_{\cS}$ introduced in Section~3.3.2.  

\begin{Pp}  
\label{Pp_GWP_and_BP}
For $K\in\D(\Bun_{\tilde H})$ there is a functorial isomorphism
$$
F_{\cS}(K)\mid_{^{rss}\cS^r_P}\,\iso\, (\gp_R)_!\gq_R^*K[\dimrel(\gq_R)]
$$
\end{Pp}
\begin{Prf} 
\Step 1
Define the stack $\cV_{4,P}$ and the maps $\tilde\gq$, $\tilde\gp$ by the diagram
$$
\begin{array}{ccccc}
\Bun_4 & \getsup{\tilde\gq} & \cV_{4,P} \\
\downarrow\lefteqn{\scriptstyle \rho} && \downarrow & \searrow\lefteqn{\scriptstyle \tilde\gp}\\
\Bun_{\tilde H} & \getsup{\gq_{\cV}} & \cV_{\tilde H,P} & \toup{\gp_{\cV}} & \cS_P,
\end{array}
$$
where the square is cartesian. The stack $\cV_{4,P}$ classifies: $L\in\Bun_2$, $W\in\Bun_4$ and a map $t: L\otimes V\to\Omega$ with $V=\wedge^2 W$. 

 Let $\zeta_2: \Bun_{2,r}\to \cV_{4,P}$ be the map sending $(\phi:Y\to X, U)$ to $(W=\phi_*U, L=\phi_*((\wedge^2 U)^{\star}), t)$, where $t:V\to L^*\otimes\Omega$ is the following map. We have 
$L^*\otimes\Omega\,\iso\, \phi_*(\wedge^2 U)$. The exterior square of the natural map $\phi^*\phi_* U\to U$ is a map $\phi^*(\wedge^2 W)\to \wedge^2 U$,  by adjointness it yields a map $t:V\to \phi_*(\wedge^2 U)$.

 We have a commutative diagram
\begin{equation}
\label{diag_sect7_cartesian}
\begin{array}{ccccc}
 && \Bun_{2,r} & \toup{\gp_2} & {^{rss}\cS^r_P}\\
 & \swarrow\lefteqn{\scriptstyle \gq_2} & \downarrow\lefteqn{\scriptstyle \zeta_2} && \downarrow\\
 \Bun_4 & \getsup{\tilde\gq} & \cV_{4,P} & \toup{\tilde\gp} & \cS_P
 \end{array}
\end{equation}
Let us show that the square in this diagram is cartesian. 
To do so, consider a $k$-point $(L, W, t)$ of $\cV_{4,P}$ whose image under $\tilde\gp$ is given by a $k$-point $(\phi: Y\to X, \cB, D_X)$ of $^{rss}\cS^r_P$. So, we are given isomorphisms $L\,\iso\,\phi_*\cB$, $\cC\,\iso\,\det W$, and $N(\cB)\,\iso\, \Omega^2\otimes \cC^{-1}$ such that the diagram commutes
\begin{equation}
\label{diag_sect_7_V4}
\begin{array}{cccc}
\Sym^2  V & \toup{t\otimes t} & \Sym^2(L^*\otimes\Omega)\\
\uparrow && \uparrow\\
\cC & \iso  & N(\cB^{-1})\otimes\Omega^2 & \,\iso\, N(\cB^{\star})(-D_X)
\end{array}
\end{equation}
Write $\sigma_{\phi}$ for the nontrivial automorphism of $Y$ over $X$, let $D_Y$ be the ramification divisor of $\phi:Y\to X$, so $D_X=\phi_*(D_Y)$. 

 The map $t: V\to L^*\otimes\Omega\,\iso\, \phi_*(\cB^{\star})$ can be seen as $t:\phi^*V\to \cB^{\star}$. The latter map is nonzero, because the symmetric form on $L$ is generically nondegenerate. Applying $\phi^*$ to (\ref{diag_sect_7_V4}), we get a commutative diagram
$$
\begin{array}{ccc}
\phi^*\Sym^2 V & \to  & (\cB^{\star})^2\oplus (\sigma_{\phi}^*\cB^{\star})^2\oplus (\cB^{\star}\otimes \sigma_{\phi}^*\cB^{\star})\\
\uparrow && \uparrow\lefteqn{\scriptstyle (0,0,1)}\\
\phi^*\det W & \iso & (\cB^{\star}\otimes \sigma_{\phi}^*\cB^{\star})(-2D_Y)
\end{array}
$$

 The transpose $\cB\otimes\Omega^{-1}\to \phi^* V^*$ to $t$ is an isotropic subscheaf in $\phi^*V^*$. So, there is a rank 2 vector bundle $U$ on $Y$ and  a surjection $\phi^* W\toup{v} U$ such that $t$ factors as a composition
$$
\phi^*V\to \wedge^2 U\hook{u} \cB^{\star}
$$ 
We are going to check that $\wedge^2 U\hook{u}\cB^{\star}$ is actually and isomorphism, and the map $W\toup{v} \phi_*U$ is also an isomorphism. 

 Indeed, the maps $V\to \phi_*(\wedge^2 U)\to \phi_*(\cB^{\star})$ yield a commutative diagram
$$
\begin{array}{ccccc}
\Sym^2(V^*\otimes\Omega) & \gets & \Sym^2(\phi_*((\wedge^2 U)^{\star})) & \gets & \Sym^2 L\\
\downarrow && \downarrow && \downarrow\\
\cC^{-1}\otimes\Omega^2 & \gets & N((\wedge^2 U)^{\star}) & \hookleftarrow & N(\cB),
\end{array}
$$
and the composition of maps in the low row is an isomorphism. It follows that the transpose $\cB\to (\wedge^2 U)^{\star}$ to $u$ is an isomorphism. 

 Now consider the diagram
$$
V\toup{\wedge^2 v} \wedge^2(\phi_* U)\to \phi_*(\wedge^2 U),
$$
where the second map is induced by the natural map $\phi^*\phi_*U\to U$. Applying symmetric squares, one gets a commutative diagram
$$
\begin{array}{ccccc}
\cC & \iso & N(\wedge^2 U)(-D_X)\\
\downarrow &&& \searrow\\
\Sym^2 V & \to & \Sym^2(\wedge^2(\phi_* U)) & \to & \Sym^2(\phi_*(\wedge^2 U))\\
\downarrow && \downarrow && \downarrow\\
\cC & \to & \det(\phi_* U) & \hook{\xi} & N(\wedge^2 U)
\end{array}
$$
It is easy to see that $\xi$ induces an isomorphism $\det(\phi_* U)\,\iso\, N(\wedge^2 U)(-D_X)$. Thus, $\det v: \det W\,\iso\, \cC\,\iso\, \det(\phi_*U)$ is an isomorphism, so $v: W\to \phi_* U$ is an isomorphism. 

\medskip
\Step 2  The diagram (\ref{diag_sect7_cartesian}) gives rise to a commutative diagram
$$
\begin{array}{ccccc}
&& \Bun_R & \toup{\gp_R} & {^{rss}\cS^r_P}\\
& \swarrow\lefteqn{\scriptstyle \gq_R} & \downarrow && \downarrow\\
\Bun_{\tilde H} & \getsup{\gq_{\cV}} & \cV_{\tilde H,P} & \toup{\gp_{\cV}} & \cS_P
\end{array}
$$
We have to show that the square in this diagram is cartesian.
By Step 1, this is true after the base change $_0\Bun_{\tilde H}\hook{}\Bun_{\tilde H}$. For the components of $_1\Bun_{\tilde H}$ the argument is similar. 
\end{Prf}

\begin{Def} Let $K\in\D(\Bun_{\tilde H})$ be a complex with central character $\chi^{-1}$. Then $F_{\cS}(K)$ has central character $\chi$. Let $\cJ$ be a rank one local system on $Y$ equipped with $N(\cJ)\,\iso\, \chi^{-1}$. Then for $\Pic Y\hook{} {^{rss}\cS^r_P}\subset \cS_P$ the $*$-restriction 
$A\cJ\otimes F_{\cS}(K)\mid_{\Pic Y}$ is equipped with natural descent data for $e_{\phi}:\Pic Y\to\Bun_{U_{\phi}}$. Assume that
\begin{itemize}
\item[($C_{G}$)]  $\cK_K$ is a complex on $\Bun_{U_{\phi}}$ equipped with 
$$
e_{\phi}^*\cK_K[\dimrel(e_{\phi})]\,\iso\, A\cJ\otimes (\gp_{R_{\phi}})_!\gq_{R_{\phi}}^*(K)[\dimrel(\gq_{R_{\phi}})]
$$
\end{itemize}
For $a\in\ZZ/2\ZZ$ the \select{generalized Waldspurger period} of $K$ is 
$$
\GWP^a(K,\cJ)=\RG_c(\Bun^a_{U_{\phi}}, \cK_K)
$$
\end{Def}

\bigskip\bigskip
\centerline{\scshape 8. Towards a construction of automorphic sheaves on $\Bun_{GSp_4}$}

\bigskip\noindent
8.1 In Section~8 we assume the ground field $k$ algebraically closed of characteristic zero and work with  $\cD$-modules instead of $\ell$-adic sheaves. 
A local system on $X$ is now a vector bundle $E$ with connection $\nabla: E\to E\otimes\Omega$. 
 
\begin{Lm} 
\label{Lm_det_RG_DR}
For a local system $E$ on $X$ there is a canonical $\ZZ/2\ZZ$-graded isomorphism
$$
\det\RG_{DR}(X,E)\,\iso\, \det\RG(X, \det E)\otimes\det\RG(X, \Omega\otimes\det E)^{-1}
$$
\end{Lm} 
\begin{Prf}
Let $\DR(E)=(E\toup{\nabla} E\otimes\Omega)$ be the De Rham complex of $E$ placed in degrees 0 and 1.
The exact triangle $\DR(E)\to E\to E\otimes\Omega$ yields 
$$
\det\RG_{DR}(X,E)\,\iso\, \det\RG(X,E)\otimes\det\RG(X, E\otimes\Omega)^{-1}
$$ 
Our assertion follows now from Lemma~\ref{Lm_general_detRG}.
\end{Prf}  
 
\medskip  
  
 Set $n=2$, so $G=\GSp_4$. Write $\LocSys_{\check{G}}$ for the moduli stack of $\check{G}$-local systems on $X$. View a $\check{G}$-local system $E_{\check{G}}$ as a pair $(E,\chi)$, where $E$ (resp., $\chi$) is a rank 4 (resp., rank 1) local system on $X$ with symplectic form $\wedge^2 E\to\chi^{-1}$. 
 
 Let $\phi:Y\to X$ be a (possibly ramified) two-sheeted covering. Write $\LocSys_{Y,r}$ for the moduli stack of rank $r$ local systems on $Y$. 
 
 Let $\cM_Y$ denote the stack classifying: a $\check{G}$-local system $E_{\check{G}}=(E,\chi)$ on $X$, $\cJ\in \LocSys_{Y,1}$ equipped with $N(\cJ)\,\iso\,\chi$. The following is an immediate consequence of Lemma~\ref{Lm_det_RG_DR}. 
 
\begin{Lm} 
\label{Lm_second_detRG_DR}
The ($\ZZ/2\ZZ$-graded) line bundle on $\cM_Y$ with fibre $\det\RG_{DR}(Y, \cJ\otimes\phi^*E)$ at $(\cJ, E_{\check{G}})$ is canonically trivialized. \QED
\end{Lm}
 
 Let $^0\cM_Y\subset \cM_Y$ denote the open substack given by the condition $\H^i_{DR}(Y, \cJ\otimes\phi^*E)=0$ for $i=0,2$. We have a vector bundle $V$ on $^0\cM_Y$ whose fibre at $(\cJ, E_{\check{G}})$ is
$$
\H^1_{DR}(Y, \cJ\otimes\phi^* E)
$$
The rank of $V$ is $c:=8(g_Y-1)$, where $g_Y$ is the genus of $Y$.  As in Remark~\ref{Rem_sym_form_on_V}, one equips $V$ with a nondegenerate symmetric form $\Sym^2 V\to \cO$. Moreover, the trivialization $\det V\,\iso\,\cO$ given by Lemma~\ref{Lm_second_detRG_DR} is compatible with this symmetric form. 
 
\begin{Con} 
\label{Con_lifting_Spin_D-modules}
The $\SO_c$-torsor $V$ lifts naturally to a $\Spin_c$-torsor $\cF$ on $^0\cM_Y$.
\end{Con}
 
\smallskip\noindent
8.2 Let $E_{\check{G}}$ be a $\check{G}$-local system on $X$ viewed as a pair $(E,\chi)$, where $E$ and $\chi$ are local systems on $X$ of ranks $4$ and $1$ respectively, and $\wedge^2 E\to \chi^{-1}$ is a symplectic form.  Assume $E$ irreducible. In this situation we propose the following conjectural construction of an automorphic sheaf $K_{E_{\check{G}}}$ on $\Bun_G$ corresponding to $E_{\check{G}}$. 

 Let $r\ge 0$, remind the stacks $\RCov^r$ and $^{rss}\cS^r_P\subset \cS^r_P$ (cf. Section~6.1.1). A point of $\RCov^r$ is given by a divisor $D_X\in {^{rss}X^{(r)}}$ and a two-sheeted covering $\phi: Y\to X$ ramified exactly at $D_X$ with $Y$ smooth. Let $Y_{univ}\to X\times \RCov^r$ denote the universal two-sheeted covering. For a morphism of stacks $\alpha: S\to\RCov^r$ denote by $Y_S\to X\times S$ the two-sheeted covering obtained from the universal one by the base change $\id\times\alpha: X\times S\to X\times\RCov^r$. 
 
 Let $\cM$ be the stack classifying collections: a point of $\RCov^r$ given by $D_X\in {^{rss}X^{(r)}}$, $Y\toup{\phi} X$, and a rank one local system $\cJ$ on $Y$ equipped with an isomorphism $N(\cJ)\,\iso\,\chi$.  By definition, for a scheme $S$, an $S$-point of $\cM$ is given by a map $S\to\RCov^r$ and a rank one local system $\cJ$ (relative to $S$) over $Y_S$ equipped with a trivialization of $N_{X\times S}(\cJ)$. 

 Let $V$ be the vector bundle on $\cM$ whose fibre at the above point is $\H^1_{DR}(Y, \cJ\otimes\phi^* E)$. As in 8.1, we equip it with a nondegenerate symmetric form and a compatible trivialization $\det V\,\iso\,\cO$. 
Assuming Conjecture~\ref{Con_lifting_Spin_D-modules}, we get a $\Spin_c$-torsor on $\cM$. For the half-spin fundamental weights $\alpha_a$ ($a\in\ZZ/2\ZZ$) of $\Spin_c$ write $V^{\alpha_a}$ for the corresponding vector bundles on $\cM$ induced from our $\Spin_c$-torsor. 

 The projection $\cM\to\RCov^r$ should be equipped with an integrable connection along $\RCov^r$ making $\cM$ into a $\cD_{\RCov^r}$-stack (in the sense of \cite{BD}, Section~2.3.1). Then $\cM$ carries a sheaf of algebras $\cO_{\cM}[\cD_{\RCov^r}]$ (in the notation of \select{loc.cit.}, Section~2.3.4). We expect that $V^{\alpha_a}$ is naturally a module over $\cO_{\cM}[\cD_{\RCov^r}]$. 

 Consider the two-sheeted covering 
$Y_{\cM}\to X\times\cM$ 
obtained from $Y_{univ}\to X\times\RCov^r$ by the base change $\id\times \pr: X\times\cM\to X\times\RCov^r$. Let $\cJ_{univ}$ denote the universal local system (relative to $\cM$) over $Y_{\cM}$, its norm on $X\times\cM$ is trivialized. Let $Y_{\cM}^{(d)}$ denote the $d$-th symmetric power of $Y_{\cM}$ (relative to $\cM$). That is, $Y_{\cM}^{(d)}$ is the quotient of the $d$-th power $Y_{\cM}\times_{\cM} \ldots\times_{\cM} Y_{\cM}$ by the symmetric group on $d$ elements. Let $\cJ_{univ}^{(d)}$ denote the corresponding local system (relative to $\cM$) on $Y_{\cM}^{(d)}$. 
  
 Remind that, for a scheme $S$, an $S$-point of $^{rss}\cS^r_P$ is given by a map $S\to\RCov^r$ and an invertible sheaf $\cB$ on $Y_S$. For $a\in\ZZ/2\ZZ$ write $^{rss}\cS^r_{a,P}$ for the substack of $^{rss}\cS^r_P$ given by $a=(\deg\cB)\!\mod 2$. 
 
  An $S$-point of $Y^{(d)}_{\cM}$ is given by a collection: a map $S\to\RCov^r$, a rank one local system $\cJ$ (relative to $S$) over $Y_S$ with a trivialization of $N_{X\times S}(\cJ)$, and an effective Cartier divisor $D_S$ on $Y_S$ flat over $S$ of degree $d$. For $d\ge 0$ consider the Abel-Jacobi map
$$
\jac: Y^{(d)}_{\cM} \to \cM\times_{\RCov^r}{^{rss}\cS^r_P}
$$
over $\cM$, it is given by $\cB=\cO_{Y_S}(D_S)$.  
 
 There is a unique local system $\cP$ (relative to $\cM$) over $\cM\times_{\RCov^r}{^{rss}\cS^r_P}$ with the following properties. For any $d\ge 0$ one has $\jac^*\cP\,\iso\, \cJ^{(d)}_{univ}$ canonically, and $\cP$ satisfies the usual automorphic property with respect to the group structure of $\Pic Y$. More precisely, $\cM\times_{\RCov^r}{^{rss}\cS^r_P}$ is a commutative group stack over $\cM$, and the automorphic property of $\cP$ is required for this group structure. 

 In more concrete terms, $\cM\times_{\RCov^r}{^{rss}\cS^r_P}$ classifies: $D_X\in {^{rss}X^{(r)}}$, $\phi: Y\to X$, a rank one local system $\cJ$ on $Y$ with a trivialization of $N_X(\cJ)$, and an invertible sheaf $\cB$ on $Y$. Then the fibre of $\cP$ at this point identifies with $(A\cJ)_{\cB}$. 
 
 Consider the diagram of projections
$$
\cM \;\getsup{q_{\cM}}  \;\cM\times_{\RCov^r}{^{rss}\cS^r_P} \; \toup{q_{\cS}}\;
{^{rss}\cS^r_P}
$$
For $a\in\ZZ/2\ZZ$ define a complex $K_a$ on $^{rss}\cS^r_{a,P}$ by
$$
(q_{\cS})_*(q_{\cM}^*V^{\alpha_a}\otimes \cP)[\dim],
$$ 
where the direct image with respect to $q_{\cS}$ is understood in the (derived) quasi-coherent sense. 

 We expect that, for a suitable shift, $K_a$ has a natural structure of a $\cD$-module on $^{rss}\cS^r_P$. Let then $K$ be the $\cD$-module on $^{rss}\cS^r_P$ whose restriction to $^{rss}\cS^r_{a,P}$ is $K_a$. Let $\tilde K$ denote the intermediate extension of $K$ under $^{rss}\cS^r_P\subset \cS^r_P$. 
 
 Remind that $\cS_P$ and $\Bun_P$ are dual (generalized) vector bundles over $\Bun_2\times\Pic X$, let $\Four(\tilde K)$ denote the Fourier transform of $\tilde K$. Remind the projection $\nu_P:\Bun_P\to\Bun_G$, let $^0\Bun_P\subset\Bun_P$ be the open substack, where $\nu_P$ is smooth. We expect that $\nu_P^*K_{E_{\check{G}}}$ idenitifies over $^0\Bun_P$ with $\Four(\tilde K)$ as $\cD$-modules.

\bigskip\bigskip
\centerline{\scshape Appendix A. Prym varieties}

\bigskip\noindent
A.1 Let $\pi:\tilde X\to X$ be a two-sheeted covering ramified at some divisor $D_{\pi}$ on $X$ with $\deg D_{\pi}=d$. Let $\sigma$ be the nontrivial automorphism of $\tilde X$ over $X$. Let $\cE$ be the $\sigma$-anti-invariants in $\pi_*\cO$, it is equipped with $\cE^2\,\iso\,\cO(-D_{\pi})$. Let $\cE_0$ be the 
$\sigma$-anti-invariants in $\pi_*\Qlb$, it is equipped with $\cE_0^2\,\iso\,\Qlb$ over $X-D_{\pi}$. 
 
 The norm map $N:\Pic\tilde X\to\Pic X$ is given by $N(\cB)=\cE^{-1}\otimes\det(\pi_*\cB)$, this is a homomorphism of group stacks. We write $N_X(\cB)=N(\cB)$ when we need to express the dependence on $X$.
For $\cC\in\Pic X$ we have canonically $N(\pi^*\cC)\,\iso\,\cC^2$. 
 
 Let $\tilde E$ be a rank one local system on $\tilde X$. Then $\tilde E\otimes\sigma^*\tilde E$ is equipped with natural descent data for $\pi$, so there is a rank one local system $N(\tilde E)$ on $X$ equipped with $\pi^*N(\tilde E)\,\iso\,  \tilde E\otimes\sigma^*\tilde E$. (In the nonramified case we have $N(\tilde E)\,\iso\, \cE_0\otimes\det(\pi_*\tilde E)$ canonically).
Remind that $A\tilde E$ denotes the automorphic local system on $\Pic\tilde X$ corresponding to $\tilde E$. The restriction of $A\tilde E$ under $\pi^*:\Pic X\to\Pic\tilde X$ identifies canonically with $AN(\tilde E)$. For a rank one local system $E$ on $X$ we have canonically $N^*(AE)\,\iso\, A(\pi^* E)$. 
 
 Write $\uPic X$ for the Picard scheme of $X$, so we have a $\Gm$-gerbe $\Pic X\to\uPic X$. 
Write $\Pic^r X$ for the connected component classifying line bundles of degree $r$. Let $\und{P}$ denote the connected component of unity of $\Ker \und{N}$, where $\und{N}: \uPic\tilde X\to\uPic X$ is the norm map. This is the Prym variety (\cite{M}). We need the following results proved in \select{loc.cit}. Assume $\tilde X$ connected. 
 
\medskip\noindent
 \select{Case of ramified $\pi$}. The group scheme $\Ker \und{N}$ is connected, $\pi^*: \uPic X\hook{}\uPic\tilde X$ is a closed immersion. For each $r$ we have a surjection $$
\uPic^{2r}\tilde X/\uPic^r X\to\und{P}
$$ 
sending $\cB$ to $\cB^{-1}\otimes\sigma^*\cB$. For $r=0$ its kernel is a finite group isomorphic to $P_2/\phi(J_2)$ for some inclusion $\phi: J_2\to P_2$. Here $P_2$ and $J_2$ are the groups of order 2 points of $\und{P}$ and $\uPic^0 X$ respectively. Remind that $\dim\uPic^0 X=g$, $\dim\uPic^0\tilde X=2g+\frac{d}{2}-1$ and $\dim\und{P}=g+\frac{d}{2}-1$. So, 
$$
J_2\,\iso\, (\ZZ/2\ZZ)^{2g}\hskip 1em \mbox{and}\hskip 1 em
P_2\,\iso\, (\ZZ/2\ZZ)^{2g+d-2}
$$ 

\smallskip\noindent 
\select{Case of unramified $\pi$}. The group scheme $\Ker \und{N}$ has two connected components, say $\Ker^r \und{N}$ for $r\in\ZZ/2\ZZ$. We denote by $\Ker^0\und{N}$ the connected component of unity. The kernel of $\pi^*: \uPic X\to \uPic\tilde X$ is $H_0:=\{\cO, \cE\}$, and we have an isomorphism $\uPic^0\tilde X/(\uPic^0 X/ H_0)\,\iso\, \und{P}$ sending $\cB$ to $\cB^{-1}\otimes\sigma^*\cB$.
 
 The following is probably well-known, but we have not found a proof of it, so we give one.
 
\begin{Lm} Assume $\pi$ nonramified. Both $\uPic\tilde X/(\uPic X/H_0)$ and $\Ker\und{N}$ have two connected components indexed by $\ZZ/2\ZZ$. The map $\delta: \uPic\tilde X/(\uPic X/H_0)\to \Ker\und{N}$ sending $\cB$ to  $\cB^{-1}\otimes\sigma^*\cB$ is an isomorphism, so 
sends the odd connected component to the odd one.
\end{Lm}
\begin{Prf}
Consider the map $\pi_*\Gm\toup{\xi}\pi_*\Gm$ sending $f$ to $\sigma(f)f^{-1}$. Let $\cN:\pi_*\Gm\to\Gm$ be the norm map. The sequences 
$$
\pi_*\Gm\toup{\xi}\pi_*\Gm\toup{\cN}\Gm\to 1
$$ 
and
$$
1\to \Gm\to\pi_*\Gm\toup{\xi}\pi_*\Gm
$$
are exact in \'etale topology (indeed, it suffices to check this after base change $\tilde X\toup{\pi} X$, what is easy). Taking the \'etale cohomology of $X$, we get an exact sequence
$$
\H^1(X,\Im\xi)\to \uPic\tilde X\toup{\und{N}} \uPic X\to 1
$$ 
and the map induced by $\xi$ is $\delta$. The map $\uPic\tilde X\,\iso\,\H^1(X, \pi_*\Gm)\toup{\xi} \H^1(X,\Im\xi)$ is surjective, because $\H^2(X,\Gm)=0$. 
\end{Prf}

\medskip

 Let $\tilde D_{\pi}\in \tilde X^{(d)}$ be the ramification divisor of $\pi$, so $D_{\pi}=\pi_*\tilde D_{\pi}$. Define the group scheme $U_{\pi}$ on $X$ by the exact sequence $1\to \Gm\to \pi_*\Gm\to U_{\pi}\to 1$. The stack $\Bun_{U_{\pi}}$ classifies $\cB\in\Pic\tilde X$ together with a trivialization 
\begin{equation}
\label{triv_appendix_A} 
N(\cB)\,\iso\, \cO_X
\end{equation}
and a compatible isomorphism $\gamma: \cB\mid_{\tilde D_{\pi}}\,\iso\cO_{\tilde D_{\pi}}$. This means that $\gamma^{\otimes 2}: N(\cB)\mid_{D_{\pi}}\,\iso\, \cO_{D_{\pi}}$ is the isomorphism 
induced by (\ref{triv_appendix_A}).
We have used that $\pi$ induces an isomorphism of reduced divisors $\tilde D_{\pi}\,\iso\, D_{\pi}$.  
The corresponding extension of scalars map $e_{\pi}: \Pic \tilde X\to \Bun_{U_{\pi}}$ is smooth and surjective, it sends $\cB$ to $\cC=\cB^{-1}\otimes \sigma^*\cB$ with natural trivializations $\cC\mid_{\tilde D_{\pi}}\,\iso\, \cO_{\tilde D_{\pi}}$ and $N(\cC)\,\iso\,\cO_X$.  
  
  In both cases ($\pi$ ramified or not) the stack $\Bun_{U_{\pi}}$ has two connected components indexed by $a\in\ZZ/2\ZZ$, here $\Bun_{U_{\pi}}^0$ is the connected component of unity. The image of $\Pic^r \tilde X$ under $e_{\pi}$ equals $\Bun_{U_{\pi}}^a$ with $a=r\mod 2$. 
  
  If $\pi$ is ramified then $\Bun_{U_{\pi}}$ is a scheme, and the projection $\Bun^a_{U_{\pi}}\to\und{P}$ is a Galois covering with Galois group $(\ZZ/2\ZZ)^{d-2}$. 
 
\bigskip\bigskip
\centerline{Appendix B. Group schemes and Hecke operators}

\bigskip\noindent
B.1.1 Let $\pi:\tilde X\to X$ be an \' etale Galois covering with Galois group $\Sigma=\Aut_X(\tilde X)$. The category of affine group schemes over $X$ is canonically equivalent to the category of affine group schemes over $\tilde X$ equipped with an action of $\Sigma$.

 Assume that $\GG$ is an affine group scheme over $\Spec k$ (viewed as constant group scheme on $\tilde X$). Then an action of $\Sigma$ on $\GG$ is a datum of a homomoprhism $\Sigma\to \Aut(\GG)$. The corresponding group scheme $G$ over $X$ is then obtained as the twisting of $\GG$ by the $\Sigma$-torsor $\pi:\tilde X\to X$. 
 
 The action of $\Sigma$ on $\GG$ gives rise to the semi-direct product $\GG\ltimes \Sigma$ included into an exact sequence $1\to \GG\to \GG\ltimes \Sigma\to\Sigma\to 1$. 
 
 Let us describe the category of $G$-torsors on $X$. For $\sigma\in\Sigma$ and a $\GG$-torsor $\cF_{\GG}$ on a scheme $S$ denote by $\cF_{\GG}^{\sigma}$ the $\GG$-torsor on $S$ obtained from $\cF_{\GG}$ by the extension of scalars $\sigma: \GG\to \GG$. 

\begin{Lm} Let $S$ be a scheme. The category of $G$-torsors on $S\times X$ is canonically equivalent to the category of pairs $(\cF_{\GG}, \alpha)$, where $\cF_{\GG}$ is a $\GG$-torsor on $S\times\tilde X$, and $\alpha=(\alpha_{\sigma})_{\sigma\in\Sigma}$ is a collection of isomorphisms
$$
\alpha_{\sigma}: \sigma^*\cF_{\GG}\,\iso\, \cF_{\GG}^{\sigma}
$$
such that for any $\sigma,\tau\in\Sigma$ the diagram commutes
$$
\begin{array}{ccc}
\tau^*\sigma^*\cF_{\GG} & \toup{\alpha_{\sigma}} \tau^*(\cF_{\GG}^{\sigma})= & (\tau^*\cF_{\GG})^{\sigma}\\ 
\downarrow\lefteqn{\scriptstyle \alpha_{\sigma\tau}} && \downarrow\lefteqn{\scriptstyle \alpha_{\tau}}\\
(\cF_{\GG})^{\sigma\tau} & = & (\cF^{\tau}_{\GG})^{\sigma}
\end{array}
$$
\end{Lm}
\begin{Prf}
Let $\tilde F$ be an affine scheme over $\tilde X$. Assume that the action of $\Sigma$ on $\tilde X$ is lifted to an action on $\tilde F$. Let $F$ denote the affine scheme over $X$, the descent of $\tilde F$. 

 Assume that $\tilde F$ is, in addition, a $\GG$-torsor. 
Then for $F$ to be a $G$-torsor, the actions of $\GG$ and of $\Sigma$ on $\tilde F$ should come from an action of $\GG\ltimes \Sigma$ on $\tilde F$. 
 \end{Prf}

\medskip\noindent
Example 1. Take $\HH$ to be a group scheme over $\Spec k$, set $\GG=\Hom(\Sigma, \HH)$ (the group structure on $\GG$ comes from that of $\HH$). Let $\Sigma$ act on $\GG$ via its action on $\Sigma$ by translations. Then $G\,\iso\, \pi_*\HH$. 

\medskip\noindent
B.1.2 Let $\GG$ be a connected reductive group over $\Spec k$ equipped with a homomoprhism $\Sigma\to \Aut(\GG)$. Let $G$ be the twisting of $\GG$ by the $\Sigma$-torsor $\pi:\tilde X\to X$. Assume that $\TT\subset \GG$ is a maximal torus invariant under $\Sigma$. Let $\Lambda$ (resp., $\check{\Lambda}$) be the coweights (resp., weights) lattices of $\TT$, so $\Sigma$ acts on the corresponding root datum $(\Lambda, R, \check{\Lambda}, \check{R})$. Here $R$ (resp. $\check{R}$) are the coroots (resp., roots) of $\GG$.

 Let $W=N_{\GG}(\TT)/\TT$ be the Weil group. Since $\Sigma$ preserves both $\TT$ and $N_{\GG}(\TT)$, $\Sigma$ acts on $W$ by group automorphisms. For $\lambda\in\Lambda$, $w\in W$ and $\sigma\in\Sigma$ we have $\sigma(w\lambda)=(\sigma w)(\sigma\lambda)$, so $\Sigma$ acts on the set $\Lambda/W$ of dominant coweights of $\GG$. 

Write $\Bun_G$ for the stack of $G$-torsors on $X$. Let $\cH_G$ denote the Hecke stack classifying: $\tilde x\in\tilde X$, $G$-torsors $\cF_G,\cF'_G$ on $X$, and an isomorphism $\cF_G\,\iso\,\cF'_G\mid_{X-\pi(\tilde x)}$. We have a diagram
$$
\tilde X\times \Bun_G \getsup{\supp\times p}\cH_G\toup{p'}\Bun_G,
$$ 
where $p$ (resp., $p'$) sends the above point to $\cF_G$ (resp., to $\cF'_G$). 

 A choice of a Borel subgroup in $\GG$ containing $\TT$ identifies $\Lambda/W$ with the corresponding set of dominant coweights, hence yields a usual order on $\Lambda/W$. This order does not depend on a choice of such Borel subgroup. 

 Let $D_x$ (resp., $D_{\tilde x}$ denote the formal neighbourhood of $x\in X$ (resp., of $\tilde x\in\tilde X$).  The map $\pi$ induces $D_{\tilde x}\,\iso\, D_x$ for $x=\pi(\tilde x)$. For $\lambda\in\Lambda/W$ denote by $\ov{\cH}_G^{\lambda}\hook{} \cH_G$ the closed substack given by the condition that $\cF'_G\mid_{D_{\tilde x}}$ is in a position $\le\lambda$ w.r.t. $\cF_G\mid_{D_{\tilde x}}$, here we view them as $\GG$-torsors using the canonical isomorphism $G\mid_{\tilde X}\,\iso\,\GG\times X$. Given $\sigma\in\Sigma$, this condition is equivalent to requiring that $\cF'_G\mid_{D_{\sigma\tilde x}}$ is in a position $\le \sigma\lambda$ w.r.t. $\cF_G\mid_{D_{\sigma\tilde x}}$.
  
 The Hecke functor 
$$
\H^{\lambda}_G:\D(\Bun_G)\to \D(\tilde X\times\Bun_G)
$$ 
is given by 
$$
\H^{\lambda}_G(K)=(\supp\times p)_!((p')^*K\otimes \IC_{\ov{\cH}_G^{\lambda}})[-\dim\Bun_G]
$$
For each $\sigma\in\Sigma$ we have a commutative diagram
$$
\begin{array}{ccccc}
\tilde X\times \Bun_G  & \getsup{\supp\times p} & \ov{\cH}_G^{\sigma^{-1}\lambda} & \toup{p'} & \Bun_G\\
\downarrow\lefteqn{\scriptstyle \sigma\times\id} && \downarrow && \downarrow\lefteqn{\scriptstyle \id}\\
 \tilde X\times \Bun_G  & \getsup{\supp\times p} & \ov{\cH}_G^{\lambda} & \toup{p'} & \Bun_G,
\end{array}
$$
where the vertical middle arrow is an isomorphism. 
This yields a compatible system of isomorphisms for $\sigma\in\Sigma$
\begin{equation}
\label{iso_action_Sigma_Hecke_sectionD}
(\sigma\times\id)^*\comp \H^{\lambda}_G\,\iso\, \H^{\sigma^{-1}\lambda}_G
\end{equation}

\noindent
Example 2. Given a homomorphism $\Sigma\to N_{\GG}(\TT)\subset \GG$, consider the corresponding action of $\Sigma$ on $\GG$ by conjugation. Then $G$ is an inner form of $\GG$, and $\Sigma$ acts trivially on $\Lambda/W$. Let $\cF^1_{\GG}$ be the $\GG$-torsor on $X$ obtained from the $\Sigma$-torsor $\tilde X$ by the extensioon of scalars via $\Sigma\to \GG$. Then $G$ identifies with the group scheme (over $X$) of automorphisms of the $\GG$-torsor $\cF^1_{\GG}$. 
In this case we identify canonically $\Bun_{\GG}\,\iso\,\Bun_G$ sending $\cF_{\GG}$ to the $G$-torsor $\Isom(\cF_{\GG}, \cF^1_{\GG})$. Then $\H^{\lambda}$ becomes the usual Hecke functor followed by restriction under $\pi\times\id: \tilde X\times\Bun_{\GG}\to X\times\Bun_{\GG}$. 

\bigskip\noindent
B.1.3 The map $\supp\times p: \ov{\cH}^{\lambda}_G\to \tilde X\times\Bun_G$ identifies with the twisted projection 
$$
(\tilde X\times\Bun_G)\ttimes \,\ov{\Gr}^{\lambda}_{\GG}\to \tilde X\times\Bun_G
$$
Similarly, $\supp\times p': \ov{\cH}^{\lambda}_G\to \tilde X\times\Bun_G$ identifies with $(\tilde X\times\Bun_G)\ttimes \, \ov{\Gr}^{-w_0(\lambda)}_{\GG}\to \tilde X\times\Bun_G$, where $w_0$ is the longuest element of $W$. Note that
$\IC_{\ov{\cH}^{\lambda}_G}$ is ULA with respect to $\supp\times p'$. This implies that $\H^{\lambda}_G$ commutes with the Verdier duality. 

\bigskip\noindent
B.1.4 Let $\check{\GG}$ denote the Langlands dual group to $\GG$, it comes equipped with a maximal torus $\check{\TT}$.  The group $\Sigma$ acts naturally on the root datum $(\check{\Lambda}, \check{R}, \Lambda, R)$ of $(\check{\GG}, \check{\TT})$. Remind that we have an exact sequence 
$$
1\to W\to \Aut(\check{\Lambda}, \check{R}, \Lambda, R)\to \Out(\GG)\to 1,
$$
where $\Out(\GG)$ is the group of exteriour automorphisms of $\GG$. Assume given a lifting of 
$$
\Sigma\to \Aut(\check{\Lambda}, \check{R}, \Lambda, R)
$$ 
to a homomorphism $\mu: \Sigma\to \Aut(\check{\GG}, \check{\TT})$. (Such lifting exists under the additional assumption that the $\Sigma$-action on $(\GG,\TT)$ preserves an epinglage of $\GG$ containing $\TT$). This lifting is uniquely defined up to inner automorphisms by elements of $\check{\TT}$. 
 
 Then we have the semi-direct product $G^L:=\check{\GG}\ltimes \Sigma$ included into an exact sequence $1\to \check{\GG}\to \check{\GG}\ltimes \Sigma\to\Sigma\to 1$. This is a version of the $L$-group associated to $G_F$. Here $G_F$ denotes the restriction of the group scheme $G$ to the generic point $\Spec F\in X$ of $X$ (cf. \cite{B}).
 
\bigskip\noindent
B.2 Let now $\GG_1$ be another reductive connected group over $\Spec k$ equipped with an action $\Sigma\to \Aut(\GG_1)$, let $G_1$ be the group scheme on $X$ obtained as the twisting of $\GG_1$ by the $\Sigma$-torsor $\pi:\tilde X\to X$. 

 Assume that $\GG_1$ satisfies the same conditions as $\GG$ in B.1. (The subscript $1$ denotes the corresponding objects for $G_1$). 
So, we have a maximal torus $\TT_1\subset \GG_1$ stable under $\Sigma$, and we assume given a homomorphism $\mu_1: \Sigma\to\Aut(\check{\GG}_1)$ as above. Assume given a $\Sigma$-equivariant homomorphism $\check{\GG}\to \check{\GG}_1$ sending $\check{\TT}$ to $\check{\TT}_1$. It yields a homomorphism $G^L\to  G_1^L$. 

 \select{The functoriality problem} is to find a family of functors 
$$
_SF: \D(S\times\Bun_G)\to \D(S\times\Bun_{G_1})
$$ 
for each scheme $S$ with the following property. Write $V^{\lambda_1}_1$ for the irreducible representation of $\check{\GG}_1$ with h.w. $\lambda_1\in \Lambda_1/W_1$. Similarly, $V^{\lambda}$ denotes the irreducible representation of $\check{\GG}$ with h.w. $\lambda\in \Lambda/W$ (this notion does not depend on a choice of a Borel subgroup in $\GG$ containing $\TT$). We would like to have for each $\lambda_1\in \Lambda_1/W_1$ isomorphisms of functors 
$$
\H_{G_1}^{\lambda_1}\comp {_SF}\,\iso\, \oplus_{\lambda} \;\, {_{\tilde X\times S}F}\comp \H^{\lambda}_G \otimes \Hom_{\check{\GG}} (V^{\lambda}, V^{\lambda_1}_1)
$$
from $\D(S\times\Bun_G)$ to $\D(\tilde X\times S\times\Bun_{G_1})$. It is required that these isomorphism are compatible with the action of $\Sigma$ on both sides. 
 
\bigskip\noindent
{\bf Acknowledgements.} It is a pleasure to thank G. Laumon for constant support. I am also grateful to V. Lafforgue and A. Braverman for nice and stimulating discussions.


\begin{thebibliography}{99}
\bibitem{A} J.~Adams, L-functoriality for dual pairs, Soc. Math. de France, Ast\'erique 171-172 (1989), p. 85-129
\bibitem{B} A. Borel, Automorphic L-functions, Proceedings of Symposia in Pure Math., vol. 33 (1979), part 2, p. 27-61
\bibitem{BD} A. Beilinson, V. Drinfeld, Chiral algebras, AMS Colloquium publications, 51, Providence, RI (2004)
\bibitem{BG} A. Braverman, D. Gaitsgory, Geometric Eisenstein series, Inv. Math. 150 (2002), 287-384
\bibitem{BFF} D. Bump, S. Friedberg, M. Furusawa, Explicit formulas for the Waldspurger and Bessel models, Israel J. Math. 102 (1997), 125--177
\bibitem{FGV} E. Frenkel, D. Gaitsgory, K. Vilonen, On the geometric Langlands conjecture. 
J. Amer. Math. Soc.  15  (2002),  no. 2, 367--417
\bibitem{FH} W. Fulton, J. Harris, Representation theory, a first course. Graduate texts in Math. 129, Springer-Verlag (1991)

\bibitem{G} D. Gaitsgory, On a vanishing conjecture appearing in the geometric Langlands correspondence,  Ann. of Math. (2) 160 (2004), no. 2, 617--682

\bibitem{Ku} S.~Kudla, Notes on the local theta correspondence, Notes from 10 lectures given at the European School on Group Theory in September 1996, downloadable at http://www.math.umd.edu/$\sim$ssk/castle.pdf

\bibitem{Ly} S. Lysenko, Moduli of metaplectic bundles on curves and Theta-sheaves, Ann. Scient. \'Ec. Norm. Sup. 4 s\'erie, t.39 (2006), 415-466

\bibitem{Ly2} S. Lysenko, Whittaker and Bessel functors for $\GSp_4$, Ann. Institut Fourier, t.56, No. 5 (2006), 1505-1565
\bibitem{Ly3} S. Lysenko, Local geometrized Rankin-Selberg method for $\GL(n)$, Duke Math.~J., 111 (2002), no. 3, 451--493.

\bibitem{Ly4} S. Lysenko, Geometric theta-lifting for the dual pair $\SO_{2m}$, $\Sp_{2n}$, math.RT/0701170

\bibitem{MV} I. Mirkovi\'c, K. Vilonen, Geometric Langlands duality and representations of algebraic groups over commutative rings, math.RT/0401222, to appear in Ann. of Math.
\bibitem{MVW} C. Moeglin, M.-F. Vigneras, J.L. Waldspurger, Correspondence de Howe sur un corps $p$-adique, Lecture Notes in Math. 1291 (1987)  
\bibitem{M} D. Mumford, Prym varieties I, Contribution to analysis , pp. 325-350. Academic Press (1974)
\bibitem{Pr} D. Prasad, Weil Representation, Howe duality, and the Theta
correspondence (lectures given in Montreal)
http://www.mri.ernet.in/~mathweb/dprasad.html
\bibitem{R} S. Rallis, Langlands functoriality and the Weil representation, Amer. J. Math., vol. 104, No. 3 (1982), p. 469-515
\bibitem{Wa} J.-L. Waldspurger, Sur les valeurs de certaines fonctions L automorphes en leur centre de symetrie, Comp. Math. 54 (1985), p. 173-242
\bibitem{We} A.~Weil, Sur la formule de Siegel dans la th\'eorie des groupes classiques, Acta Math. 113 (1965), p. 1-87, ou Oeuvres Sc., vol. III, p. 71-157.
\end{thebibliography}
\end{document}